\documentclass{amsart}
\usepackage{etex}
\usepackage{latexsym}
\usepackage{amsfonts,amsmath,amssymb,amscd,amsxtra,amsthm,verbatim,calc}
\usepackage{MnSymbol}
\usepackage[all]{xy}
\usepackage[a4paper,top=2cm,bottom=2cm,left=3.4cm,right=3.4cm]{geometry}

\usepackage{eucal}

\usepackage{tikz}
\usetikzlibrary{math} 
\usetikzlibrary{cd}
\usetikzlibrary{patterns}
\usetikzlibrary{backgrounds}

\usepackage{amsmath}
\usepackage{amstext}
\usepackage{amssymb}
\usepackage{amsthm}
\usepackage{amscd}

\usepackage{lscape}
\usepackage{longtable}

\usepackage{epsfig}
\input txdtools
\let\et=\etexdraw
\def\etexdraw{\drawbb\et}

\usepackage{arydshln}

\theoremstyle{plain}
\newtheorem{thm}{Theorem}[section]
\newtheorem{thm*}{Theorem}
\newtheorem{lem}[thm]{Lemma}
\newtheorem{prop}[thm]{Proposition}
\newtheorem{prop*}[thm*]{Proposition}
\newtheorem{cor}[thm]{Corollary}

\theoremstyle{definition}
\newtheorem{defn}[thm]{Definition}

\newtheorem{ex}[thm]{Example}
\newtheorem{qu}[thm]{Question}
\newtheorem{notation}[thm]{Notation}

\theoremstyle{remark}

\newtheorem{rem}[thm]{Remark}

\DeclareMathOperator{\Supp}{Supp}

\DeclareMathOperator{\Hred}{\widetilde{H}}
\DeclareMathOperator{\link}{link}
\DeclareMathOperator{\Sym}{Sym}
\DeclareMathOperator{\Skel}{Skel}

\newcommand{\ZZ}{\mathbb{Z}}

\newcommand{\RR}{\mathbb{R}}

\newcommand{\KK}{\mathbb{K}}

\newcommand{\PP}{\mathbb{P}}

\newcommand{\hh}{\widehat{h}}

\DeclareMathOperator{\calI}{\mathfrak{I}}
\DeclareMathOperator{\calP}{\mathcal{P}}
\DeclareMathOperator{\barP}{\widehat{\mathcal{P}}}
\DeclareMathOperator{\calB}{\mathcal{B}}

\newcommand{\bc}{\mathbf{c}}
\newcommand{\bd}{\mathbf{d}}

\newcommand{\bm}{\mathbf{m}}
\newcommand{\bp}{\mathbf{p}}
\newcommand{\bs}{\mathbf{s}}

\newcommand{\tDelta}{\widetilde{\Delta}}

\newcommand{\tg}{\widetilde{g}}

\newcommand{\BDelta}{\mathcal{B}(\Delta)}
\newcommand{\barDelta}{\widehat{\Delta}}

\newcommand{\eseq}[2]{\mathbf{e}^{#1}_{#2}}
\newcommand{\nth}[2][th]{{#2}^\text{#1}}

\usepackage{xparse}
\NewDocumentCommand{\lkds}{O{\sigma} O{\Delta}}{\link_{#2} #1}
\NewDocumentCommand{\Idstar}{O{\Delta}}{I_{{#1}^*}}

\begin{document}
	
	\title{Stanley-Reisner ideals with pure resolutions}
		\author{David Carey}
		\email{davidcarey277@gmail.com}
		\address{School of Mathematics,
			University of Sheffield, Hicks Building, Sheffield S3 7RH, United Kingdom}
		\author{Mordechai Katzman}
		\email{M.Katzman@sheffield.ac.uk}
		\address{School of Mathematics,
			University of Sheffield, Hicks Building, Sheffield S3 7RH, United Kingdom}
			
	

	
	\begin{abstract}
		This paper investgates Stanley-Reisner ideals with pure resolutions. We first describe two infinite families of such ideals associated to highly symmetric complexes. We then prove a partial analogue to the first Boij-Söderberg Conjecture for Stanley-Reisner ideals, by detailing an algorithm for constructing Stanley-Reisner ideals with pure Betti diagrams of any given shape, save for an initial shift.
	\end{abstract}
	
	\maketitle
	\section{Introduction}\label{Section: Introduction}
	
	This paper studies Stanley-Reisner ideals whose minimal free resolutions are pure.
	
	Recall that for any standard $\mathbb{N}$-graded ring $R$ and finitely generated graded $R$-module $M$, we can construct a
	graded minimal free resolution
	$$ \dots \rightarrow F_i \rightarrow \dots \rightarrow F_{1} \rightarrow F_0 \rightarrow M \rightarrow 0 \label{eq1}$$
	with $F_i=\oplus_{j\geq i} R(-j)^{\beta_{i j}}$
	where $R(-d)$ denotes the $R$-module $R$ with its degree shifted so that $R(-d)_\delta=R_{\delta +d}$ and
	$\beta_{i j}$ is a positive integer for all $i\geq 0$ and finitely many $j\geq i$.
	While graded modules may afford different minimal free resolutions, the \emph{graded Betti numbers} $\beta_{i j}$ are invariants of $M$.
	
	The \emph{Betti diagram} of $M$, denoted $\beta(M)$, is a matrix containing the graded Betti numbers of $M$. In order to reduce the number of rows of this matrix, the standard notational convention is to write $\beta(M)$ as a matrix $(a_{ij})$ with $a_{ij}=\beta_{j,i+j}$, as shown.
	$$\begin{array}{c | cccc}
		0&\beta_{0,0}&\beta_{1,1}&\dots&\beta_{n,n}\\
		1&\beta_{0,1}&\beta_{1,2}&\dots&\beta_{n,n+1}\\
		\vdots&\vdots&\vdots&\vdots&\vdots
	\end{array}$$
	So for example
	$$\begin{array}{l | rrr }
		2 & 5 & 5 & .\\ 
		3 & . & . & 1\\ 
	\end{array}$$
	denotes a pure resolution 
	$$0\rightarrow R(-5) \rightarrow R(-3)^5 \rightarrow R(-2)^5 \rightarrow 0.$$
	
	Throughout the rest of this paper $R$ will denote a polynomial ring $\mathbb{K}[x_1, \dots, x_n]$ where $\KK$ is an unspecified field of arbitrary characteristic. Hilbert's Syzygy Theorem now implies that $\beta_{i j}=0$ for all $i>n$. We use $\Hred_i(\Delta)$ to denote $\Hred_i(\Delta;\KK)$.
	
	Minimal free resolutions of $R$-modules were introduced by David Hilbert in 1890 as a tool for the study of invariants, and gave an early impetus for the
	development of homological methods in commuative algebra, now so ubiquitous in the field.
	
	Despite more than a century of study of free resolutions, many questions about these remain open, and perhaps the most basic of all is: which values of
	$\{ \beta_{i j} \}$ can occur as Betti numbers of a finitely generated graded $R$-module? Rather than looking at individual Betti-diagrams $\{\beta_{i j}\}$ one can consider
	the rational cones generated by all Betti diagrams of modules from a certain category and one can aim to describe these cones, e.g., it terms of their
	extremal rays. This was the approach taken in  \cite{BoijSoderberg08} where the authors formulated the famous
	Boij-S\"{o}derberg conjectures (now theorems), one of them being that
	the Betti diagram of Cohen-Macaulay $R$-module is a
	non-negative rational linear combination
	of \emph{pure} Betti diagrams.
	A pure Betti diagram is the Betti diagram $\{\beta_{i j}\}$ where for each $i\geq 0$, $\beta_{i j}\neq 0$ for at most one $j$, in other words the Betti diagram arising from a graded minimal free resolution of the form
	$$ \dots \rightarrow  R(-d_i)^{\beta_{i d_i}} \rightarrow \dots \rightarrow  R(-d_1)^{\beta_{1 d_1}} \rightarrow R(-d_0)^{\beta_{0 d_0}} \rightarrow 0 \label{eq2} .$$
	The precise formulation of the conjecture in terms of rational cones as follows.
	Fix $0\leq c \leq n$;  for $\mathbf{a}=(a_0, \dots, a_c)$ and
	$\mathbf{b}=(b_0, \dots, b_c)$ in $\mathbb{N}^{c+1}$
	let $\mathbb{D}(\mathbf{a}, \mathbf{b})$ denote the set of Betti diagrams $\{ \beta_{ i j} \}$
	where $\beta_{ i j} \neq 0$ only when $0\leq i\leq c$ and  $a_i \leq j \leq b_i$.
	Consider the rational cone $C$ generated by Betti diagrams  $\beta_{i j}\in \mathbb{D}(\mathbf{a}, \mathbf{b})$ of graded Cohen-Macaulay $R$-modules of codimension $c$.
	The content of \cite[Conjecture 2.4]{BoijSoderberg08} is that the extremal rays of $C$ are given by pure Betti diagrams.
	For an account of these conjectures and their eventual proofs see \cite{Floystad12}.
	
	In view of the great power of Boij-S\"{o}derberg's theory to describe possible Betti diagrams, one may wonder whether a similar theory exists if, instead of looking at Betti cones of
	graded Cohen-Macaulay $R$ modules of given codimension $c$, one were to look at Betti cones of other fixed class of graded $R$-modules.
	The motivation for this paper is our desire to discover a Boij-S\"{o}derberg's-like theory  where that class of graded $R$-modules consists of direct sums of Stanley-Reisner ideals (i.e., ideals generated by square-free monomials; see section \ref{Section: elementary concepts and tools}.)
	The guiding question behind our research is  the following:
	
	\begin{qu}
		What are all the possible Betti diagrams of Stanley-Reisner ideals?
		Consider the rational cone $\mathcal{C}$ generated by Betti diagrams  
		$\beta_{i j}\in \mathbb{D}(\mathbf{a}, \mathbf{b})$ of Stanley-Reisner ideals;
		what are the extremal rays of $\mathcal{C}$?
	\end{qu}
	
	Note that we do not restrict our attention to \emph{Cohen-Macaulay} Stanley-Reisner ideals, and as a result pure Betti diagrams are not necessarily 
	extremal rays in $\mathcal{C}$ (but the extremal rays of the rational cone generated by all pure Betti diagrams of the same same shape are extremal rays of $\mathcal{C}$). In addition, some extremal rays in $\mathcal{C}$ do not correspond to pure Betti diagrams. See section \ref{Section: Pure resolutions of Stanley-Reisner ideals in small number of variables} for concrete examples with small vertex sets.
	
	This stands in contrast to the original Boij--S\"{o}derberg setting: there, every Betti diagram of a Cohen--Macaulay module is a positive combination of pure ones, so pure diagrams are exactly the extremal rays of the relevant cone. Boij and S\"{o}derberg later showed this phenomenon persists without the Cohen--Macaulay hypothesis \cite{BoijSoderberg12}: every Betti diagram of a graded module decomposes as a positive rational combination of pure diagrams, so pure diagrams remain the extremal building blocks for Betti diagrams of graded modules in general. Since a Stanley-Reisner ideal is in particular a graded module, its pure Betti diagrams, when they exist, are still extremal rays of this larger ambient cone, even on the occasions where they fail to be extremal once we restrict attention to the smaller cone $\mathcal{C}$ generated by Stanley-Reisner ideals alone. Constructing pure Betti diagrams of every degree type for Stanley-Reisner ideals tests how this extremal behaviour changes in the more rigid monomial setting.
	
	There is probably no hope to answer this general question, but in view of the central role of pure Betti diagrams as extremal rays, we ask a slightly less ambitious question:
	
	\begin{qu}
		What are all the possible \emph{pure} Betti diagrams of Stanley-Reisner ideals?
	\end{qu}
	

	Special cases of this question have been raised in other contexts. For example, consider a special type of pure resolutions,
	namely \emph{$t$-linear resolutions}. These are pure resolutions with Betti diagram $\beta_{0, t}\neq 0, \beta_{1, t+1}\neq 0, \dots \beta_{p,t+p}\neq 0$.
	Stanley-Reisner ideals with $2$-linear resolutions were first classified by Ralf Fr\"{o}berg in \cite{Froberg88} and later
	those with $t$-linear resolutions were described by John Eagon and Victor Reiner in \cite{EagonReiner98}.
	At about the same time Winfried Bruns and Takayuki Hibi classified Cohen–Macaulay partially ordered sets whose corresponding
	Stanley-Reisner ideal has a pure minimal resolution (\cite{BrunsHibi95}) and later Stanley-Reisner ideals of dimension 1 and 2 with pure resolutions (\cite{BrunsHibi98}).
	
	As observed in \cite{BrunsHibi98} the complete classification of all pure Betti diagrams of Stanley-Reisner ideals
	might be too ambitious a goal, but some understanding of the constraints on these might be feasible. This paper makes a contribution to this project by
	describing two infinite families of Stanley-Reisner ideals with pure resolutions of arbitrary length, and also detailing an algorithm for generating Stanley-Reisner ideals with pure resolutions of any given degree type (as defined in the next section). 
	
	The tools used in this paper are combinatorial and are mostly applications of Hochster's formula for the Betti numbers of Stanley-Reisner ideals.
	
	\subsection{A Partial Analogue to the First Boij-S\"{o}derberg Conjecture for Stanley-Reisner Ideals}\label{Section: Intro Part 2}
	
	The main result of this paper can be seen as a partial analogue to the First Boij-S\"{o}derberg Conjecture, for Stanley-Reisner ideals. In order to state it, we need to introduce the following concepts.
	
	\begin{defn}
		As in \cite{BrunsHibi95}, given a pure graded free resolution $F_\bullet$
		\begin{equation*}\label{pure-res}
			0\rightarrow R(-c_p)^{\beta_{p,c_p}} \rightarrow \dots \rightarrow R(-c_1)^{\beta_{1,c_1}} \rightarrow R(-c_0)^{\beta_{0,c_0}}
		\end{equation*}
		We define the \emph{shift type} of $F_\bullet$ as $(c_p,\dots, c_0)$, and its \emph{degree type} to be the positive sequence $(c_p-c_{p-1},\dots,c_1-c_0)$. We also extend this terminology accordingly to define the shift types and degree types of pure Betti diagrams.
	\end{defn}
	\begin{rem}
		Our notion of degree type differs slightly from that found in \cite{BrunsHibi95}, which defines it as the sequence $(c_p-c_{p-1},...,c_1-c_0,c_0)$. This is the sequence of degrees of the maps in the resolution. Our version of degree type records the degrees of all of these maps except for $R(-c_0)^{\beta_{0,c_0}}\rightarrow I$. Note in particular that under our definition, The degree type of a pure resolution contains slightly less information than its shift type, but the shift type is uniquely determined by the degree type and the value of $c_0$.
	\end{rem}
	
	The first Boij-S\"{o}derberg Conjecture (\cite{Floystad12}, Theorem 1.9) states that for any strictly decreasing sequence of integers $\bc = (c_p,\dots,c_0)$ with $p\leq n$, there exists a Cohen-Macaulay $R$-module of codimension $p$ with a pure Betti diagram of shift type $\bc$. One might wonder if a similar result holds for Stanley-Reisner ideals: that is, for a strictly decreasing sequence of integers $\bc = (c_p,\dots,c_0)$, does there exist a Staney-Reisner ideal with a pure Betti diagram of shift type $\bc$? The answer to this is \textit{no}. This can be seen immediately, because no ideal in $R$ can have a pure Betti diagram with a shift type that contains negative numbers. However, even if we restricted our attention purely to sequences of \textit{positive} integers, the answer would still be no: for example, in \cite{Carey24} Chapter 4 page 69, we show that no Stanley-Reisner ideal can have a pure diagram of shift type $(5,2)$.
	
	We can, however, broaden the question slightly: instead of asking whether there are pure Betti diagrams of Stanley-Reisner ideals corresponding to every possible \textit{shift type}, we ask whether there are diagrams corresponding to every possible \textit{degree type}. The main result of this paper is a positive answer to this question.
	\begin{thm}\label{Theorem: PR Complexes of Any Degree Type}
		Let $\bd=(d_p,\dots,d_1)$ be any sequence of positive integers. There exists a Stanley-Reisner ideal $I$ such that the Betti diagram $\beta(I)$ is pure with degree type $\bd$.
	\end{thm}
	
	The structure of the paper is as follows. Section \ref{Section: elementary concepts and tools} gives a brief exposition of some relevant background results in the literature, and then introduces the theory of \textit{PR Complexes}, simplicial complexes whose dual Stanley-Reisner ideals have pure resolutions. Sections 3 and 4 provide constructions for two families of PR complexes with particular degree types. Section 5 is devoted to the proof of Theorem \ref{Theorem: PR Complexes of Any Degree Type}. And finally Section 6 contains a list of extremal rays for the Betti cones on Stanley-Reisner ideals in small numbers of variables. It is worth highlighting that the results of Sections 3, 4 and 5 are independent of each other. In particular, the proof of Theorem \ref{Theorem: PR Complexes of Any Degree Type} in Section 5 does not depend on results in Sections 3 or 4.
	
	\section{Elementary concepts and tools}\label{Section: elementary concepts and tools}
	In this section we note some elementary tools used throughout this paper. We begin with a (very) brief overview of simplicial complexes and Stanley-Reisner Theory.
	
	A \textit{simplicial complex} on vertex set $V$ is a set of subsets of $V$ such that for any $G \subseteq F \subseteq V$, if $F$ is in $\Delta$ then $G$ is also in $\Delta$. We refer to the elements of $\Delta$ as \textit{faces} and the maximal elements as \textit{facets}. We define the dimension of a face $\sigma$ of $\Delta$ to be $|\sigma|-1$, and the dimension $\dim \Delta$ of $\Delta$ to be the dimension of its largest facet. If all facets of $\Delta$ have the same dimension we say $\Delta$ is \textit{pure}.
	
	We often write simplicial complexes in \textit{facet notation}, as follows. 
	\begin{notation}\label{Notation: Facet Notation for Complexes}
		For subsets $F_1$,...,$F_m$ in a vertex set $V$, we use $\langle F_1,...,F_m\rangle$ to denote the smallest simplicial complex on $V$ containing $F_1$,...,$F_m$. Note that if the set $\{F_1,...,F_m\}$ is irredundant, then this is the same as the complex on $V$ with facets $F_1,...,F_m$.
	\end{notation}
	
	The following constructions will be of particular importance in this paper.
	\begin{defn}\label{Definition: Simplicial Complex Constructions}
		Let $\Delta$ be a simplicial complexes on vertex set $V$.
		\begin{itemize}
			\item [(a)] The \textit{Alexander Dual} of $\Delta$ is the complex $$\Delta^* = \{F \subseteq V : V - F \notin \Delta\}.$$
			
			\item [(b)] The \textit{link} of a face $f\in \Delta$ is the complex
			$$\link_\Delta f = \{ \sigma\in \Delta \,|\, f\cap \sigma = \emptyset \text{ and } f\cup \sigma\in \Delta \}.$$
			
			\item [(c)] For any $r\geq -1$, the \textit{$r$-skeleton} of $\Delta$ is the complex $$\Skel_r(\Delta) = \{\sigma \in \Delta : \dim \sigma \leq r\}.$$
			If $\Delta$ is a full simplex on vertex set $V$ we often denote this as $\Skel_r(V)$.
			
			\item [(d)] For a complex $\Gamma$ on vertex set $U$ with $U\cap V = \emptyset$, the \textit{join} of $\Delta$ and $\Gamma$ is the complex $$\Delta \ast \Gamma = \{\sigma_1 \sqcup \sigma_2 : \sigma_1 \in \Delta, \sigma_2 \in \Gamma\}.$$
			
			\item [(e)] The \textit{induced subcomplex} of $\Delta$ on a subset $U\subseteq V$ is the complex $$\Delta_U=\{\sigma \in \Delta : \sigma \subseteq U\}.$$
			
			\item [(g)] The \textit{deletion} of a face $f$ from $\Delta$ is the complex $$\Delta-f = \{\sigma \in \Delta : \sigma \nsupseteq f \}.$$
		\end{itemize}
	\end{defn}
	\begin{rem}
		Note that, except in the case where $f$ is a vertex, $\Delta-f$ is larger than the induced subcomplex $\Delta_{V-f}$, because it contains all faces of $\Delta$ which intersect with $f$ strictly along its boundary. The following example illustrates this.
		\begin{center}
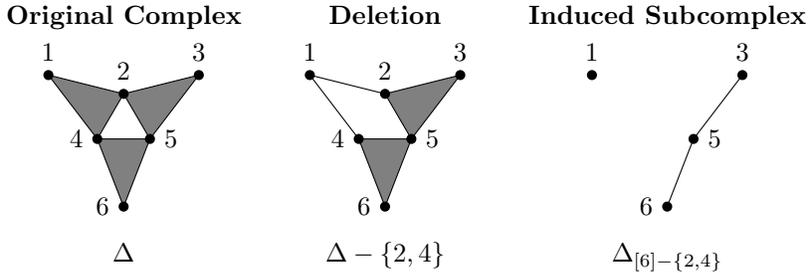
			\begin{tabular}{ c c c c c }
				\textbf{Original Complex}&&\textbf{Deletion}&&\textbf{Induced Subcomplex} \\
				\begin{tikzpicture}[scale = 0.5]
					\tikzstyle{point}=[circle,thick,draw=black,fill=black,inner sep=0pt,minimum width=3pt,minimum height=3pt]
					\node (a)[point,label=above:$1$] at (0,3.5) {};
					\node (b)[point,label=above:$2$] at (2,3) {};
					\node (c)[point,label=above:$3$] at (4,3.5) {};
					\node (d)[point,label=left:$4$] at (1.3,1.8) {};
					\node (e)[point,label=right:$5$] at (2.7,1.8) {};
					\node (f)[point,label=left:$6$] at (2,0) {};	
					
					\begin{scope}[on background layer]
						\draw[fill=gray] (a.center) -- (b.center) -- (d.center) -- cycle;
						\draw[fill=gray] (b.center) -- (c.center) -- (e.center) -- cycle;
						\draw[fill=gray]   (d.center) -- (e.center) -- (f.center) -- cycle;
					\end{scope}
				\end{tikzpicture}&&
				\begin{tikzpicture}[scale = 0.5]
					\tikzstyle{point}=[circle,thick,draw=black,fill=black,inner sep=0pt,minimum width=3pt,minimum height=3pt]
					\node (a)[point,label=above:$1$] at (0,3.5) {};
					\node (b)[point,label=above:$2$] at (2,3) {};
					\node (c)[point,label=above:$3$] at (4,3.5) {};
					\node (d)[point,label=left:$4$] at (1.3,1.8) {};
					\node (e)[point,label=right:$5$] at (2.7,1.8) {};
					\node (f)[point,label=left:$6$] at (2,0) {};	
					
					\begin{scope}[on background layer]
						\draw (d.center) -- (a.center) -- (b.center);
						\draw[fill=gray] (b.center) -- (c.center) -- (e.center) -- cycle;
						\draw[fill=gray]   (d.center) -- (e.center) -- (f.center) -- cycle;
					\end{scope}
				\end{tikzpicture} && 	\begin{tikzpicture}[scale = 0.5]
					\tikzstyle{point}=[circle,thick,draw=black,fill=black,inner sep=0pt,minimum width=3pt,minimum height=3pt]
					\node (a)[point,label=above:$1$] at (0,3.5) {};
					\node (c)[point,label=above:$3$] at (4,3.5) {};
					\node (e)[point,label=right:$5$] at (2.7,1.8) {};
					\node (f)[point,label=left:$6$] at (2,0) {};	
					
					\begin{scope}[on background layer]
						\draw(c.center) -- (e.center);
						\draw (e.center) -- (f.center);
					\end{scope}
				\end{tikzpicture} \\
				$\Delta$&& $\Delta-\{2,4\}$ && $\Delta_{[6]-\{2,4\}}$
			\end{tabular}
		\end{center}
	\end{rem}
	
	\begin{defn}\label{Definition: SR Ideals}
		Let $\Delta$ be a simplicial complex on vertex set $[n]=\{1,\dots,n\}$.
		The \emph{Stanley-Reisner} ideal $I_\Delta$ is the ideal of $R$ generated by the square-free monomials
		$$\{ x_{i_1} \dots x_{i_k} \,|\, \{ i_1, \dots,  i_k \}\notin \Delta \} .$$
	\end{defn}
	Every ideal generated by square free monomials is a Stanley-Reisner ideal (and vice-versa).
	
	\begin{rem}
		Most of the simplicial complexes presented in this paper are not on vertex set $[n]$. We define their Stanley-Reisner ideals by first relabelling the vertices. Note that any two relabellings give isomorphic Stanley-Reisner ideals.
	\end{rem}
	
	The main tool for investigating Betti numbers of a Stanley-Reisner ideal $I_\Delta$ is Hochster's Formula, which expresses the Betti numbers as a sum of dimensions
	of the reduced homologies of sub-simplicial complexes of $\Delta$.
	
	\begin{thm}[Hochster's Formula (Theorem 5.1 in \cite{Hochster77})] \label{Theorem: Hochster's formula}
		Let $\Delta$ be a simplicial complex on a vertex set $V$ of size $n$.
		\[\beta_{i d} (I_\Delta) = \sum_{W \subseteq V, \#W=d} \dim_\mathbb{K} \Hred_{d-i}(\Delta_W; \mathbb{K})\]
		where for any
		$W \subseteq V$, $\Delta_W$ denotes the simplicial complex
		with vertex set $W$ and whose faces are the faces of $\Delta$ containing only vertices in $W$.
	\end{thm}
	
	Hochster's formula can be reformulated in terms of the Alexander dual of $\Delta$ as follows.

\begin{thm}[Alexander Dual version of Hochster's Formula (ADHF), see Corollary 1.40 in \cite{MillerSturmfels05}] \label{Theorem: Hochster's formula  -- dual version}
 		Let $\Delta$ be a simplicial complex on a vertex set $V$ of size $n$.
		$$\beta_{i d}(\Idstar)=\sum_{f \in \Delta, |f|=n-d} \dim_{\mathbb{K}} \Hred_{i-1}({\link_{\Delta}} f; \mathbb{K})$$
	\end{thm}

	ADHF plays a central role in many of the proofs in this paper (and hence the acronym.) Thus many of our results revolve around the computation of simplicial homology, and in particular, make use of deformation retractions. The following lemma will be of particular importance in this regard. It provides sufficient conditions for ensuring that a complex $\Delta$ deformation retracts on to the complex $\Delta - g$ (the \textit{deletion of $g$ from $\Delta$}, as given in Definition \ref{Definition: Simplicial Complex Constructions}). We suspect this lemma exists in the literature, but have been unable to find it.
	\begin{lem}\label{Lemma: Deformation Retract}
		Let $\Delta$ be a simplicial complex, and let $g\subsetneqq f$ be faces of $\Delta$ such that every facet of $\Delta$ that contains $g$ also contains $f$. There is a deformation retraction $\Delta \rightsquigarrow \Delta - g$, obtained by 
		collapsing $g$ onto 
		a vertex in $f-g$.
	\end{lem}
	\begin{proof}
		Suppose $\Delta$ has $n$ vertices $v_1,\dots,v_n$ with $g=\{v_1,\dots,v_k\}$ and $v_{k+1}$ in $f-g$.
		
		Let $\{e_1,\dots,e_n\}$ be the canonical basis for $\RR^n$, and let $X=X_{\Delta}$ be the geometric realization of $\Delta$ in $\RR^n$, with $v_i$ realized as $e_i$. Similarly, let $X_{\Delta-g}$ be the geometric realization of $\Delta-g$. For each nonempty face $\sigma =\{v_{i_1},\dots,v_{i_r}\}$ of $\Delta$, we define $X_\sigma$ to be the set $\left\{\lambda_1 e_{i_1}+\dots \lambda_r e_{i_r} : \lambda_1, \dots, \lambda_r > 0, \sum_{j=1}^r \lambda_j = 1\right\} \subset X$. Note that $X=\bigcup_{\sigma \in \Delta} X_\sigma$, while $X_{\Delta-g}=\bigcup_{\sigma \in \Delta, \sigma \nsupseteq g} X_\sigma$.
		
		Let $\bp$ be a point in $X$. We may write it as $\bp = \sum_{j=1}^n \lambda_j e_j$ where the coefficients $\lambda_j$ are all nonnegative and sum to $1$. In particular, $\bp$ lies inside $X_\sigma$ for some $\sigma\in \Delta$ containing $g$ if and only if all of $\lambda_1,\dots,\lambda_k$ are positive.
		
		We define $\lambda = \min \{\lambda_1, \dots, \lambda_k\}$. This allows us to rewrite the point $\bp$ as
		\begin{equation}\label{Equation: bp}
			\bp = \lambda(e_1 + \dots + e_k) + \sum_{i=1}^{k} (\lambda_i - \lambda) e_i + \lambda_{k+1} e_{k+1} + y
		\end{equation}
		for some $y$ in the span of $\{e_{k+2},\dots,e_n\}$. Note that the coefficients $(\lambda_i-\lambda)$ are all nonnegative, and (by the definition of $\lambda$) at least one of them is zero. Note also that $\lambda$ itself is nonzero if and only if all of $\lambda_1,\dots,\lambda_k$ are positive, which occurs if and only if $\bp$ is in $X_\sigma$ for some $\sigma$ containing $g$. In other words, we have $\lambda = 0$ if and only if $\bp$ lies inside $X_{\Delta-g}$.
		
		Using this notation for the points in $X$, we can define a function $\varphi: X \times [0,1] \rightarrow \RR^n$ as follows. For $\bp$ as in Equation (\ref{Equation: bp}) and $0\leq t \leq 1$ we define 
		\begin{equation}\label{Equation: varphi}
			\varphi(\bp, t) = \lambda (1-t)(e_1 + \dots + e_k) + \sum_{i=1}^{k} (\lambda_i - \lambda) e_i + (\lambda_{k+1} + k \lambda t) e_{k+1} + y.
		\end{equation}
		
		We claim that $\varphi$ is a deformation retraction from $X$ to $X_{\Delta-g}$.
		
		First, note that $\varphi$ is continuous. Indeed, the function $\lambda=\min \{\lambda_1,\dots,\lambda_k\}$ is continuous in the variables $\lambda_i$, and hence so is each summand in Equation (\ref{Equation: varphi}). Each summand is also continuous in $t$.
		
		Next, we note that $\varphi(*,0)$ is the identity on $X$. Moreover, for $\bp$ in $X_{\Delta-g}$ we have $\lambda=0$, and hence $\varphi(\bp,t)=\bp$ for every $0\leq t \leq 1$.
		
		It remains to show that the image of $\varphi(*,t)$ lies inside $X$ for every value of $t$, and in particular that the image of $\varphi(*,1)$ is $X_{\Delta-g}$.
		
		For the former claim, we note that the sum of the coefficients of $e_1,\dots,e_n$ in the decomposition of $\varphi(\bp,t)$ is the same as the sum of these coefficients in the decomposition of $\bp$ (which is $1$), and all of these coefficients are nonnegative. We may assume that $\bp$ lies inside $X_{\sigma}$ for some face $\sigma$ in $\Delta$ containing $g$ (otherwise $\bp$ lies inside $X_{\Delta-g}$ and we are already done). By our assumption on $g$ and $f$ we have that $\sigma\cup \{v_{k+1}\}$ is also a face of $\Delta$. If $\sigma = \{v_1,\dots,v_k,v_{i_1},\dots,v_{i_r}\}$ for some $k<i_1<\dots<i_r\leq n$, then the vectors in the decomposition of $\varphi(\bp,t)$ with strictly positive coefficients are all contained in $\{e_1,\dots,e_k,e_{i_1},\dots,e_{i_r}\}\cup\{e_{k+1}\}$. We conclude that $\varphi(\bp,t)$ lies inside $X_{\tau}$ for some $\tau \subseteq \sigma \cup \{v_{k+1}\}$, and hence inside $X$.
		
		For the latter claim, we note that the coefficients of the vectors $e_1,\dots, e_k$ in $\varphi(\bp,1)$ are $\lambda_1-\lambda,\dots,\lambda_k-\lambda$, and at least one of these must be zero by the definition of $\lambda$. Thus $\varphi(\bp,1)$ lies inside $X_{\Delta-g}$.
	\end{proof}
	
	In particular we have the following corollary, which will be sufficient for our needs in most (but not all) cases.
	\begin{cor}\label{Corollary: Deformation Retract}
		Let $\Delta$ be a simplicial complex, and let $a$ and $b$ be distinct vertices of $\Delta$ satisfying the following conditions.
		\begin{enumerate}
			\item There is at least one facet of $\Delta$ containing $a$.
			\item Every facet of $\Delta$ containing $a$ also contains $b$.
		\end{enumerate}
		There is a deformation retraction $\Delta \leadsto \Delta - \{a\}$ given by the map $a\mapsto b$.
	\end{cor}
	\begin{proof}
		This comes from Lemma \ref{Lemma: Deformation Retract}, setting $g=\{a\}$ and $f=\{a,b\}$. Note that $g$ is a face of $\Delta$ by assumption (1), and $f$ is a face of $\Delta$ by assumption (2).
	\end{proof}
	
	We will also need the K\"{u}nneth Formula for computing the homologies of joins of simplicial complexes, which is the following.
	\begin{prop}[K\"{u}nneth Formula for Joins of Simplicial Complexes, see \cite{Hatcher02} page 276, Corollary 3B.7]\label{Proposition: Kunneth Formula}
		Let $\Delta$ be a simplicial complex and suppose we have $\Delta = A \ast B$ for two subcomplexes $A$ and $B$. For any integer $r\geq -1$ we have an isomorphism
		\begin{equation*}
			\Hred_{r+1}(\Delta)\cong \sum_{i+j=r}\Hred_i(A)\otimes\Hred_j(B).
		\end{equation*}
	\end{prop}
	
	\subsection{PR Complexes} Our primary aim is to construct Stanley-Reisner ideals with pure resolutions of varying degree types. In particular we wish to construct simplicial complexes $\Delta$ for which the diagram $\beta(I_{\Delta^*})$ is pure.
	
	The Alexander Dual version of Hochster's Formula (Theorem \ref{Theorem: Hochster's formula  -- dual version}) reduces this to an entirely combinatorial problem, as shown below.
	
	\begin{cor}\label{Corollary: PR Complexes}
		Let $\Delta$ be a simplicial complex with vertices $x_1,\dots,x_n$. The diagram $\beta=\beta(I_{\Delta^*})$ is pure if and only if $\Delta$ satisfies the following condition:
		
		For every $i$, and every simplex $\sigma,\tau \in \Delta$, if $\Hred_i(\link_\Delta \sigma)\neq 0 \neq \Hred_i(\link_\Delta \tau)$ then $|\sigma|=|\tau|$.
	\end{cor}
	\begin{proof}
		The diagram $\beta$ is pure if and only if for every $i$, there exists at most one $c_i$ such that $\beta_{i,c_i}\neq 0$. By ADHF, this holds if and only if there are no two simplices of different sizes in $\Delta$ whose links both have non-trivial homology at the same degree.
	\end{proof}
	
	\begin{defn}\label{Definition: PR Complexes}
		We refer to complexes which satisfy the condition in Corollary \ref{Corollary: PR Complexes} as \textit{PR Complexes} (over $\KK$), where \textit{PR} stands for \textit{Pure Resolution}.
	\end{defn}
	\begin{rem}
		All of our work in this paper is done over the arbitrary field $\KK$. Hence, for the rest of this paper we will simply use the phrase `\textit{PR}' to mean `\textit{PR over} $\KK$'.
	\end{rem}
	
	\begin{lem}\label{Lemma: PR complexes are pure}
		All PR complexes are pure (i.e. their facets all have the same dimension).
	\end{lem}
	\begin{proof}
		The link of a facet is the irrelevant complex $\{\emptyset\}$ which has nontrivial $\nth[st]{(-1)}$ homology. Thus a PR complex cannot have two facets of different sizes.
	\end{proof}
	
	The following lemma will be particularly useful for us in our study of PR complexes.
	\begin{lem}\label{Lemma: Homology of descending links}
		Let $\Delta$ be any simplicial complex. Suppose we have some face $\sigma$ in $\Delta$ such that $\Hred_j(\lkds)\neq 0$ for some $j\geq 0$. There exists a chain of simplices $\sigma = \tau_j \subsetneqq \tau_{j-1} \subsetneqq \dots \subsetneqq \tau_0 \subsetneqq \tau_{-1}$ in $\Delta$ such that for each $-1\leq i\leq j$, we have $\Hred_i(\lkds[\tau]_i)\neq 0$.
	\end{lem}
	\begin{proof}
		We prove this algebraically, setting $\beta=\beta(\Idstar)$. It suffices to show that there is some face $\sigma \subsetneqq \tau \in \Delta$ for which $\Hred_{j-1}(\lkds[\tau])\neq 0$. The result then follows by induction on $j$.
		
		By replacing $\Delta$ with $\lkds$, we may assume that $\sigma = \emptyset$, and thus we need only find a nonempty face $\tau$ in $\Delta$ whose link has $\nth[st]{(j-1)}$ homology. To find a candidate for $\tau$, note that by ADHF we have that $\beta_{j+1,n}\neq 0$. Hence there must be some $d<n$ such that $\beta_{j,d}\neq 0$. This means (again, by ADHF) that there exists some nonempty face $\tau$ in $\Delta$ of size $n-d$ for which $\Hred_{j-1}(\lkds)\neq 0$, as required.
	\end{proof}
	
	\begin{cor}\label{Corollary: PR Complex links have single homology}
		Let $\Delta$ be a PR complex. Every link in $\Delta$ has at most one nontrivial homology.
	\end{cor}
	\begin{proof}
		Let $\sigma$ be a face of $\Delta$ and suppose for contradiction that $\Hred_i(\lkds)\neq 0\neq \Hred_j(\lkds)$ for some $i < j$. By Lemma \ref{Lemma: Homology of descending links}, there exists some $\tau \supsetneqq \sigma$ such that $\Hred_i(\link_\Delta \tau)\neq 0$. But this contradicts the fact that $\Delta$ is PR, because $|\tau|>|\sigma|$. 
	\end{proof}
	\begin{rem}
		In particular, the complex $\Delta$ is equal to $\link_\Delta \emptyset$, so it must have at most one nontrivial homology itself.
	\end{rem}
	
	\begin{cor}\label{Corollary: Homology index sets decreasing sequence}
		Let $\Delta$ be a PR complex, and suppose $\sigma_1$ and $\sigma_2$ are faces of $\Delta$ with $\Hred_{i_1}(\lkds_1)\neq 0 \neq \Hred_{i_2} (\lkds_2)$. If $|\sigma_1|<|\sigma_2|$ then $i_1>i_2$.
	\end{cor}
	\begin{proof}
		We cannot have $i_1=i_2$ as this directly contradicts the PR condition. If $i_1< i_2$, then by Lemma \ref{Lemma: Homology of descending links} we can find some face $\sigma_3$ of $\Delta$ strictly containing $\sigma_2$ such that $\Hred_{i_1}(\lkds_3)\neq 0$, which also contradicts the PR property because $|\sigma_3|>|\sigma_1|$.
	\end{proof}
	
	Using these results about the homology index sets of PR complexes, we are now able to give an entirely combinatorial description of the degree type of a PR complex $\Delta$, which agrees with the degree type of the Betti diagram $\beta(\Idstar)$.
	
	\begin{defn}\label{Definition: PR Complex Degree Types}
		Let $\Delta$ be a PR Complex, and let $p$ be the maximum index for which there exists some face $\sigma$ in $\Delta$ such that $\Hred_{p-1}(\lkds)\neq 0$. For each $0\leq i\leq p$, we define $s_i$ to be the size $|\sigma|$ of the faces $\sigma$ of $\Delta$ for which $\Hred_{i-1}(\lkds)\neq 0$, and for each $1\leq i\leq p$, we define $d_i=s_{i-1}-s_i$. We call the sequence $(d_p,\dots,d_1)$ the \textit{degree type} of $\Delta$.
	\end{defn}
	\begin{rem}
		The integers $s_{p-1},\dots,s_0$ are well-defined by Lemma \ref{Lemma: Homology of descending links}, and form a strictly decreasing sequence by Corollary \ref{Corollary: Homology index sets decreasing sequence}. Thus the degree type $(d_p,\dots,d_1)$ must consist of positive integers.
		
		To see why this notion of the degree type is the same as the degree type of the pure diagram $\beta(\Idstar)$, suppose that $$0\rightarrow R(-c_p)^{\beta_{p,c_p}} \rightarrow ... \rightarrow R(-c_1)^{\beta_{1,c_1}} \rightarrow R(-c_0)^{\beta_{0,c_0}}\rightarrow \Idstar$$
		is a minimal graded free resolution of $\Idstar$. By definition the degree type of this resolution is the sequence $(d_p,\dots,d_1)$ where for each $1\leq i \leq p$ we define $d_i=c_i-c_{i-1}$. By ADHF, for each $1\leq i \leq p$ we have $c_i=n-s_i$ and $c_{i-1}=n-s_{i-1}$, and hence we have $d_i=(n-s_i)-(n-s_{i-1})=s_{i-1}-s_i$.
	\end{rem}

	\begin{defn}\label{Definition: offset}
		Let $\Delta$ be a PR complex and let $s_0,\dots,s_p$ be as in Definition \ref{Definition: PR Complex Degree Types} above. We call the value $s_p$ (i.e. the minimum size of a face of $\Delta$ whose link has homology) the \textit{offset} of $\Delta$.
	\end{defn}
	
	Note that a PR complex $\Delta$ has offset 0 if and only if it has nontrivial homology itself (because $\lkds[\emptyset]$ is equal to $\Delta$). Taken together, the degree type and offset of a PR complex $\Delta$ determine its dimension, as shown below.
	
	\begin{prop}\label{Proposition: Dimension of PR Complexes}
		Let $\Delta$ be a PR complex with degree type $\bd$ and offset $s$. We have $\dim \Delta = s + \sum \bd - 1$.
	\end{prop}
	\begin{proof}
		Using the notation of Definition \ref{Definition: PR Complex Degree Types} we have $\sum \bd = \sum_{i=1}^p d_i = \sum_{i=1}^p (s_i - s_{i-1}) = s_0 - s_p$. The value $s_p$ is the offset of $\Delta$ (i.e. $s=s_p$), and the value $s_0$ is the size of the facets of $\Delta$ (all of which are the same by Lemma \ref{Lemma: PR complexes are pure}). Thus $\dim \Delta = s_0 - 1 = s_p + (s_0 - s_p) - 1 = s + \sum \bd - 1$.
	\end{proof}

	\subsection{PR Complexes and Cohen-Macaulay Complexes}
	The PR condition is in fact a generalisation of Reisner's criterion for Cohen-Macaulay complexes, which is the following.
	\begin{thm}[Reisner's Cohen-Macaulay Criterion, see \cite{MillerSturmfels05} Theorem 5.53]
		Let $\Delta$ be a simplicial complex. The Stanley-Reisner ring $\KK[\Delta]$ is Cohen-Macaulay if and only if for any $\sigma$ in $\Delta$ and any $i\neq \dim \Delta - |\sigma|$ we have $\Hred_i(\lkds) = 0$. In this case we call $\Delta$ a \textit{Cohen-Macaulay complex} (over $\KK$).
	\end{thm}
	
	In particular, all Cohen-Macaulay complexes are PR. From an algebraic perspective, this is a result of the following theorem.
	\begin{thm}[Eagon-Reiner Theorem, see \cite{EagonReiner98} Theorem 3]\label{Theorem: CM iff Linear}
		Let $\Delta$ be a simplicial complex. The following are equivalent.
		\begin{enumerate}
			\item $\KK[\Delta]$ is Cohen-Macaulay.
			\item $I_{\Delta^*}$ has a linear resolution.
		\end{enumerate}
	\end{thm}
	A \textit{linear resolution} is a pure resolution of degree type $(1,...,1)$. Hence we can rephrase the Eagon-Reiner Theorem as follows.
	\begin{cor}\label{Corollary: CM iff PR of deg type (1...1)}
		Let $\Delta$ be a simplicial complex. The following are equivalent.
		\begin{enumerate}
			\item $\Delta$ is Cohen-Macaulay.
			\item $\Delta$ is PR with degree type $(1,\dots,1)$.
		\end{enumerate}
	\end{cor}
	
	\subsection{Some Motivating Examples}
	We now present two examples of PR complexes. The symmetries and combinatorial properties of these complexes motivate the families of complexes constructed in the next three sections.
	
	For both examples we use the Alexander Dual version of Hochster's Formula (ADHF) to demonstrate why the corresponding Betti diagram is pure.
	
	\begin{ex}\label{Example: PR Simplex}
		Let $\Delta$ be the boundary of the $(n-1)$-simplex on vertex set $[n]$, and let $\beta=\beta(\Idstar)$. For any integer $1\leq d \leq n$ there are ${n \choose d}$ faces of $\Delta$ of size $n-d$, and for any such face $\sigma$, the complex $\lkds$ is the boundary of the $(n-1-|\sigma|)$-simplex on vertex set $[n]-\sigma$, which has homology only at degree $n-|\sigma|-2= d-2$ (and this homology has dimension $1$). Thus, $\beta$ is equal to
		$$\begin{array}{c | cccc }
			1 & {n \choose 1} & \dots & {n \choose n-1} & {n \choose n}\\ 
		\end{array}$$
		which means $\Delta$ is PR with degree type $(\underbrace{1,\dots,1}_{n-1})$ (i.e. it is Cohen-Maculay).
	\end{ex}
	Hence we can obtain PR complexes with degree type $(\underbrace{1,\dots,1}_{n})$ for any $n$.
	\begin{rem}
		The algebraic interpretation of the above example is as follows. If $\Delta$ is the the boundary of the $(n-1)$-simplex on vertex set $[n]$, then its dual Stanley-Reisner ideal $\Idstar$ is the maximal ideal $(x_1,\dots,x_n)$ in the polynomial ring $R=\KK[x_1,\dots,x_n]$. Thus the Koszul complex
		\[
		\begin{tikzcd}[row sep=1.5em,column sep=1.5em]
			0 \arrow{r} &R(-n) \arrow{r} &R(-(n-1))^{{n \choose n-1}} \arrow{r} &\dots \arrow{r} &R(-1)^{n\choose 1} \arrow{r} &\Idstar \arrow{r} &0.
		\end{tikzcd}
		\]
		is a free resolution of $\Idstar$.
	\end{rem}
	
	\begin{rem}\label{Remark: Simplicial Spheres are CM}
		In fact Example \ref{Example: PR Simplex} is an instance of a more general result, namely that every simplicial sphere is Cohen-Macaulay (see \cite{BrunsHerzog98} Section 5.4), and more specifically, any simplicial $n$-sphere has degree type $(\underbrace{1,\dots,1}_{n})$. 
	\end{rem}
	
	\begin{ex}\label{Example: PR Zelda Symbol}
		Let $\Delta$ be the complex
		
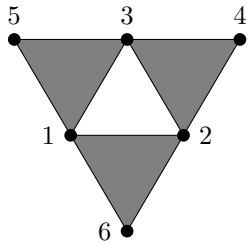
\begin{center}
			\begin{tikzpicture}[scale = 0.75]
				\tikzstyle{point}=[circle,thick,draw=black,fill=black,inner sep=0pt,minimum width=4pt,minimum height=4pt]
				\node (a)[point,label=above:$5$] at (0,3.4) {};
				\node (b)[point,label=above:$3$] at (2,3.4) {};
				\node (c)[point,label=above:$4$] at (4,3.4) {};
				\node (d)[point,label=left:$1$] at (1,1.7) {};
				\node (e)[point,label=right:$2$] at (3,1.7) {};
				\node (f)[point,label=left:$6$] at (2,0) {};	
				
				\begin{scope}[on background layer]
					\draw[fill=gray] (a.center) -- (b.center) -- (d.center) -- cycle;
					\draw[fill=gray] (b.center) -- (c.center) -- (e.center) -- cycle;
					\draw[fill=gray]   (d.center) -- (e.center) -- (f.center) -- cycle;
				\end{scope}
			\end{tikzpicture}
		\end{center}
		and let $\beta=\beta(\Idstar)$. This complex has only $1^\text{st}$ homology. It has three vertices with links of the form \begin{tikzpicture}[scale = 0.2]
			\tikzstyle{point}=[circle,thick,draw=black,fill=black,inner sep=0pt,minimum width=2pt,minimum height=2pt]
			\node (a)[point] at (0,0) {};
			\node (b)[point] at (0,2) {};
			\node (c)[point] at (1.5,0) {};
			\node (d)[point] at (1.5,2) {};
			
			\draw[fill=gray] (a.center) -- (b.center);
			\draw[fill=gray] (c.center) -- (d.center);
		\end{tikzpicture}, and three 2-dimensional facets. Using ADHF we get that $\beta$ is equal to $$\begin{array}{l | rrr }
			3 & 3 & . & .\\ 
			4 & . & 3 & 1\\ 
		\end{array}$$
		which means it is pure, with shift type $(6,5,3)$ and degree type $(1,2)$.
	\end{ex}

	\subsection{Homology Index Sets}
	The PR property is a condition on the homology of links. For this reason it will be useful for us to introduce the following concept and notation.
	\begin{defn}\label{Definition: Homology Index Sets With Links}
		Let $\Delta$ be a simplicial complex.
		\begin{itemize}
			\item We define we define the \textit{homology index set of} $\Delta$ to be the set $h(\Delta)=\{i\in \ZZ : \Hred_i(\Delta) \neq 0\}$.
			\item For a face $\sigma$ in $\Delta$ we define the \textit{homology index set of} $\Delta$ \textit{at} $\sigma$ to be the set $h(\Delta, \sigma)=\{i\in \ZZ : \Hred_i(\lkds) \neq 0\}$. Note that under this definition, we have $h(\Delta)=h(\Delta,\emptyset)$.
			\item For a natural number $m$ we define the \textit{complete homology index set} of $\Delta$ at $m$ as $\hh(\Delta,m)=\bigcup_{\sigma\in \Delta, |\sigma|=m} h(\Delta,\sigma)$.
		\end{itemize}
	\end{defn}
	In many of our later proofs we will add homology index sets $A$ and $B$ together using the following construction.
	\begin{defn}\label{Definition: Boxplus}
		For sets $A$ and $B$, we define $A\boxplus B = \{i+j: i\in A, j\in B\}$.
	\end{defn}
	In particular, note that if $A$ and $B$ are both singletons, then this is the same as adding the elements of those singletons together; and if either set is empty then the resulting sum is also empty.
	
	Using the notation for homology index sets, we can present an alternate definition for PR complexes.
	\begin{prop}\label{Proposition: Alternate PR Definition}
		Let $\Delta$ be a simplicial complex. The following are equivalent.
		\begin{enumerate}
			\item $\Delta$ is PR.
			\item For any $\sigma$ and $\tau$ in $\Delta$ with $|\sigma|\neq |\tau|$, we have $h(\Delta,\sigma)\cap h(\Delta,\tau)=\emptyset$.
			\item For any distinct integers $m_1\neq m_2$, we have $\hh(\Delta,m_1)\cap \hh(\Delta,m_2) = \emptyset$.
		\end{enumerate}
	\end{prop}
	\begin{proof}
		This is a rephrasing of the PR condition in terms of homology index sets.
	\end{proof}
	
	In fact, the complete homology sets of a PR complex can be computed directly from its degree type and offset.
	\begin{prop}\label{Proposition: Alternate PR Definition With Degree Type}
		Let $\Delta$ be a simplicial complex. The following are equivalent.
		\begin{enumerate}
			\item $\Delta$ is a PR complex with degree type $(d_p,\dots,d_1)$ and offset $s$.
			\item $\hh(\Delta,m) = \begin{cases}
				\{r-2\} & \text{ if } m = s + \sum_{j=r}^p d_j \text{ for some } 1\leq r \leq p+1\\
				\emptyset & \text{ otherwise.}\\
			\end{cases}$
		\end{enumerate}
	\end{prop}
	\begin{proof}	
		Let $\Delta$ be a PR complex of the specified degree type and offset. Let $s_p<\dots<s_0$ be as in Definition \ref{Definition: PR Complex Degree Types} (so that $s_p = s$). We know from Corollary \ref{Corollary: PR Complex links have single homology} that every homology index set of $\Delta$ is either a singleton or empty, and from Corollary \ref{Corollary: Homology index sets decreasing sequence} we deduce that the same is true of the \textit{complete} homology index sets. By definition of $s_p,\dots,s_0$, the nonempty complete homology index sets are precisely the sets $\hh(\Delta,s_p),\dots,\hh(\Delta,s_0)$.
		
		Specifically, for each $1\leq r \leq p+1$, the set $\hh(\Delta,s_{r-1})$ is the singleton $\{r-2\}$, and we have $s_{r-1} = s_p + (s_{p-1} - s_p) + \dots +(s_{r-1} - s_{r}) = s + \sum_{j=r}^p d_j$.
		
		Conversely, if $\Delta$ satisfies condition (2), then $\Delta$ must be PR by Proposition \ref{Proposition: Alternate PR Definition}, and we can recover its degree type and offset from the values of $m$ for which $\hh(\Delta,m)$ are nonempty.
	\end{proof}
	
	The following corollary of Proposition \ref{Proposition: Alternate PR Definition With Degree Type} will be particularly crucial to our proof of Theorem \ref{Theorem: PR Complexes of Any Degree Type}.
	\begin{cor}\label{Corollary: PR (1...1) Definition}
		If $\Delta$ has offset $s$ and degree type $\bd = (d_p,\dots,d_i,1,\dots,1)$, then for any natural number $m$ with $s+\sum_{j=i}^p d_j \leq m \leq \dim \Delta + 1$ we have $\hh(\Delta,m)=\{\dim \Delta - m\}$.
	\end{cor}
	\begin{proof}
		For convenience we write $d_1=\dots=d_{i-1}=1$. Note that $\sum \bd = \sum_{j=i}^p d_j +i-1$, and hence by Proposition \ref{Proposition: Dimension of PR Complexes} we have $\dim \Delta = s+ \sum_{j=i}^p d_j + i-2$.
		
		Now fix some $s+\sum_{j=i}^p d_j \leq m \leq \dim \Delta + 1$. In particular, we have $m= s+\sum_{j=i}^p d_j + l$ for some $0\leq l \leq i-1$, and we can rewrite this as $m=s+\sum_{j=i-l}^p d_j$. Proposition \ref{Proposition: Alternate PR Definition With Degree Type} tells us that $\hh(\Delta,m)=\{i-l-2\}$. The result follows because $\dim \Delta - m = i-l-2$.
	\end{proof}

	\section{Intersection Complexes}\label{Section: Intersection Complexes}
	
	In this section we construct an infinite family of PR complexes (and hence a corresponding infinite family of pure resolutions). We begin by noting the following important lemma.
	
	\begin{lem}\label{Lemma: link facet intersections}
		Suppose $\sigma\in \Delta$ is such that $\link_\Delta \sigma$ is not acyclic. Then $\sigma$ is an intersection of facets in $\Delta$.
	\end{lem}
	\begin{proof}
		We prove the contrapositive. Suppose that $\sigma$ is not an intersection of facets of $\Delta$ and let $F_1,\dots,F_m$ be all the facets of $\Delta$ that contain $\sigma$.
		
		We must have $\sigma \subseteq \bigcap_{i=1}^m F_i$. Because $\sigma$ is not an intersection of facets, there must be some vertex $v\in \bigcap_{i=1}^m F_i - \sigma$. This means that $\link_\Delta \sigma$ is a cone over $v$, and is therefore acyclic.
	\end{proof}
	In particular, Lemma \ref{Lemma: link facet intersections} tells us that the only faces of a complex $\Delta$ that contribute towards the Betti numbers of $\Idstar$ are the ones that are intersections of facets. Due to this observation, we might expect many PR complexes to exhibit some kind of symmetry around the points where their facets intersect. The following definition gives us one such symmetry condition on facet intersections.
	\begin{defn}\label{Definition: FIS}
		Let $\Delta$ be a complex with facets $F_1,\dots,F_n$. We say that $\Delta$ is \textit{intersectionally symmetric} if for any $1\leq k\leq n$ and any permutation $\alpha$ in $S_n$, we have
$|F_1\cap \dots \cap F_k|=|F_{\alpha(1)}\cap \dots \cap F_{\alpha(k)}|$.
	\end{defn}
	

\begin{prop}
Let $\Delta$ be an intersectionally symmetric
complex with facets $F_1, \dots, F_n$.
For any $1\leq k\leq n$, and for any permutation $\alpha$ in $S_n$, we have $\link_\Delta (F_1\cap \dots \cap F_k)\cong \link_\Delta (F_{\alpha(1)}\cap \dots \cap F_{\alpha(k)})$.
\end{prop}
\begin{proof}
For each $1 \leq k \leq n$, write $\mu_k$ for the common size of the intersection
of any $k$ facets.

Fix $k$ and a subset  $\{i_1,\ldots,i_k\} \subseteq [n]$;
set $P = F_{i_1} \cap \cdots \cap F_{i_k}$, so $|P| = \mu_k$.
We must show that $\link_\Delta(P) \cong \link_\Delta(Q)$ for 
$Q = F_{1} \cap \cdots \cap F_{k}$.

For any $\ell \in [n] \setminus \{i_1,\ldots,i_k\}$, we have
$|F_{i_1} \cap \cdots \cap F_{i_k} \cap F_\ell| = \mu_{k+1}$.
Hence $F_\ell \supseteq P$ if and only if $\mu_{k+1} = |P| = \mu_k$,
a condition depending only on $k$, not on the specific subset
$\{i_1,\ldots,i_k\}$ or on the choice of $\ell$.

Now if $k = n$, or $k < n$ and $\mu_{k+1} = \mu_k$, 
then every facet contains $P$, so $P = F_1 \cap \cdots \cap F_n$; 
the identical argument applied to $Q$ gives $Q = F_1 \cap \cdots \cap F_n$ as well, so $P=Q$ and we are done.

Assume henceforth that $k < n$ and $\mu_{k+1} < \mu_k$:
now $P$ is contained only in facets  $F_{i_1}, \dots, F_{i_k}$
and $\link_{\Delta} P$ has facets 
$F_{i_1} \setminus P, \dots, F_{i_k} \setminus P$.
Define $A_j= F_{i_j}  \setminus P$ for all $1\leq j\leq k$.

For any vertex $v$ in the set of vertices $V$ of
$A_1 \cup \cdots \cup A_k$ 
define $a(v)=\{ j\in [k] : v\in A_j\}$
and for any non-empty $S\subseteq [k]$ define
$T_S=\{ v\in V : a(v) = S\}$.
These pairwise disjoint sets 
cover
$A_1 \cup \cdots \cup A_k$.

Also, for any non-empty $S\subseteq [k]$ 
$v\in V$ is a vertex of $\bigcap_{j\in S} A_j$
if and only if $a(v)\supseteq S$, hence
\begin{eqnarray*}
 \sum_{S^\prime \supseteq S} \vert T_{S^\prime} \vert  &=&
 \vert \bigcap_{j\in S} A_j \vert\\
 &=& \vert \bigcap_{j\in S} F_{i_j} \setminus P \vert\\
 &=& \left\vert \left(\bigcap_{j\in S} F_{i_j}\right) \right\vert 
 - \vert P \vert\\
 &=&\mu_{\vert S \vert} - \mu_k .
\end{eqnarray*}

Now repeat this argument with
$Q = F_{1} \cap \cdots \cap F_{k}$
with $B_j= F_{j}  \setminus Q$ for all $1\leq j\leq k$.
For any vertex $w$ in the set of vertices $W$ of
$B_1 \cup \cdots \cup B_k$ 
define $b(w)=\{ j\in [k] : w\in B_j\}$
and for any non-empty $S\subseteq [k]$ define
$T^\prime_S=\{ w\in W : b(w) = S\}$.

Now $|T^\prime_S| = |T_S|$ for
every nonempty $S \subseteq [k]$ and we choose
any bijections $T_S \xrightarrow{\;\sim\;} T^\prime_S$ for each $S$ 
resulting in 
 a bijection
$\phi \colon A_1 \cup \cdots \cup A_k \xrightarrow{\;\sim\;} B_1 \cup \cdots \cup B_k$.

Since
\[
	\phi(A_j)  =  \phi\Bigl(\bigcupdot_{S \ni j} T_S\Bigr)
    = \bigcupdot_{S \ni j} \phi(T_S)
    = \bigcupdot_{S \ni j} T^\prime_S
    = B_j ,
\]
the map $\phi$ sends each facet $A_j$ of $\langle A_1,\ldots,A_k\rangle$
to the facet $B_j$ of $\langle B_1,\ldots,B_k\rangle$, and is therefore an
isomorphism of simplicial complexes. Hence
$\link_\Delta(P) = \langle A_1,\ldots,A_k\rangle
\cong \langle B_1,\ldots,B_k\rangle = \link_\Delta(Q)$.
\end{proof}
	
	\subsection{Defining the Complexes}
	In this section, we define the \textit{intersection complexes}, the family of all complexes with intersectional symmetry, along with some examples.
	
	\begin{defn}\label{Definition: Intersection Complexes}
		Let $\bm=(m_1,\dots,m_n)$ be a sequence of nonnegative integers.
		We define the \emph{intersection complex} $\calI(\bm)$ as follows.
		\begin{enumerate}
			\item The vertices of $\calI(\bm)$ are all symbols of the form $v_S^r$ where $S$ is a subset of $[n]$ and $1\leq r \leq m_{|S|}$.
			\item The facets of $\calI(\bm)$ are the sets $F_1,\dots,F_n$ where 
			$$F_j=\left \{v_S^r \,:\, S\subseteq [n], j\in S, 1 \leq r \leq m_{|S|} \right\}.$$
		\end{enumerate}
	\end{defn}
	
	\begin{rem}\label{Remark: Intersection Complex Observations}
		\begin{enumerate}
			\item The integer $m_i$ is equal to the number of vertices contained in the intersection $F_1\cap\dots \cap F_i$ (or any other intersection of $i$ facets) that are not contained in an intersection of $i+1$ facets.
			\item For any subset $S\subseteq[n]$ of size $j$ and any $1\leq i \leq n$, the vertices $v_S^1,\dots,v_S^{m_j}$ are in the facet $F_i$ if and only if we have $i\in S$.
		\end{enumerate}
	\end{rem}
	
	By construction, all intersection complexes have intersectional symmetry. In fact, the family of intersection complexes contains all simplicial complexes with intersectional symmetry.
	
	To see this, suppose that $\Delta$ is a complex with intersectional symmetry and let $F_1,...,F_n$ be the facets of $\Delta$. Define $m_n$ to be the number of vertices in all $n$ facets, and label these vertices $v_{[n]}^1,\dots, v_{[n]}^{m_n}$. Now  define $m_{n-1}$ to be the number of vertices in $F_1\cap \dots \cap F_{n-1}$ but not in $F_n$, and label these vertices $v_{[n-1]}^1,\dots, v_{[n-1]}^{m_{n-1}}$. By symmetry we know that any intersection of $n-1$ facets contains $m_{n-1}$ vertices outside of $F_1\cap \dots \cap F_n$, and we can label each of these accordingly. Now we define $m_{n-2}$ to be the number of vertices in the intersection $F_1\cap \dots \cap F_{n-2}$ which are not in the intersection of a larger number of facets. Proceeding in this way we can find a sequence $(m_1,\dots,m_n)$ of nonnegative integers such that $\Delta$ is isomorphic to the intersection complex $\calI(m_1,\dots,m_n)$.

	\begin{ex}\label{Example: Intersection Complex 1}
		The boundary of the n-simplex in Example \ref{Example: PR Simplex} can be thought of as the intersection complex $\calI(\underbrace{0,\dots,0,1}_{n},0)$. For instance the boundary of the $3$-simplex is the intersection complex $\calI(0,0,1,0)$ as shown below.
		
		\begin{center}
			\begin{tikzpicture}[line join = round, line cap = round]
				
				\coordinate [label=above:{$v_{\{1,2,3\}}^1$}] (4) at (0,{sqrt(2)},0);
				\coordinate [label=left:{$v_{\{1,2,4\}}^1$}] (3) at ({-.5*sqrt(3)},0,-.5);
				\coordinate [label=below:{$v_{\{1,3,4\}}^1$}] (2) at (0,0,1);
				\coordinate [label=right:{$v_{\{2,3,4\}}^1$}] (1) at ({.5*sqrt(3)},0,-.5);
				
				\begin{scope}
					\draw (1)--(3);
					\draw[fill=lightgray,fill opacity=.5] (2)--(1)--(4)--cycle;
					\draw[fill=gray,fill opacity=.5] (3)--(2)--(4)--cycle;
					\draw (2)--(1);
					\draw (2)--(3);
					\draw (3)--(4);
					\draw (2)--(4);
					\draw (1)--(4);
				\end{scope}
			\end{tikzpicture}
		\end{center}
	\end{ex}
	
	\begin{ex}\label{Example: Intersection Complex 2}
		The complex in Example \ref{Example: PR Zelda Symbol} can be thought of as the intersection complex $\calI(1,1,0)$, as shown below.
		
		\begin{center}
			\begin{tikzpicture}[scale = 0.75]
				\tikzstyle{point}=[circle,thick,draw=black,fill=black,inner sep=0pt,minimum width=4pt,minimum height=4pt]
				\node (a)[point,label=above:{$v_{\{1\}}^1$}] at (0,3.4) {};
				\node (b)[point,label=above:{$v_{\{1,2\}}^1$}] at (2,3.4) {};
				\node (c)[point,label=above:{$v_{\{2\}}^1$}] at (4,3.4) {};
				\node (d)[point,label=left:{$v_{\{1,3\}}^1$}] at (1,1.7) {};
				\node (e)[point,label=right:{$v_{\{2,3\}}^1$}] at (3,1.7) {};
				\node (f)[point,label=left:{$v_{\{3\}}^1$}] at (2,0) {};	
				
				\node (g)[label=above: {$F_1$}] at (1,2.2) {};
				\node (h)[label=above: {$F_2$}] at (3,2.2) {};
				\node (i)[label=above: {$F_3$}] at (2,0.5) {};
				
				\begin{scope}[on background layer]
					\draw[fill=gray] (a.center) -- (b.center) -- (d.center) -- cycle;
					\draw[fill=gray] (b.center) -- (c.center) -- (e.center) -- cycle;
					\draw[fill=gray]   (d.center) -- (e.center) -- (f.center) -- cycle;
				\end{scope}
			\end{tikzpicture}
		\end{center}
	\end{ex}
	
	\begin{ex}\label{Example: Intersection Complex (0,1,0,0)}
		The intersection complex $\calI(0,1,0,0)$ is
		$$\begin{tikzpicture}[scale = 1.3]
			\tikzstyle{point}=[circle,thick,draw=black,fill=black,inner sep=0pt,minimum width=3pt,minimum height=3pt]
			\node (a)[point,label=above:$v^1_{\{1,2\}}$] at (0,3.5) {};
			\node (b)[point,label=above:$v^1_{\{1,4\}}$] at (2,2.7) {};
			\node (c)[point,label={[label distance = 0mm]above:$v^1_{\{1,3\}}$}] at (4,3.5) {};
			\node (d)[point,label={[label distance = -2.7mm]below left:$v^1_{\{2,4\}}$}] at (1.6,2.2) {};
			\node (e)[point,label={[label distance = -2.8mm]below right:$v^1_{\{3,4\}}$}] at (2.4,2.2) {};
			\node (f)[point,label=left:$v^1_{\{2,3\}}$] at (2,0) {};
			
			\begin{scope}[on background layer]
				\draw[fill=gray, fill opacity = .8] (d.center) -- (e.center) -- (b.center) -- cycle;
				\draw[fill=gray, fill opacity = .8] (a.center) -- (b.center) -- (c.center) -- cycle;
				\draw[fill=gray, fill opacity = .8] (a.center) -- (f.center) -- (d.center) -- cycle;
				\draw[fill=gray, fill opacity = .8] (c.center) -- (f.center) -- (e.center) -- cycle;
			\end{scope}
		\end{tikzpicture}$$
	\end{ex}
	
	\begin{ex}\label{Example: Intersection Complex Trivial Cases}
		There are two trivial cases of intersection complexes where the defining facets $F_1,\dots,F_n$ are all equal.
		\begin{enumerate}
			\item The complex $\Delta=\calI(\underbrace{0,\dots,0}_n)$ contains no vertices so the sets $F_1,\dots,F_n$ are all empty, which means $\Delta$ is the complex $\{\emptyset\}$.
			\item For any $m>0$, the complex $\Delta=\calI(\underbrace{0,\dots,0}_{n-1},m)$ contains only the $m$ vertices $v_{[n]}^1,...,v_{[n]}^m$, each of which is in every set $F_1,\dots,F_n$. Thus $\Delta$ is the full simplex on these $m$ vertices.
		\end{enumerate}
	\end{ex}
	
	We can find the size of the intersection of any $i$ facets of an intersection complex as follows.
	
	\begin{lem}\label{Lemma: size of sigma_i}
		Let $\bm=(m_1,\dots,m_n)$ be a sequence of non-negative integers, and let $\Delta = \calI(\bm)$ be the corresponding intersection complex with facets $F_1,\dots, F_n$. We have $|F_1 \cap \dots \cap F_i|=\sum_{j=0}^{n-i} \binom{n-i}{j} m_{i+j}$.
	\end{lem}
	\begin{proof}
		Let $S$ be a subset of $[n]$ and $1\leq k \leq m_{|S|}$. From Remark \ref{Remark: Intersection Complex Observations} (2), the vertex $v_S^k$ is contained in the intersection $F_1\cap \dots \cap F_i$ if and only if $S\supseteq [i]$. For each $0\leq j \leq n-i$, there are ${n-i \choose k}$ subsets of $[n]$ containing $[i]$ of size $i+j$ (one for each choice of $j$ elements from the set $\{i+1,...,n\}$). The result follows.
	\end{proof}
	
	We now state our key theorems about intersection complexes, concerning the purity of their corresponding Betti diagrams, and their degree types and Betti numbers. In what follows, we fix a sequence of nonnegative integers $\bm=(m_1,\dots,m_n)$ such that $m_n= 0$ but $\bm\neq 0$. We also let $p$ denote the maximum value of $1\leq i\leq n-1$ for which $m_i\neq 0$. We define $\Delta=\calI(\bm)$ and $\beta=\beta(\Idstar)$.
	
	Note that if $m_n > 0$, then every facet of $\calI(\bm)$ contains the vertices $v_{[n]}^1,\dots,v_{[n]}^{m_n}$, which means $\calI(\bm)$ is a multi-cone over $\calI(m_1,\dots,m_{n-1},0)$, and so the two complexes have the same Betti diagram. The condition $m_n=0$ is therefore harmless.
	
	We wish to prove the following two theorems.
	
	\begin{thm}\label{Theorem: Intersection Complexes Degree Type}
		%
		Let $\bm=(m_1,\dots,m_n)$ be a nonzero sequence of nonnegative integers with $m_n=0$, and define $p=\max\{j\in [n] : m_j \neq 0 \}$. The intersection complex $\Delta=\calI(\bm)$ is PR with degree type $(d_p,\dots,d_1)$ where for each $1\leq i\leq p$, we have $d_i = \sum_{j=0}^{n-i-1} {n-i-1 \choose j} m_{i+j}$.
	\end{thm}
	
	\begin{thm}\label{Theorem: Intersection Complexes Betti Numbers}
		Let $\Delta$ and $\beta$ be as in Theorem \ref{Theorem: Intersection Complexes Degree Type}, and suppose $\beta$ has nonzero Betti numbers $\beta_{0,c_0},\dots.,\beta_{p,c_p}$. We have the following result. 
		\begin{equation*}
			\beta_{i,c_i}=\begin{cases}
				{n \choose i+1}& \text{if } 1\leq i \leq p-1\\
				{n-1 \choose p}& \text{if } i=p\, . 
			\end{cases}
		\end{equation*}
	\end{thm}
	
	As the next lemma shows, the degree types of intersection complexes are prescisely the positive integer sequences $\bs$ for which every difference sequence $\bs^{(r)}$ of $\bs$ is 
	weakly increasing. Here, the $\nth[th]{r}$ difference sequence $\bs^{(r)}=(s^{(r)}_{k-r},\dots,s^{(r)}_{1})$ of a sequence $\bs=(s_k,\dots,s_1)$ is defined recursively via
	\begin{align*}
		s_i^{(0)}&=s_i\\
		s^{(r)}_i&=s^{(r-1)}_{i}-s^{(r-1)}_{i+1}
	\end{align*}
	for each $1 \leq r \leq k-1$ and $1\leq i \leq k-r$.
	\begin{lem}\label{Lemma: Difference Sequences}
		Let $\bs = (s_k,\dots,s_1)$ be a sequence of positive integers. The following are equivalent.
		\begin{enumerate}
			\item $\bs^{(r)}$ is 
			weakly increasing for every $0\leq r \leq k-1$
			\item The leading terms of the difference sequences $\bs^{(1)},\dots, \bs^{(k-1)}$ are nonnegative.
			\item There exists a sequence $(a_1,\dots,a_k)$ of nonnegative integers with $a_k\neq 0$ such that $s_i^{(r)}=\sum_{j=0}^{k-r-i} {k-r-i \choose j} a_{i+j}$ for each $0\leq r \leq k-1$ and each $1\leq i \leq k-r$.
			\item There exists a sequence $(a_1,\dots,a_k)$ of nonnegative integers with $a_k\neq 0$ such that $s_i=\sum_{j=0}^{k-i} {k-i \choose j} a_{i+j}$ for each $1\leq i \leq k$.
		\end{enumerate}
	\end{lem}
	\begin{rem}  
		The equivalence of conditions (1) and (4) in Lemma \ref{Lemma: Difference Sequences} shows that the degree types of intersection complexes given in Theorem \ref{Theorem: Intersection Complexes Degree Type} are precisely those positive integer sequences for which every difference sequence is 
		weakly increasing. Indeed, for any such sequence $\bs=(s_k,\dots,s_1)$ with corresponding nonnegative sequence $(a_1,\dots,a_k)$, the complex $\calI(a_1,\dots,a_k,0)$ has degree type $\bs$. Conversely, let $\bm=(m_1,\dots,m_n)$ be a nonzero sequence of nonnegative integers with $m_n=0$, and define a sequence $\widetilde{\bd}=(d_{n-1},\dots,d_1)$ via $d_i=\sum_{j=0}^{n-i-1} {n-i-1 \choose j} m_{i+j}$ for each $1\leq i \leq n-1$. By Lemma \ref{Lemma: Difference Sequences}, $\widetilde{\bd}^{(r)}$ is 
		weakly increasing for each $r$, and hence the same is true for the degree type $\bd=(d_p,\dots,d_1)$ of the complex $\calI(\bm)$.
	\end{rem}
	\begin{proof}
		(1) $\Rightarrow$ (2) is immediate, because if any term in $\bs^{(r)}$ is negative, then the sequence $\bs^{(r-1)}$ cannot be 
		weakly increasing.
		
		We now move on to (2) $\Rightarrow$ (3). For each $1\leq r \leq k$ we set $a_r$ to be the leading term $s_{r}^{(k-r)}$ of the difference seqence $\bs^{(k-r)}$. By assumption each of these terms is nonnegative, and because $\bs$ consists of positive integers, we have $a_k=s_k\neq 0$. We aim to show that $s_i^{(r)}= \sum_{j=0}^{k-r-i} {k-r-i \choose j} a_{i+j}$.
		
		We proceed by reverse induction on both $k-1\geq r \geq 0$ and $k-r\geq i \geq 0$. The base cases $r=k-1$ and $i=k-r$ come from the definitions of $a_1,\dots,a_k$. For the inductive case $r < k-1$ and $i < k-r$, we have
		\begin{align*}
			s_{i}^{(r)}&= s_{i}^{(r+1)}+s_{i+1}^{(r)}\\
			&= \sum_{j=0}^{k-r-i-1} {k-r-i-1 \choose j} a_{i+j} + \sum_{j=0}^{k-r-i-1} {k-r-i-1 \choose j} a_{i+j+1}\\
			&= \sum_{j=0}^{k-r-i}\left({k-r-i-1 \choose j} + {k-r-i-1 \choose j-1}\right) a_{i+j}\\
			&= \sum_{j=0}^{k-r-i} {k-r-i \choose j} a_{i+j}.
		\end{align*}
		
		For (3) $\Rightarrow$ (1), we note that for any $0\leq r \leq k-2$ and any $1\leq i \leq k-r$, the term $s^{(r+1)}_i$ is a sum of nonnegative integers, and thus nonnegative itself. This means that $s^{(r)}_i\geq s^{(r)}_{i+1}$.

		It only remains to show that conditions (3) and (4) are equivalent. (3) $\Rightarrow$ (4) is immediate. For (4) $\Rightarrow$ (3), we proceed by induction on $0\leq r \leq k$. The base case $r=0$ is immediate, and the inductive case $r>0$ follows from the equation $s_{i}^{(r)} = s_i^{(r-1)}-s_{i+1}^{(r-1)}$.

	\end{proof}
	
	Our proof for Theorems \ref{Theorem: Intersection Complexes Degree Type} and \ref{Theorem: Intersection Complexes Betti Numbers} proceeds as follows. First, we make use of some deformation retractions to allow us to restrict our attention to intersection complexes of a particularly simple form. Next we show that these simple intersection complexes each have only a single nontrivial homology group. And then we assemble these pieces together to show that all intersection complexes are PR complexes with the desired degree types and Betti numbers.
	
	\subsection{Deformation Retractions and Links}
	For the rest of this section we let $\eseq{n}{i}$ denote the sequence $(\overbrace{\underbrace{0,\dots,0, 1}_{i},0,\dots,0}^n)$ (that is, the sequence of length $n$ whose only nonzero term is a $1$ at position $i$). We will also use $\eseq{n}{0}$ to denote the sequence of $n$ zeroes.
	
	The following result shows that all intersection complexes deformation retract on to an intersection complex of the form $\calI(\eseq{j}{i})$ for some $i$ and $j$.
	
	\begin{prop}\label{Proposition: Def Retract of Intersection Complex}
		Let $\bm=(m_1,\dots,m_n)$ be a nonzero sequence in $\ZZ_{\geq 0}^n$ with $m_n=0$, and define $p=\max \{i\in [n-1]: m_i \neq 0\}$. We have a deformation retraction $\calI(\bm)\leadsto \calI(\eseq{n}{p})$.
	\end{prop}
	\begin{proof}
		Let $v_S^r$ be a vertex of $\calI(\bm)$, with $S$ a subset of $[n]$ of size less than or equal to $p$, and $1 \leq r \leq m_{|S|}$. Choose any subset $S'\subseteq [n]$ of size $p$ which contains $S$. The vertex $v_{S'}^1$ lies inside $\calI(\eseq{n}{p})$, and any facet of $\calI(\bm)$ containing $v_S^r$ must also contain $v_{S'}^1$. Thus by Corollary \ref{Corollary: Deformation Retract} there is a deformation retraction $\calI(\bm) \leadsto \calI(\eseq{n}{p})$ obtained by 
		collapsing every vertex $v_S^r$ in $\calI(\bm)$ onto an appropriate vertex $v_{S'}^1$ in $\calI(\eseq{n}{p})$.
	\end{proof}
	\begin{ex}\label{Example: Intersection complex def retract}
		For the intersection complex $\calI(1,1,0)$ in Example \ref{Example: Intersection Complex 2} we have the deformation retraction
		\begin{center}
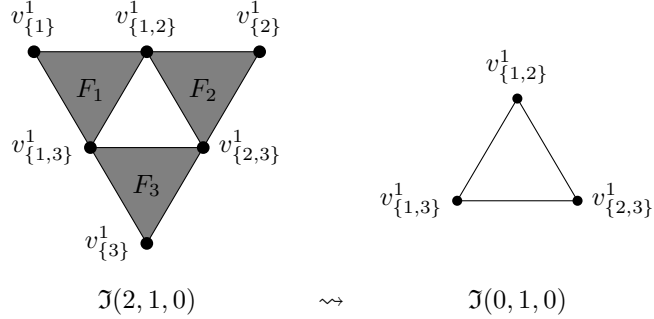
			\begin{tabular}{ccc}
				\begin{tikzpicture}[scale = 0.75]
					\tikzstyle{point}=[circle,thick,draw=black,fill=black,inner sep=0pt,minimum width=4pt,minimum height=4pt]
					\node (a)[point,label=above:{$v_{\{1\}}^1$}] at (0,3.4) {};
					\node (b)[point,label=above:{$v_{\{1,2\}}^1$}] at (2,3.4) {};
					\node (c)[point,label=above:{$v_{\{2\}}^1$}] at (4,3.4) {};
					\node (d)[point,label=left:{$v_{\{1,3\}}^1$}] at (1,1.7) {};
					\node (e)[point,label=right:{$v_{\{2,3\}}^1$}] at (3,1.7) {};
					\node (f)[point,label=left:{$v_{\{3\}}^1$}] at (2,0) {};	
					
					\node (g)[label=above: {$F_1$}] at (1,2.2) {};
					\node (h)[label=above: {$F_2$}] at (3,2.2) {};
					\node (i)[label=above: {$F_3$}] at (2,0.5) {};
					
					\begin{scope}[on background layer]
						\draw[fill=gray] (a.center) -- (b.center) -- (d.center) -- cycle;
						\draw[fill=gray] (b.center) -- (c.center) -- (e.center) -- cycle;
						\draw[fill=gray]   (d.center) -- (e.center) -- (f.center) -- cycle;
					\end{scope}
				\end{tikzpicture} & & \begin{tikzpicture}[scale=0.8]
					\tikzstyle{point}=[circle,thick,draw=black,fill=black,inner sep=0pt,minimum width=3pt,minimum height=3pt]
					\node (a)[point, label=left:$v_{\{1,3\}}^1$] at (0,0) {};
					\node (b)[point, label=right:$v_{\{2,3\}}^1$] at (2,0) {};
					\node (c)[point, label=above:$v_{\{1,2\}}^1$] at (1,1.7) {};
					
					\node at (0,-1) {};
					
					\draw (a.center) -- (b.center) -- (c.center) -- cycle;
				\end{tikzpicture}\\ 
				$\calI(2,1,0)$& $\rightsquigarrow$ & $\calI(0,1,0)$\\
			\end{tabular}
		\end{center}
	\end{ex}
	
	In fact, not only is there a deformation retraction from the intersection complex $\calI(\bm)$ itself onto a complex of the form $\calI(\eseq{j}{i})$, but the link of any intersections of facets in $\calI(\bm)$ also deformation retracts onto a complex of this form. By symmetry it suffices to consider facet intersections of the form $\sigma_i=F_1\cap \dots \cap F_i$.
	\begin{prop}\label{Proposition: Def Retract of Intersection Complex Link}
		Let $\bm=(m_1,\dots,m_n)$ be a nonzero sequence in $\ZZ_{\geq 0}^n$ with $m_n=0$, and let $\Delta = \calI(\bm)$ be the corresponding intersection complex with facets $F_1,\dots, F_n$. Suppose $\sigma_i = F_1 \cap \dots \cap F_i$ for some $1\leq i \leq p$. We have a deformation retraction $\lkds_i\leadsto \calI(\eseq{i}{i-1})$.
	\end{prop}
	\begin{proof}
		The complex $\calI(\eseq{i}{i-1})$ contains precisely those vertices of the form $v^1_{[i]-\{j\}}$ for some $1 \leq j \leq i$. All of these vertices are contained in exactly $i-1$ of the facets $F_1,\dots,F_i$, which means none of them is contained in $\sigma_i$ and all are therefore contained in $\lkds_i = \langle F_1 - \sigma_i,\dots, F_i-\sigma_i \rangle$.
		
		Let $v^r_S$ be any vertex in $\lkds_i$, with $S$ a subset of $[n]$ and $1 \leq r \leq m_{|S|}$. Because $v^r_S$ does not lie in $\sigma_i$, there must be some $1\leq j \leq i$ for which $v^r_S$ does not lie in the facet $F_j$. This means that we have $j\notin S$, and hence $S\cap [i] \subseteq [i]-\{j\}$. Thus every facet of $\lkds_i$ containing $v^r_S$ must also contain $v^1_{[i]-\{j\}}$.
		
		Thus by Corollary \ref{Corollary: Deformation Retract} there is a deformation retraction $\lkds_i\leadsto \calI(\eseq{i}{i-1})$ given by
		collapsing every vertex $v_S^r$ in $\lkds_i$ onto an appropriate vertex $v^1_{[i]-\{j\}}$ in $\calI(\eseq{i}{i-1})$.
	\end{proof}
	\begin{ex}\label{Example: Intersection complex link def retract}
		For the intersection complex $\Delta=\calI(1,1,0)$ in Example \ref{Example: Intersection Complex 2} we have $\sigma_2=F_1\cap F_2 = \{v^1_{\{1,2\}}\}$, and we get the following deformation retraction.
		\begin{center}
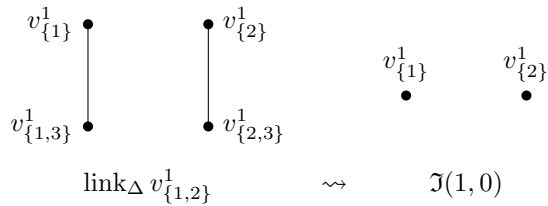
			\begin{tabular}{ccc}
				\begin{tikzpicture}[scale = 0.8]
					\tikzstyle{point}=[circle,thick,draw=black,fill=black,inner sep=0pt,minimum width=3pt,minimum height=3pt]
					\node (a)[point,label=left:$v_{\{1\}}^1$] at (1,3.4) {};
					\node (d)[point,label=left:$v_{\{1,3\}}^1$] at (1,1.7) {};
					
					\node (c)[point,label=right:$v_{\{2\}}^1$] at (3,3.4) {};
					\node (e)[point,label=right:$v_{\{2,3\}}^1$] at (3,1.7) {};
					
					\draw(a.center)--(d.center);
					\draw (c.center)--(e.center);
				\end{tikzpicture}  & & \begin{tikzpicture}[scale=0.8]
					\tikzstyle{point}=[circle,thick,draw=black,fill=black,inner sep=0pt,minimum width=3pt,minimum height=3pt]
					\node (a)[point, label=above:$v_{\{1\}}^1$] at (1,2.5) {};
					\node (b)[point, label=above:$v_{\{2\}}^1$] at (3,2.5) {};
					
					\node at (1,1.7) {};
					
				\end{tikzpicture}\\ 
				$\lkds[v^1_{\{1,2\}}]$& $\rightsquigarrow$ & $\calI(1,0)$\\
			\end{tabular}
		\end{center}
	\end{ex}
	
	These two deformation retractions allow us to restrict our attention solely to the homologies of the complexes $\calI(\eseq{n}{p})$ for $p=1,\dots,n-1$. To find the homologies, we will make use of the following lemma.
	
	Note that the vertex set of $\calI(\eseq{n}{p})$ contains exactly one vertex $v_S^1$ for each subset $S$ in $[n]$ of size $p$. For ease of notation, we label this vertex as $v_S$ rather than $v_S^1$.
	\begin{lem}\label{Lemma: subcomplexes e^n_p}
		Let $\Delta$ be the complex $\calI(\eseq{n}{p})$ for some $1\leq p \leq n-1$, with facets $F_1,...,F_n$.
		
		Define $A$ to be the subcomplex of $\Delta$ with facets $F_1,...,F_{n-1}$ and $B$ to be the subcomplex of $\Delta$ with facet $F_n$. We have the following results.
		\begin{enumerate}
			\item There is a deformation retraction from $A$ on to $\calI(\eseq{n-1}{p})$.
			\item $A\cap B$ is homeomorphic to $\calI(\eseq{n-1}{p-1})$.
		\end{enumerate}
	\end{lem}
	\begin{proof}
		For part (1), note that $A$ contains the complex $\calI(\eseq{n-1}{p})$, and the only vertices in $A$ outside of $\calI(\eseq{n-1}{p})$ are those of the form $v_S$ where $S$ is a subset in $[n]$ of size $p$ containing $n$. For any such $S$, we may choose some $l\in [n] - S$. All of the facets $F_1,...,F_{n-1}$ containing the vertex $v_S$ must also contain the vertex $v_{S\cup\{l\}-\{n\}}$, which lies inside $\calI(\eseq{n-1}{p})$, and hence we may 
		collapse the vertex $v_S$ onto the vertex $v_{S\cup\{l\}-\{n\}}$. 
		This gives us a deformation retraction from $\Delta$ on to its subcomplex $\calI(\eseq{n-1}{p})$.
		
		For part (2), note that $A \cap B = \langle F_1\cap F_n,....,F_{n-1} \cap F_n \rangle$. Hence the vertices in $A\cap B$ are all those of the form $v_S$ where $S$ is a subset in $[n]$ of size $p$ containing $n$. Thus there is a bijection from the vertex set of $A\cap B$ to the vertex set of $\calI(\eseq{n-1}{p-1})$ given by $v_S\mapsto v_{S-\{n\}}$. This bijection takes each facet $F_i\cap F_n$ to a facet of $\calI(\eseq{n-1}{p-1})$, and hence it is a homeomorphism.
	\end{proof}
	
	Using Lemma \ref{Lemma: subcomplexes e^n_p}, we can find the homology of the complex $\calI(\eseq{n}{p})$.
	\begin{prop}\label{Proposition: homology of e^n_p}
		For all $0\leq p\leq n-1$, the complex $\calI(\eseq{n}{p})$ has only $\nth[st]{(p-1)}$ homology, of dimension ${n-1 \choose p}$.
	\end{prop}
	\begin{proof}
		We proceed by induction on $n\geq 1$ and $p\geq 0$. For the base cases, if $p=0$, the complex $\calI(\eseq{n}{p})=\{\emptyset\}$ has only $\nth[st]{(-1)}$ homology, of dimension $1$. Similarly if $n=1$, then we must have $p=0$, so the same reasoning applies.
		
		Now suppose $n\geq 2$ and $1\leq p \leq n-1$, and let $\Delta=\calI(\eseq{n}{p})$, with facets $F_1,...,F_n$. We define the subcomplexes $A$ and $B$ of $\Delta$ as in Lemma \ref{Lemma: subcomplexes e^n_p}. Note that $B$ consists of a single facet and is therefore acyclic.
		
		By Lemma \ref{Lemma: subcomplexes e^n_p}, $A$ deformation retracts on to $\calI(\eseq{n-1}{p})$. If $p=n-1$, then $\calI(\eseq{n-1}{p})$ is a full simplex and is therefore acyclic. Otherwise, by the inductive hypothesis, it has only $(p-1)^\text{st}$ homology, of dimension ${n-2 \choose p}$. Either way we have that $\dim_\KK \Hred_{p-1} (\calI(\eseq{n-1}{p})) = {n-2\choose p}$ and all other homologies are zero.
		
		We also have that $A\cap B$ is homeomorphic to $\calI(\eseq{n-1}{p-1})$. By the inductive hypothesis, this has only $(p-2)^\text{nd}$ homology, of dimension ${n-2 \choose p-1}$.
		
		Thus the Mayer-Vietoris Sequence yields an exact sequence
		$$0\rightarrow \Hred_{p-1}(A)\rightarrow \Hred_{p-1}(\Delta)\rightarrow\Hred_{p-2}(A\cap B)\rightarrow 0$$
		which means that $\Delta$ has only $(p-1)^\text{st}$ homology, and this homology has dimension ${n-2 \choose p}+{n-2\choose p-1}={n-1 \choose p}$.
	\end{proof}
	
	\subsection{Proving Theorems \ref{Theorem: Intersection Complexes Degree Type} and \ref{Theorem: Intersection Complexes Betti Numbers}}
	
	We now have all the ingredients we need to prove Theorems \ref{Theorem: Intersection Complexes Degree Type} and \ref{Theorem: Intersection Complexes Betti Numbers}. We prove both theorems together below.
	
	\begin{proof}[Proof of Theorems \ref{Theorem: Intersection Complexes Degree Type} and \ref{Theorem: Intersection Complexes Betti Numbers}]
		As in Proposition \ref{Proposition: Def Retract of Intersection Complex Link}, we define $\sigma_i$ to be the intersection of facets $F_1\cap \dots \cap F_i$ for $1\leq i \leq p$. We also define $\sigma_{p+1}=F_1\cap \dots \cap F_{p+1}=\emptyset$. Propositions \ref{Proposition: Def Retract of Intersection Complex}, \ref{Proposition: Def Retract of Intersection Complex Link} and \ref{Proposition: homology of e^n_p} show that for $1\leq i\leq p+1$, the link $\lkds_i$ has only $(i-1)^\text{st}$ homology.
		
		We wish to show that $\Delta=\calI(\bm)$ is a PR complex with the desired degree type and Betti numbers.
		
		To see that $\Delta$ is a PR complex, suppose we have two faces $\tau_1$ and $\tau_2$ of $\Delta$ whose links both have nontrivial homology at the same degree. By Lemma \ref{Lemma: link facet intersections}, $\tau_1$ and $\tau_2$ must be intersections of facets. By the intersectional symmetry of $\Delta$ we may assume that $\tau_1=\sigma_{j_1}$ and $\tau_2 = \sigma_{j_2}$ for some $1\leq j_1,j_2\leq p+1$, and thus their links have $(j_1-1)^\text{st}$ homology and $(j_2-1)^\text{st}$ homology respectively. This means $j_1=j_2$, and in particular $|\tau_1|=|\tau_2|$.
		
		We now move on to the degree type. For any $1\leq i \leq p$, the value $d_i$ as in Definition \ref{Definition: PR Complex Degree Types} is given by $d_i = |\sigma_i|-|\sigma_{i+1}|$. Proposition \ref{Lemma: size of sigma_i} tells us that $|\sigma_i| = \sum_{j=0}^{n-i} {n-i \choose j} m_{i+j}$, and hence we have
		\begin{align*}
			d_i &= \sum_{j=0}^{n-i} {n-i \choose j} m_{i+j} - \sum_{j=0}^{n-i-1} {n-i-1 \choose j} m_{i+j+1}\\
			&= \sum_{j=0}^{n-i} {n-i \choose j} m_{i+j} - \sum_{j=0}^{n-i} {n-i-1 \choose j-1} m_{i+j}\\
			&= \sum_{j=0}^{n-i} ({n-i \choose j}-{n-i-1 \choose j-1}) m_{i+j}\\
			&= \sum_{j=0}^{n-i} {n-i-1 \choose j} m_{i+j}
		\end{align*}
		
		Finally we look at the Betti numbers $\beta_{0,c_0},\dots,\beta_{p,c_p}$. First we note that $\beta_{p,c_p}$ is equal to $\dim_{\KK}\Hred_{p-1}(\Delta)$, which is ${n-1\choose p}$ by Propositions \ref{Proposition: Def Retract of Intersection Complex} and \ref{Proposition: homology of e^n_p}. Let $0\leq i \leq p-1$. For any face $\tau$ of $\Delta$, the complex $\link_\Delta \tau$ has nontrivial $(i-1)^\text{st}$ homology if and only if it is an intersection of $i+1$ facets, in which case this homology has dimension ${i-1\choose i-1}=1$ by Propositions \ref{Proposition: Def Retract of Intersection Complex Link} and \ref{Proposition: homology of e^n_p}. Thus, by ADHF, the Betti number $\beta_{i,c_i}$ is equal to the number of intersections of $i+1$ facets in $\Delta$, which is ${n\choose i+1}$.
	\end{proof}

	\section{Partition Complexes}\label{Section: Partition Complexes}
	
	In this section we construct another infinite family of PR complexes, which we call \textit{partition complexes}. Just as every intersection complex has a corresponding sequence of nonnegative integers $\bm$, every partition complex has three corresponding integers $a$, $p$ and $m$, and we denote them accordingly as $\calP(a,p,m)$. We show that when $a\geq 2$, and $1\leq m\leq p$ the partition complex $\calP(a,p,m)$ has degree type $(\overbrace{1,\dots,1\underbrace{a,1,\dots,1}_{m}}^{p})$.
		
		Just like intersection complexes, partition complexes can be seen as generalisations of the boundary complexes of simplices. Specifically, we saw in the last section that for any sequence of nonnegative integers $\bm=(m_1,\dots,m_{p-1},0)$ with $m_{p-1}\neq 0$, the intersection complex $\calI(\bm)$ deformation retracts onto the boundary of the $p$-simplex. The same is also true for the partition complex $\calP(a,p,m)$ for any integers $a\geq 2$ and $1\leq m \leq p$.
		
		\subsection{Defining the Complexes}
		In this section we define the family of partition complexes along with some examples, and present some preliminary results about them. For any integer $p\geq 0$ the partition complex $\calP(a,p,m)$ admits a natural symmetry under the action of the symmetric group $\Sym\{0,\dots,p\}$ (which is isomorphic to $S_{p+1}$). For notational convenience we will denote the group $\Sym\{0,\dots,p\}$ by $S^0_p$.
		
		We also introduce the following notation for complexes generated by group actions.
		\begin{defn}\label{Definition: Complexes from Groups}
			Suppose we have a group $G$ acting on a set $V$. For any subset $F\subseteq V$ and any element $g\in  G$, we define the subset $gF\subseteq V$ to be the set $\{gx:x\in F\}$. For subsets $F_1,\dots,F_m$ in $V$, and a group $G$ acting on $V$, we use $\langle F_1,\dots, F_m\rangle_G$ to denote the complex $\langle gF_i|1\leq i\leq m, g\in G\rangle$.
		\end{defn}
		
		It is worth noting that our main theorem on the degree types of the partition complex $\calP(a,p,m)$ (Theorem \ref{Theorem: Partition Complex Degree Type}) holds only for integers $a\geq 2$ and $1\leq m \leq p$. However, we define our complexes below slightly more broadly, to include the additional cases $p\in \{-1,0\}$ and $m\in \{0,p+1\}$. We do not care about these fringe cases of partition complexes for their own sake, but their construction will be crucial to our proof of Theorem \ref{Theorem: Partition Complex Degree Type}, because they occur as links in the partition complexes that we do care about.
		
		\begin{defn}\label{Definition: Partition Complex Vertex Set}
			For two integers $a\geq 2$ and $p\geq -1$. We define the vertex set $V_p^a$ to be the set consisting of vertices of the form $x_i$ and $y_i^j$ for $0\leq i\leq p$ and $1\leq j \leq a-1$. For convenience, we often write $y_i^1$ simply as $y_i$.
			
			We will sometimes partition $V_p^a$ into subsets $X_p \sqcup Y_p^a$ with $X_p=\{x_i\}_{0\leq i \leq p}$ and $Y_p^a =\{y_i^j\}_{0\leq i \leq p, 1\leq j \leq a-1}$. For reasons that will become apparent we refer to the vertices in $X_p$ as \textit{boundary vertices} and the vertices in $Y_p^a$ as \textit{partition vertices}. We sometimes make a further distinction between those partition vertices $y_i^j$ for which $j\geq 2$ and those for which $j=1$, by referring to them respectively as \textit{upper} and \textit{lower} partition vertices.
		
		For $p\geq 0$ the symmetry group $S^0_p=\Sym\{0,\dots,p\}$ acts on $V_p^a$ via the action
		\begin{align*}
			\sigma(x_i)&=x_{\sigma(i)}\\
			\sigma(y_i^j)&=y_{\sigma(i)}^j\\
		\end{align*}
	\end{defn}
	
	\begin{defn}\label{Definition: Partition Complex Generating Sets}
		Let $a$, $p$ and $i$ be integers with $a\geq 2$, $p\geq -1$ and $1\leq i \leq p+1$.
		
		Let $\PP(a,i)$ denote the set of partitions $\lambda = (\lambda_1,...,\lambda_i)$ of $a+i-2$ into $i$ parts (i.e. $\lambda_1\geq \dots \geq \lambda_i>0$ and $\sum_{j=1}^i\lambda_j=a+i-2$).
		
		For $\lambda = (\lambda_1,...,\lambda_i) \in \PP(a,i)$, we define the $\lambda$-\textit{generating set} $G^\lambda_{p,i}$ as follows.
		\begin{itemize}
			\item $G_{p,i}=\{x_i,\dots,x_p\}\subset X_p$.
			\item $G_{p,i}^\lambda = G_{p,i}\sqcup \{y_r^j:0\leq r\leq i-1, 1\leq j\leq \lambda_{r+1}\}\subset V_p^a$.
		\end{itemize}
	\end{defn}
	
	\begin{rem}\label{Remark: Size of Facets of Partition Complexes}
		Note that \begin{align*}
			|G^\lambda_{p,i}|&=|G_{p,i}|+\sum_{j=1}^1 \lambda_j\\
			&=(p-i+1)+(a+i-2)\\
			&=a+p-1.
		\end{align*}
	\end{rem}
	
	\begin{defn}\label{Definition: Partition Complexes}
		Let $a$ and $p$ be integers with $a\geq 2$, $p\geq -1$. For a third integer $1\leq m \leq p+1$ we define the \textit{partition complex} $\calP(a,p,m)$ on vertex set $V_p^a$ to be
		$$\calP(a,p,m) = \left \langle G^\lambda_{p,i} : 1\leq i \leq m, \lambda \in \PP(a,i) \right \rangle_{S^0_p},$$
		using the notation laid out in Definition \ref{Definition: Complexes from Groups}. We also define the partition complex $\calP(a,p,0)$ to be the complex $\{\emptyset\}$ on vertex set $V_p^a$.
	\end{defn}
	
	By definition, partition complexes are symmetric under the action of $S^0_p$.
	
	\begin{rem}\label{Remark: Partitions of a-2}
		For every $1\leq i \leq m$ and every partition $\lambda\in \PP(a,i)$, the $\lambda$-generating set $G^\lambda_{p,i}$ contains the lower partition vertices $y^1_1,\dots,y^1_i$. This means that it must also contain exactly $a-2$ upper partition vertices. Thus for $m\geq 1$ and $p\geq 0$, every facet of the partition complex $\calP(a,p,m)$ contains exactly $a-2$ upper partition vertices.
	\end{rem}
	
	For explanatory purposes in some of our later proofs, we will sometimes represent the faces of partition complexes in \textit{grid notation} as below.
	\begin{notation}\label{Notation: Grid Notation For Partition Complex Faces}
		Let $F$ be a subset of the vertex set $V^a_p$. We sometimes denote $F$ pictorially as a partially shaded grid of size $(a-1)\times (p+1)$. The cells in the grid represent the vertices in $V^a_p$, as shown below
		\begin{center}
			\begin{tikzpicture}[xscale=1,yscale=.7]
				\draw[xstep=1cm,ystep=1,black,very thin] (0,3) grid (3,1);
				\draw[step=1cm,ystep=1,gray,dotted] (0,1) grid (3,0);
				\draw (0.5,2.5) node{$x_0$};
				\draw (1.5,2.5) node{$\dots$};
				\draw (2.5,2.5) node{$x_p$};
				\draw (0.5,1.5) node{$y^1_0$};
				\draw (1.5,1.5) node{$\dots$};
				\draw (2.5,1.5) node{$y^1_p$};
				\draw (0.5,0.5) node{$y^2_0$};
				\draw (1.5,0.5) node{$\dots$};
				\draw (2.5,0.5) node{$y^2_p$};
				\draw (0.5,-.5) node{$\vdots$};
				\draw (1.5,-.5) node{$\vdots$};
				\draw (2.5,-.5) node{$\vdots$};
			\end{tikzpicture}
		\end{center}
		and the shaded cells denote the vertices that are elements of $F$.
		
		For instance, we denote the set $G_{3,2}^{(2,1)}=\{y_0^1,y_0^2,y_1^2,x_2,x_3\}$ in $V^3_2$ by the shaded grid
		\begin{center}
			\begin{tikzpicture}[xscale=1,yscale=.7]
				\fill[gray] (0,3) rectangle (1,2);
				\draw (0.5,2.5) node{$y^1_0$};
				\fill[gray] (0,2) rectangle (1,1);
				\draw (0.5,1.5) node{$y^2_0$};
				\fill[gray] (1,3) rectangle (2,2);
				\draw (1.5,2.5) node{$y^1_1$};
				\fill[gray] (2,4) rectangle (3,3);
				\draw (2.5,3.5) node{$x_2$};
				\fill[gray] (3,4) rectangle (4,3);
				\draw (3.5,3.5) node{$x_3$};
				\draw[step=1cm,black,very thin] (0,4) grid (4,2);
				\draw[step=1cm,gray,dotted] (0,2) grid (4,1);
			\end{tikzpicture}
		\end{center}
	\end{notation}
	
	We now consider some special cases and examples to help illustrate the construction os partition complexes. We begin with the cases for small values of $p$.
	\begin{ex}
		\underline{\textbf{The case $p=-1$:}}
		
		The vertex set $V^a_{-1}$ is empty, and we have $\calP(a,-1,0)=\{\emptyset\}$, by definition.
	\end{ex} 
	\begin{ex}
		\underline{\textbf{The case $p=0$:}}
		
		The vertex set $V^a_0$ is equal to $\{x_0,y_0^1,\dots,y_0^{a-1}\}$. The complex  $\calP(a,0,0)$ is the irrelevant complex $\{\emptyset\}$ just as above, and we also have $$\calP(a,0,1)=\langle G^{(a-1)}_{0,1} \rangle = \langle \{y_0^1,\dots,y_0^{a-1}\}\rangle.$$
	\end{ex}
	
	\begin{ex}\label{Example: Partition Complex p=m=1}
		\underline{\textbf{The case $p=m=1$:}} 
		
		The partition complex $\calP(a,1,1)$ has precisely two facets, namely the facets $\{y_0^1,\dots,y_0^{a-1},x_1\}$ and $\{y_1^1,\dots,y_1^{a-1},x_0\}$. Thus it comprises of two disjoint $(a-1)$-simplices, so it is PR with degree type $(a)$.
	\end{ex}
	
	In the above example, there is a deformation retraction from the complex $\calP(a,1,1)$ on to the boundary $\partial \langle X_1 \rangle$ of the full simplex on $X_1=\{x_0,x_1\}$. In fact we will see in the next section that for \textit{any} $1\leq m \leq p$, the partition complex $\calP(a,p,m)$ always deformation retracts onto the boundary $\partial \langle X_p \rangle$ of the full simplex on $X_p$.
	\begin{ex}\label{Example: Partition Complex a=2}
		\underline{\textbf{The case $a=2$:}}
		
		The vertex set $V^2_p$ is the set $\{x_0,\dots,x_p,y_0,\dots,y_p\}$. For any positive integer $i$, there is only one partition of $2+i-2=i$ into $i$ parts, namely the partition $\underbrace{(1,\dots,1)}_i$. Thus for any integers $m$ and $p$ with $1 \leq m\leq p+1$, the generating facets of $\calP(2,p,m)$ are the sets
		\begin{align*}
			G_{p,1}^{(1)}&= \{y_0,x_1,\dots,x_p\} &= \begin{tikzpicture}[xscale=.8,yscale=.4]
				\fill[gray] (0,1) rectangle (1,0);
				\draw (0.5,0.5) node{$y_0$};
				\fill[gray] (1,2) rectangle (2,1);
				\draw (1.5,1.5) node{$x_1$};
				\fill[gray] (2,2) rectangle (3,1);
				\draw (2.5,1.5) node{$\dots$};
				\fill[gray] (3,2) rectangle (4,1);
				\draw (3.5,1.5) node{$\dots$};
				\fill[gray] (4,2) rectangle (5,1);
				\draw (4.5,1.5) node{$\dots$};
				\fill[gray] (5,2) rectangle (6,1);
				\draw (5.5,1.5) node{$\dots$};
				\fill[gray] (6,2) rectangle (7,1);
				\draw (6.5,1.5) node{$x_p$};
				\draw[step=1cm,black,very thin] (0,2) grid (7,0);
			\end{tikzpicture}&\\
			G_{p,2}^{(1,1)}&= \{y_0,y_1,x_2\dots,x_p\} &= \begin{tikzpicture}[xscale=.8,yscale=.4]
				\fill[gray] (0,1) rectangle (1,0);
				\draw (0.5,0.5) node{$y_0$};
				\fill[gray] (1,1) rectangle (2,0);
				\draw (1.5,0.5) node{$y_1$};
				\fill[gray] (2,2) rectangle (3,1);
				\draw (2.5,1.5) node{$x_2$};
				\fill[gray] (3,2) rectangle (4,1);
				\draw (3.5,1.5) node{$\dots$};
				\fill[gray] (4,2) rectangle (5,1);
				\draw (4.5,1.5) node{$\dots$};
				\fill[gray] (5,2) rectangle (6,1);
				\draw (5.5,1.5) node{$\dots$};
				\fill[gray] (6,2) rectangle (7,1);
				\draw (6.5,1.5) node{$x_p$};
				\draw[step=1cm,black,very thin] (0,2) grid (7,0);
			\end{tikzpicture}&\\
			&\text{         } \vdots\\
			G_{p,m}^{(1,\dots,1)}&=\{y_0,\dots,y_{m-1},x_m,\dots,x_p\}&= \begin{tikzpicture}[xscale=.8,yscale=.4]
				\fill[gray] (0,1) rectangle (1,0);
				\draw (0.5,0.5) node{$y_0$};
				\fill[gray] (1,1) rectangle (2,0);
				\draw (1.5,0.5) node{$\dots$};
				\fill[gray] (2,1) rectangle (3,0);
				\draw (2.5,0.5) node {$\dots$};
				\fill[gray] (3,1) rectangle (4,0);
				\draw (3.5,0.5) node{$y_{m-1}$};
				\fill[gray] (4,2) rectangle (5,1);
				\draw (4.5,1.5) node{$x_m$};
				\fill[gray] (5,2) rectangle (6,1);
				\draw (5.5,1.5) node{$\dots$};
				\fill[gray] (6,2) rectangle (7,1);
				\draw (6.5,1.5) node{$x_p$};
				\draw[step=1cm,black,very thin] (0,2) grid (7,0);;
			\end{tikzpicture}&.
		\end{align*}
		The facets generated by these sets under the action of $S^0_p$ are all the sets $F=\{\varepsilon_0,\dots, \varepsilon_p\}$ such that $\varepsilon_i\in \{x_i,y_i\}$ for each $0\leq i \leq p$ and the number of partition vertices in $F$ is less than or equal to $m$.
		
		For example, the figure below shows the partition complexes $\calP(2,2,1)$ (on the left) and $\calP(2,2,2)$ (on the right).
		
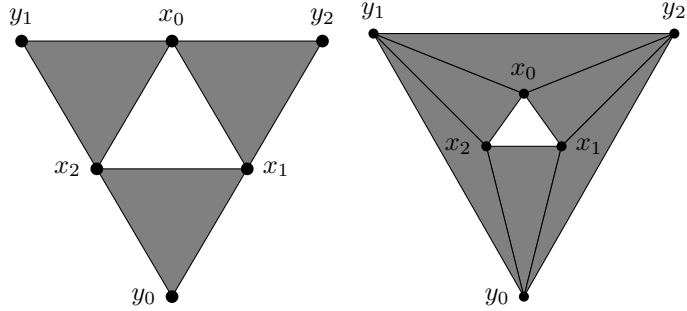
\begin{center}
			\begin{tikzpicture}[scale = 1]
				\tikzstyle{point}=[circle,thick,draw=black,fill=black,inner sep=0pt,minimum width=4pt,minimum height=4pt]
				\node (a)[point,label=above:$y_1$] at (0,3.4) {};
				\node (b)[point,label=above:$x_0$] at (2,3.4) {};
				\node (c)[point,label=above:$y_2$] at (4,3.4) {};
				\node (d)[point,label=left:$x_2$] at (1,1.7) {};
				\node (e)[point,label=right:$x_1$] at (3,1.7) {};
				\node (f)[point,label=left:$y_0$] at (2,0) {};	
				
				\begin{scope}[on background layer]
					\draw[fill=gray] (a.center) -- (b.center) -- (d.center) -- cycle;
					\draw[fill=gray] (b.center) -- (c.center) -- (e.center) -- cycle;
					\draw[fill=gray]   (d.center) -- (e.center) -- (f.center) -- cycle;
				\end{scope}
			\end{tikzpicture}
			\begin{tikzpicture}[scale = 1]
				\tikzstyle{point}=[circle,thick,draw=black,fill=black,inner sep=0pt,minimum width=3pt,minimum height=3pt]
				\node (a)[point,label=above:$y_1$] at (0,3.5) {};
				\node (b)[point,label=above:$x_0$] at (2,2.7) {};
				\node (c)[point,label=above:$y_2$] at (4,3.5) {};
				\node (d)[point,label=left:$x_2$] at (1.5,2) {};
				\node (e)[point,label=right:$x_1$] at (2.5,2) {};
				\node (f)[point,label=left:$y_0$] at (2,0) {};	
				
				\begin{scope}[on background layer]
					\draw[fill=gray] (a.center) -- (b.center) -- (d.center) -- cycle;
					\draw[fill=gray] (b.center) -- (c.center) -- (e.center) -- cycle;
					\draw[fill=gray]   (d.center) -- (e.center) -- (f.center) -- cycle;
					
					\draw[fill=gray] (a.center) -- (b.center) -- (c.center) -- cycle;
					\draw[fill=gray] (a.center) -- (f.center) -- (d.center) -- cycle;
					\draw[fill=gray] (c.center) -- (f.center) -- (e.center) -- cycle;
				\end{scope}
			\end{tikzpicture}
		\end{center}
	\end{ex}

	\begin{ex}\label{Example: Partition Complexes and Intersection Complexes}
		For a sequence $\bm = (\underbrace{a-1,0,\dots,0,1}_{p},0)$, the intersection complex $\calI(\bm)$ is equal to the partition complex $\calP(a,p,1)$. 
	\end{ex}
	
	%
	\begin{ex}\label{Example: Partition Complex m=p+1}
		
		\underline{\textbf{The case $m=p+1$:}}
		
		As mentioned at the start of this section, we are not interested in this case for its own sake; we include it only because complexes of this form occur as links in other partition complexes.
		
		Unlike the cases where $1\leq m \leq p$ we will see in the next section that the complex $\calP(a,p,p+1)$ is always acyclic. It may be helpful to think of $\calP(a,p,p+1)$ as the complex $\calP(a,p,p)$ but with some additional facets which `fill in' the homology.
		
		For example, as we saw in Example \ref{Example: Partition Complex p=m=1}, the complex $\calP(3,1,1)$ consists of two disjoint $2$-simplices: specifically, the faces $G_{2,1}^{(1)}$ and $gG_{2,1}^{(1)}$, where $g$ denotes the permutation $(0,1)$ in $S^0_1$.
		$$\begin{tikzpicture}[scale = 0.5]
			\tikzstyle{point}=[circle,thick,draw=black,fill=black,inner sep=0pt,minimum width=3pt,minimum height=3pt]
			\node (a)[point,label=above:$x_1$] at (-1.5,3) {};
			\node (b)[point,label=below:$y^2_0$] at (0,0) {};
			\node (c)[point,label=above:$y^1_0$] at (1.5,3) {};
			\node (d)[point,label=below:$y^1_1$] at (3,0) {};
			\node (e)[point,label=above:$y^2_1$] at (4.5,3) {};
			\node (f)[point,label=below:$x_0$] at (6,0) {};
			
			\node (g)[label= {$G_{2,1}^{(1)}$}] at (0,1.2) {};
			\node (j)[label= {$gG_{2,1}^{(1)}$}] at (4.5,0.04) {};
			
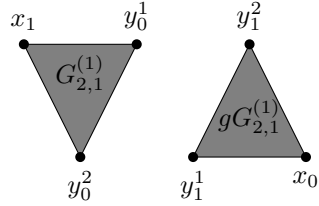
			\begin{scope}[on background layer]
				\draw[fill=gray] (a.center) -- (b.center) -- (c.center) -- cycle;
				\draw[fill=gray]   (d.center) -- (e.center) -- (f.center) -- cycle;
			\end{scope}
		\end{tikzpicture}$$
		Meanwhile, the complex $\calP(3,1,2)$ has two additional facets containing only partition vertices, which `bridge the gap' between these two disjoint simplices, thus making the resulting complex acyclic.
		$$
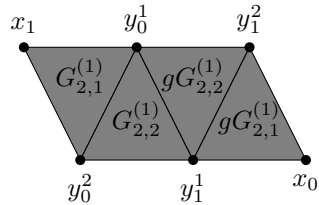
\begin{tikzpicture}[scale = 0.5]
			\tikzstyle{point}=[circle,thick,draw=black,fill=black,inner sep=0pt,minimum width=3pt,minimum height=3pt]
			\node (a)[point,label=above:$x_1$] at (-1.5,3) {};
			\node (b)[point,label=below:$y^2_0$] at (0,0) {};
			\node (c)[point,label=above:$y^1_0$] at (1.5,3) {};
			\node (d)[point,label=below:$y^1_1$] at (3,0) {};
			\node (e)[point,label=above:$y^2_1$] at (4.5,3) {};
			\node (f)[point,label=below:$x_0$] at (6,0) {};
			
			\node (g)[label= {$G_{2,1}^{(1)}$}] at (0,1.2) {};
			\node (h)[label= {$G_{2,2}^{(1)}$}] at (1.5,0.04) {};
			\node (i)[label= {$gG_{2,2}^{(1)}$}] at (3,1.2) {};
			\node (j)[label= {$gG_{2,1}^{(1)}$}] at (4.5,0.04) {};
			
			\begin{scope}[on background layer]
				\draw[fill=gray] (a.center) -- (b.center) -- (c.center) -- cycle;
				\draw[fill=gray] (b.center) -- (c.center) -- (d.center) -- cycle;
				\draw[fill=gray] (c.center) -- (d.center) -- (e.center) -- cycle;
				\draw[fill=gray]   (d.center) -- (e.center) -- (f.center) -- cycle;
			\end{scope}
		\end{tikzpicture}$$
		
		Similarly, the complex $\calP(2,2,2)$ can be drawn as
		$$
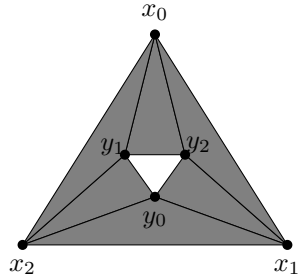
\begin{tikzpicture}[scale = 0.8]
			\tikzstyle{point}=[circle,thick,draw=black,fill=black,inner sep=0pt,minimum width=3pt,minimum height=3pt]
			\node (a)[point,label=below:$x_2$] at (-.2,0) {};
			\node (b)[point,label={[label distance = 0mm] below:$y_0$}] at (2,0.8) {};
			\node (c)[point,label=below:$x_1$] at (4.2,0) {};
			\node (d)[point,label={[label distance = -2.5mm] above left:$y_1$}] at (1.5,1.5) {};
			\node (e)[point,label={[label distance = -2.2mm] above right:$y_2$}] at (2.5,1.5) {};
			\node (f)[point,label=above:$x_0$] at (2,3.5) {};	
			
			\begin{scope}[on background layer]
				\draw[fill=gray] (a.center) -- (b.center) -- (d.center) -- cycle;
				\draw[fill=gray] (b.center) -- (c.center) -- (e.center) -- cycle;
				\draw[fill=gray]   (d.center) -- (e.center) -- (f.center) -- cycle;
				
				\draw[fill=gray] (a.center) -- (b.center) -- (c.center) -- cycle;
				\draw[fill=gray] (a.center) -- (f.center) -- (d.center) -- cycle;
				\draw[fill=gray] (c.center) -- (f.center) -- (e.center) -- cycle;
				
			\end{scope}
		\end{tikzpicture} $$
		while the complex $\calP(2,2,3)$ contains the additional facet $\{y_1,y_2,y_3\}$ which `fills in' the hole at the centre, once again making the complex acyclic.
		$$
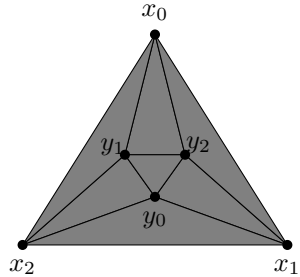
\begin{tikzpicture}[scale = 0.8]
			\tikzstyle{point}=[circle,thick,draw=black,fill=black,inner sep=0pt,minimum width=3pt,minimum height=3pt]
			\node (a)[point,label=below:$x_2$] at (-.2,0) {};
			\node (b)[point,label={[label distance = 0mm] below:$y_0$}] at (2,0.8) {};
			\node (c)[point,label=below:$x_1$] at (4.2,0) {};
			\node (d)[point,label={[label distance = -2.5mm] above left:$y_1$}] at (1.5,1.5) {};
			\node (e)[point,label={[label distance = -2.2mm] above right:$y_2$}] at (2.5,1.5) {};
			\node (f)[point,label=above:$x_0$] at (2,3.5) {};	
			
			\begin{scope}[on background layer]
				\draw[fill=gray] (a.center) -- (b.center) -- (d.center) -- cycle;
				\draw[fill=gray] (b.center) -- (c.center) -- (e.center) -- cycle;
				\draw[fill=gray]   (d.center) -- (e.center) -- (f.center) -- cycle;
				
				\draw[fill=gray] (a.center) -- (b.center) -- (c.center) -- cycle;
				\draw[fill=gray] (a.center) -- (f.center) -- (d.center) -- cycle;
				\draw[fill=gray] (c.center) -- (f.center) -- (e.center) -- cycle;
				
				\draw[fill=gray] (b.center) -- (d.center) -- (e.center) -- cycle;
			\end{scope}
		\end{tikzpicture}$$
		(Note we have adjusted our earlier picture of $\calP(2,2,2)$ in Example \ref{Example: Partition Complex a=2} here, by drawing the boundary vertices on the \textit{outside} of the complex; both depictions of this complex will be useful to us, and we frequently switch between them).
	\end{ex}
	
	The complex $\calP(a,p,p+1)$ can be obtained from $\calP(a,p,p)$ by adding faces consisting entirely of partition vertices. In a similar way, the following construction adds a face consisting entirely of boundary vertices.
	\begin{defn}\label{Definition: Closed Partition Complexes}
		Let $a$, $p$ and $m$ be integers with $a\geq 2$ and $0\leq m\leq p+1$. We define the \textit{closed partition complex} $\barP(a,p,m)$ on vertex set $V_p^a$ to be the complex $\calP(a,p,m)\cup \langle X_p\rangle$.
	\end{defn}
	\begin{rem}\label{Remark: Closed Partition Complex m=0 Cases}
		In the case $m=0$ the closed partition complex $\barP(a,p,0)$ is equal to the full simplex on the set $X_p$. Note that this is acyclic in all cases except for the case $p=-1$. In the case $p=-1$ the set $X_{-1}$ is empty so the full simplex on $X_{-1}$ is simply the irrelevant complex $\{\emptyset\}$.
	\end{rem}
	
	Just as with partition complexes of the form $\calP(a,p,p+1)$, we are interested in closed partition complexes only because they occur in the links of regular partition complexes, and our theorem on the degree types of partition complexes (Theorem \ref{Theorem: Partition Complex Degree Type}) does not extend to them. In fact, closed partition complexes are not even PR in general, except in the case where $a=2$. In all other cases the additional facet $X_p=\{x_0,\dots,x_p\}$ is of lower dimension than all the other facets of $\calP(a,p,m)$, and so $\barP(a,p,m)$ is not pure.
	
	We will see in the next section that the closure operation essentially acts as a kind of `switch' for the homology of a partition complex. Indeed, for every $1\leq m\leq p$ there is a deformation retraction from the partition complex $\calP(a,p,m)$ onto the boundary $\partial \langle X_p\rangle$, which means it has homology; whereas the closed partition complex $\barP(a,p,m)$ contains $X_p$ itself as a face, which means it deformation retracts on to the full simplex $\langle X_p\rangle$, and is hence acyclic. Meanwhile the partition complex $\calP(a,p,p+1)$ is acyclic, while its closure $\barP(a,p,p+1)$ has homology. 
	\begin{ex}\label{Example: Closed partition complex example}
		In the case $a=p=2$ we have the following partition complexes and closed partition complexes.
		\begin{center}
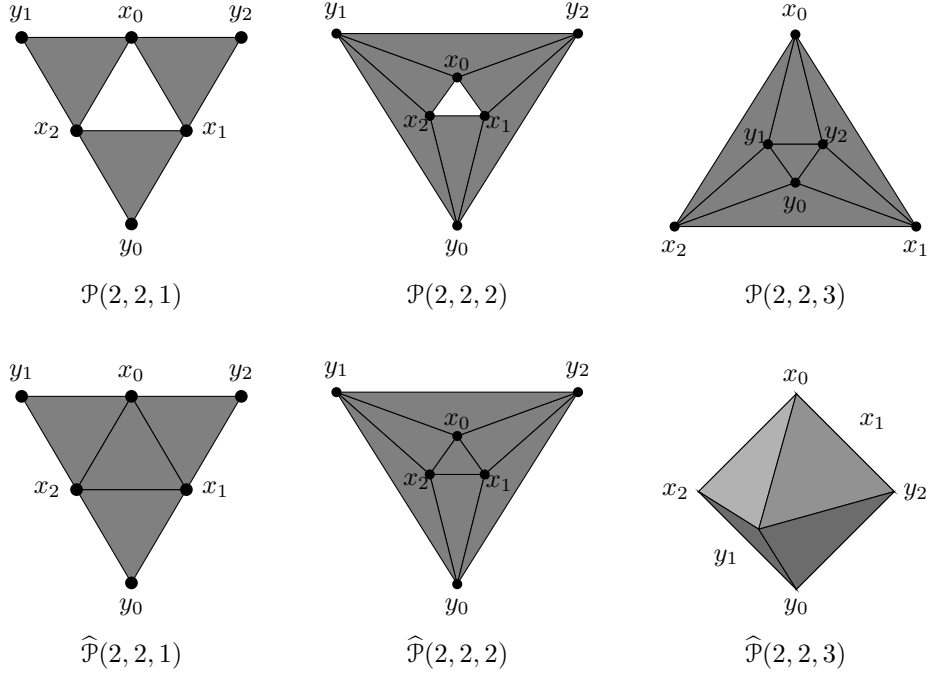
			\begin{tabular}{ c c c c c }
				\begin{tikzpicture}[scale = 0.73]
					\tikzstyle{point}=[circle,thick,draw=black,fill=black,inner sep=0pt,minimum width=4pt,minimum height=4pt]
					\node (a)[point,label=above:$y_1$] at (0,3.4) {};
					\node (b)[point,label=above:$x_0$] at (2,3.4) {};
					\node (c)[point,label=above:$y_2$] at (4,3.4) {};
					\node (d)[point,label=left:$x_2$] at (1,1.7) {};
					\node (e)[point,label=right:$x_1$] at (3,1.7) {};
					\node (f)[point,label=below:$y_0$] at (2,0) {};	
					
					\begin{scope}[on background layer]
						\draw[fill=gray] (a.center) -- (b.center) -- (d.center) -- cycle;
						\draw[fill=gray] (b.center) -- (c.center) -- (e.center) -- cycle;
						\draw[fill=gray]   (d.center) -- (e.center) -- (f.center) -- cycle;
					\end{scope}
				\end{tikzpicture}
				&& \begin{tikzpicture}[scale = 0.73]
					\tikzstyle{point}=[circle,thick,draw=black,fill=black,inner sep=0pt,minimum width=3pt,minimum height=3pt]
					\node (a)[point,label=above:$y_1$] at (-.2,3.5) {};
					\node (b)[point,label={[label distance = -1mm] above:$x_0$}] at (2,2.7) {};
					\node (c)[point,label=above:$y_2$] at (4.2,3.5) {};
					\node (d)[point,label={[label distance = -2.5mm] below left:$x_2$}] at (1.5,2) {};
					\node (e)[point,label={[label distance = -2.2mm] below right:$x_1$}] at (2.5,2) {};
					\node (f)[point,label=below:$y_0$] at (2,0) {};	
					
					\begin{scope}[on background layer]
						\draw[fill=gray] (a.center) -- (b.center) -- (d.center) -- cycle;
						\draw[fill=gray] (b.center) -- (c.center) -- (e.center) -- cycle;
						\draw[fill=gray]   (d.center) -- (e.center) -- (f.center) -- cycle;
						
						\draw[fill=gray] (a.center) -- (b.center) -- (c.center) -- cycle;
						\draw[fill=gray] (a.center) -- (f.center) -- (d.center) -- cycle;
						\draw[fill=gray] (c.center) -- (f.center) -- (e.center) -- cycle;
					\end{scope}
				\end{tikzpicture} && \begin{tikzpicture}[scale = 0.73]
					\tikzstyle{point}=[circle,thick,draw=black,fill=black,inner sep=0pt,minimum width=3pt,minimum height=3pt]
					\node (a)[point,label=below:$x_2$] at (-.2,0) {};
					\node (b)[point,label={[label distance = 0mm] below:$y_0$}] at (2,0.8) {};
					\node (c)[point,label=below:$x_1$] at (4.2,0) {};
					\node (d)[point,label={[label distance = -2.5mm] above left:$y_1$}] at (1.5,1.5) {};
					\node (e)[point,label={[label distance = -2.2mm] above right:$y_2$}] at (2.5,1.5) {};
					\node (f)[point,label=above:$x_0$] at (2,3.5) {};	
					
					\begin{scope}[on background layer]
						\draw[fill=gray] (a.center) -- (b.center) -- (d.center) -- cycle;
						\draw[fill=gray] (b.center) -- (c.center) -- (e.center) -- cycle;
						\draw[fill=gray]   (d.center) -- (e.center) -- (f.center) -- cycle;
						
						\draw[fill=gray] (a.center) -- (b.center) -- (c.center) -- cycle;
						\draw[fill=gray] (a.center) -- (f.center) -- (d.center) -- cycle;
						\draw[fill=gray] (c.center) -- (f.center) -- (e.center) -- cycle;
						
						\draw[fill=gray] (b.center) -- (d.center) -- (e.center) -- cycle;
					\end{scope}
				\end{tikzpicture} \\
				$\calP(2,2,1)$&& $\calP(2,2,2)$ && $\calP(2,2,3)$\\
				&& &&\\
				\begin{tikzpicture}[scale = 0.73]
					\tikzstyle{point}=[circle,thick,draw=black,fill=black,inner sep=0pt,minimum width=4pt,minimum height=4pt]
					\node (a)[point,label=above:$y_1$] at (0,3.4) {};
					\node (b)[point,label=above:$x_0$] at (2,3.4) {};
					\node (c)[point,label=above:$y_2$] at (4,3.4) {};
					\node (d)[point,label=left:$x_2$] at (1,1.7) {};
					\node (e)[point,label=right:$x_1$] at (3,1.7) {};
					\node (f)[point,label=below:$y_0$] at (2,0) {};	
					
					\begin{scope}[on background layer]
						\draw[fill=gray] (a.center) -- (b.center) -- (d.center) -- cycle;
						\draw[fill=gray] (b.center) -- (c.center) -- (e.center) -- cycle;
						\draw[fill=gray]   (d.center) -- (e.center) -- (f.center) -- cycle;
						
						\draw[fill=gray] (b.center) -- (d.center) -- (e.center) -- cycle;
					\end{scope}
				\end{tikzpicture}
				&& \begin{tikzpicture}[scale = 0.73]
					\tikzstyle{point}=[circle,thick,draw=black,fill=black,inner sep=0pt,minimum width=3pt,minimum height=3pt]
					\node (a)[point,label=above:$y_1$] at (-.2,3.5) {};
					\node (b)[point,label={[label distance = -1mm] above:$x_0$}] at (2,2.7) {};
					\node (c)[point,label=above:$y_2$] at (4.2,3.5) {};
					\node (d)[point,label={[label distance = -2.5mm] below left:$x_2$}] at (1.5,2) {};
					\node (e)[point,label={[label distance = -2.2mm] below right:$x_1$}] at (2.5,2) {};
					\node (f)[point,label=below:$y_0$] at (2,0) {};	
					
					\begin{scope}[on background layer]
						\draw[fill=gray] (a.center) -- (b.center) -- (d.center) -- cycle;
						\draw[fill=gray] (b.center) -- (c.center) -- (e.center) -- cycle;
						\draw[fill=gray]   (d.center) -- (e.center) -- (f.center) -- cycle;
						
						\draw[fill=gray] (a.center) -- (b.center) -- (c.center) -- cycle;
						\draw[fill=gray] (a.center) -- (f.center) -- (d.center) -- cycle;
						\draw[fill=gray] (c.center) -- (f.center) -- (e.center) -- cycle;
						
						\draw[fill=gray] (b.center) -- (d.center) -- (e.center) -- cycle;
					\end{scope}
				\end{tikzpicture} &&	\begin{tikzpicture}[scale=1.3][line join=bevel,z=-5.5]
					\coordinate [label={[label distance=3mm]above right:$x_1$}] (A1) at (0,0,-1);
					\coordinate [label=left:$x_2$] (A2) at (-1,0,0);
					\coordinate [label={[label distance=2mm]below left: $y_1$}] (A3) at (0,0,1);
					\coordinate [label=right:$y_2$] (A4) at (1,0,0);
					\coordinate [label=above:$x_0$] (B1) at (0,1,0);
					\coordinate [label=below:$y_0$] (C1) at (0,-1,0);
					
					\begin{scope}[on background layer]
						\draw (A1) -- (A2) -- (B1) -- cycle;
						\draw (A4) -- (A1) -- (B1) -- cycle;
						\draw (A1) -- (A2) -- (C1) -- cycle;
						\draw (A4) -- (A1) -- (C1) -- cycle;
						\draw [fill opacity=0.7,fill=lightgray!80!gray] (A2) -- (A3) -- (B1) -- cycle;
						\draw [fill opacity=0.7,fill=gray!70!lightgray] (A3) -- (A4) -- (B1) -- cycle;
						\draw [fill opacity=0.7,fill=darkgray!30!gray] (A2) -- (A3) -- (C1) -- cycle;
						\draw [fill opacity=0.7,fill=darkgray!30!gray] (A3) -- (A4) -- (C1) -- cycle;
					\end{scope}
				\end{tikzpicture}\\
				$\barP(2,2,1)$&& $\barP(2,2,2)$ && $\barP(2,2,3)$
			\end{tabular}
		\end{center}
	\end{ex}
	We now introduce some key notation and terminology which will help us to discuss the faces of partition complexes.
	
	\begin{notation}\label{Notation: sigma-X and sigma-Y and support}
		Let $a$, $p$ and $m$ be integers with $a\geq 2$ and $0\leq m \leq p+1$, and let $\sigma$ be a subset of $V^a_p$. We use
		\begin{enumerate}
			\item $\sigma_X$ to denote the intersection $\sigma \cap X_p$.
			\item $\sigma_Y$ to denote the intersection $\sigma \cap Y^a_p$.
		\end{enumerate}
		We also define the \textit{support} of $\sigma$ to be the set $$\Supp(\sigma) = \left \{i\in \{0,\dots, p\} : x_i \in \sigma \text{ or } y_i^j \in \sigma \text{ for some } 1\leq j\leq a-1\right \}.$$
	\end{notation}
	
	\begin{defn}\label{Definition: Totally Separated and partition complete}
		Let $a$, $p$ and $m$ be integers with $a\geq 2$ and $0\leq m \leq p+1$, and let $F$ be a subset of the vertex set $V^a_p$. We say
		\begin{enumerate}
			\item $F$ is \textit{partition complete} if for every partition vertex $y_i^j$ in $\sigma$, the vertices $y_i^1,\dots,y_i^{j-1}$ are also in $\sigma$.
			\item $F$ is \textit{separated}  if for each $0\leq i \leq p$, $F$ contains at most one of the vertices $x_i$ or $y_i$. We say it is \textit{totally separated} if for each $0\leq i \leq p$ it contains exactly one of the vertices $x_i$ or $y_i$. 
		\end{enumerate}
	\end{defn}
	\begin{ex}
		The following are examples of subsets of the vertex set $V^4_3$.
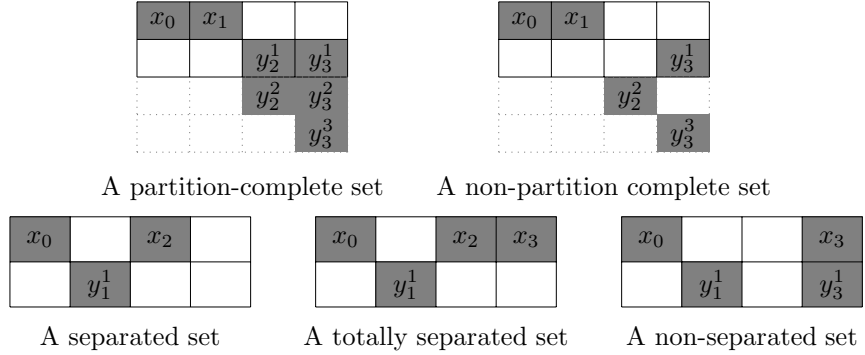
		\begin{center}
			\begin{tabular}{ c c c }
				\begin{tikzpicture}[xscale=.7,yscale=.5]
					\fill[gray] (0,4) rectangle (1,3);
					\draw (0.5,3.5) node{$x_0$};
					\fill[gray] (1,4) rectangle (2,3);
					\draw (1.5,3.5) node{$x_1$};
					\fill[gray] (2,3) rectangle (3,2);
					\draw (2.5,2.5) node{$y^1_2$};
					\fill[gray] (3,3) rectangle (4,2);
					\draw (3.5,2.5) node{$y^1_3$};
					\fill[gray] (2,2) rectangle (3,1);
					\draw (2.5,1.5) node{$y^2_2$};
					\fill[gray] (3,2) rectangle (4,1);
					\draw (3.5,1.5) node{$y^2_3$};
					\fill[gray] (3,1) rectangle (4,0);
					\draw (3.5,0.5) node{$y^3_3$};
					\draw[step=1cm,black,very thin] (0,4) grid (4,2);
					\draw[step=1cm,gray,dotted] (0,2) grid (4,0);
				\end{tikzpicture} && \begin{tikzpicture}[xscale=.7,yscale=.5]
					\fill[gray] (0,4) rectangle (1,3);
					\draw (0.5,3.5) node{$x_0$};
					\fill[gray] (1,4) rectangle (2,3);
					\draw (1.5,3.5) node{$x_1$};
					\fill[gray] (3,3) rectangle (4,2);
					\draw (3.5,2.5) node{$y^1_3$};
					\fill[gray] (2,2) rectangle (3,1);
					\draw (2.5,1.5) node{$y^2_2$};
					\fill[gray] (3,1) rectangle (4,0);
					\draw (3.5,0.5) node{$y^3_3$};
					\draw[step=1cm,black,very thin] (0,4) grid (4,2);
					\draw[step=1cm,gray,dotted] (0,2) grid (4,0);
				\end{tikzpicture}\\
				A partition-complete set && A non-partition complete set
			\end{tabular}
		\end{center}
		\begin{center}
			\begin{tabular}{ c c c c c }
				\begin{tikzpicture}[xscale=.8,yscale=.6]
					\fill[gray] (0,4) rectangle (1,3);
					\draw (0.5,3.5) node{$x_0$};
					\fill[gray] (1,3) rectangle (2,2);
					\draw (1.5,2.5) node{$y^1_1$};
					\fill[gray] (2,4) rectangle (3,3);
					\draw (2.5,3.5) node{$x_2$};
					\draw[step=1cm,black,very thin] (0,4) grid (4,2);
				\end{tikzpicture}&& \begin{tikzpicture}[xscale=.8,yscale=.6]
					\fill[gray] (0,4) rectangle (1,3);
					\draw (0.5,3.5) node{$x_0$};
					\fill[gray] (1,3) rectangle (2,2);
					\draw (1.5,2.5) node{$y^1_1$};
					\fill[gray] (2,4) rectangle (3,3);
					\draw (2.5,3.5) node{$x_2$};
					\fill[gray] (3,4) rectangle (4,3);
					\draw (3.5,3.5) node{$x_3$};
					\draw[step=1cm,black,very thin] (0,4) grid (4,2);
				\end{tikzpicture} && \begin{tikzpicture}[xscale=.8,yscale=.6]
					\fill[gray] (0,4) rectangle (1,3);
					\draw (0.5,3.5) node{$x_0$};
					\fill[gray] (1,3) rectangle (2,2);
					\draw (1.5,2.5) node{$y^1_1$};
					\fill[gray] (3,4) rectangle (4,3);
					\draw (3.5,3.5) node{$x_3$};
					\fill[gray] (3,3) rectangle (4,2);
					\draw (3.5,2.5) node{$y^1_3$};
					\draw[step=1cm,black,very thin] (0,4) grid (4,2);
				\end{tikzpicture}\\
				A separated set && A totally separated set && A non-separated set
			\end{tabular}
		\end{center}
	\end{ex}
	
	By construction every face of a partition complex is separated, and every facet is both totally separated and partition complete. In fact we have more than this.
	\begin{lem}\label{Lemma: Partition Complex Description of Facets}
		Let $a$, $p$ and $m$ be integers with $a\geq 2$ and $0\leq m \leq p+1$. A subset $F\subseteq V^a_p$ is a facet of the partition complex $\calP(a,p,m)$ if and only if it satisfies the following four conditions.
		\begin{enumerate}
			\item $F$ is partition complete.
			\item $F$ is totally separated.
			\item $|F|=a+p-1$.
			\item $1\leq |\Supp(F_Y)|\leq m$.
		\end{enumerate}
	\end{lem}
	\begin{proof}
		Every facet of $\calP(a,p,m)$ must satisfy these four conditions, because the generating facets $G_{p,i}^\lambda$ satisfy them by construction, and all four conditions are invariant under the action of $S^0_p$.
		
		Conversely, suppose $F$ is a totally separated, partition complete subset of $V^a_p$ of size $a+p-1$ with $|\Supp(F_Y)|=i$ for some $1\leq i \leq m$. We can permute the vertices of $F$ to ensure that $\Supp(F_Y)=\{0,\dots,i-1\}$. Because $F$ is totally separated we must therefore have $F_X=\{x_i,\dots,x_p\}=G_{p,i}$.
		
		This leaves us with a total of $(a+p-1)-(p+1-i)=a+i-2$ partition vertices in $F$. For each $0\leq r \leq i-1$ we let $j_r$ denote the maximum integer between $1$ and $a-1$ such that the partition vertex $y_r^{j_r}$ is in $F$. Once again, we can permute the vertices of $F$ to ensure that the sequence $\lambda = (j_0,\dots,j_{i-1})$ is monotonically decreasing. The sequence $\lambda$ is a partition of $a+i-2$ into $i$ parts, and because $F$ is partition complete we have (after our permutations of vertices) that $F$ is equal to the generating facet $G_{p,i}^\lambda$.
	\end{proof}

	Now that we have demonstrated the construction of partition complexes, and built up the tools we will need to talk about them, we present our main result.
	\begin{thm}\label{Theorem: Partition Complex Degree Type}
		Let $a$, $p$ and $m$ be positive integers with $a\geq 2$ and $1 \leq m\leq p$. The partition complex $\calP(a,p,m)$ is PR with degree type $(\overbrace{1,\dots,1\underbrace{a,1,\dots,1}_{m}}^{p})$.
	\end{thm}
	
	In the next section we will make use of some deformation retractions to find the homology of both partition complexes and their closures. We also show how all of the links in partition complexes with homology can be built out of smaller partition complexes and closed partition complexes. Then in Section \ref{Subsection: Proving Partition Complex Theorem} we will put these results together to prove Theorem \ref{Theorem: Partition Complex Degree Type}.
	
	\subsection{Deformation Retractions and Links}
	Just as with intersection complexes, we can find the homology of the links in partition complexes by using some deformation retractions. Once again, our main tools for this are Lemma \ref{Lemma: Deformation Retract} and Corollary \ref{Corollary: Deformation Retract}.
	
	\begin{prop}\label{Proposition: Deformation D(a,p,m) to D(2,p,m)}
		Let $a$, $p$ and $m$ be integers with $a\geq 2$ and $1\leq m\leq p+1$. There is a deformation retraction $\calP(a,p,m)\rightsquigarrow \calP(2,p,m)$ given by the vertex maps $y_i^j\mapsto y^1_i$. There is a similar deformation retraction $\barP(a,p,m)\rightsquigarrow \barP(2,p,m)$.
	\end{prop}
	\begin{proof}
		Let $y_i^j$ be a partition vertex of $\calP(a,p,m)$ for some $j\geq 2$. By construction, every facet of $\calP(a,p,m)$ which contains the vertex $y_i^j$ also contains the vertex $y_i^1$. Thus Corollary \ref{Corollary: Deformation Retract} gives us a deformation retraction of $\calP(a,p,m)$ on to the complex obtained by deleting every partition vertex $y_i^j$ for $j\geq 2$, by 
		collapsing each of these vertices onto the corresponding vertex $y^1_i$.
		The complex obtained from these deletions is $\calP(2,p,m)$. The proof for $\barP(a,p,m)$ is identical.
	\end{proof}
	
	\begin{prop}\label{Proposition: Deformation D(2,p,m) to D(1,p,m)}
		Let $p$ and $m$ be positive integers with $1\leq m\leq p$. There is a deformation retraction $\calP(2,p,m)\rightsquigarrow \partial \langle X_p \rangle$. There is a similar deformation retraction $\barP(2,p,m)\rightsquigarrow \langle X_p \rangle$.
	\end{prop}
	\begin{proof}
		As noted in Example \ref{Example: Partition Complex a=2}, all the facets of $\calP(2,p,m)$ contain no more than $m$ partition vertices. Thus if $F$ is any facet with exactly $m$ partition vertices $y_{i_1},\dots,y_{i_m}$, it must be the \textit{only} facet containing the face $\{y_{i_1},\dots,y_{i_m}\}$. Because $m\leq p$, we know $F$ \textit{strictly} contains $\{y_{i_1},\dots,y_{i_m}\}$, and hence Lemma \ref{Lemma: Deformation Retract} allows us to retract the face $\{y_{i_1},\dots,y_{i_m}\}$. The complex obtained by deleting \textit{all} faces consisting of $m$ partition vertices from $\calP(2,p,m)$ is $\calP(2,p,m-1)$.
		
		By induction on $m\geq 1$ we obtain a series of deformation retractions $\calP(2,p,m)\rightsquigarrow \calP(2,p,m-1) \rightsquigarrow \dots \rightsquigarrow \calP(2,p,1)$. The facets of $\calP(2,p,1)$ are all the totally separated subsets of $V^2_p$ with only a single partition vertex, and so once again we may use Lemma \ref{Lemma: Deformation Retract} to delete these partition vertices. This leaves us with a complex with facets of the form $\{x_0,\dots,x_p\}-\{x_i\}$ for $0\leq i \leq p$. This is $\partial \langle X_p\rangle$.
		
		The proof for $\barP(2,p,m)$ is identical, except because $\barP(2,p,m)$ also contains the face $X_p$, it deformation retracts onto the full simplex $\langle X_p\rangle$.
		
		
	\end{proof}
	\begin{ex}\label{Example: D(2,2,2) def retract}
		In the case $p=m=2$ we have the following deformation retraction.
		\begin{center}
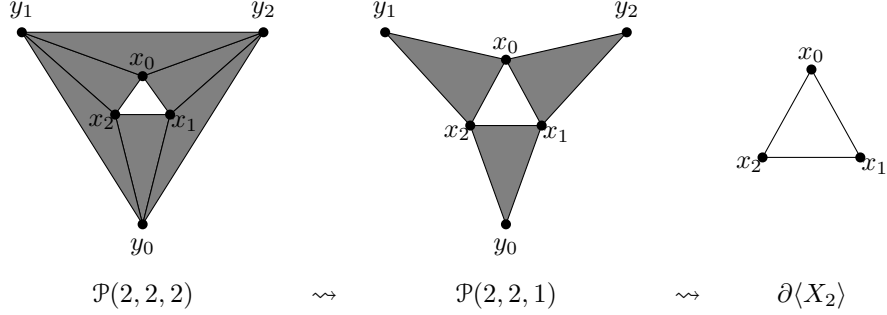
			\begin{tabular}{ c c c c c }
				\begin{tikzpicture}[scale = 0.73]
					\tikzstyle{point}=[circle,thick,draw=black,fill=black,inner sep=0pt,minimum width=3pt,minimum height=3pt]
					\node (a)[point,label=above:$y_1$] at (-.2,3.5) {};
					\node (b)[point,label={[label distance = -1mm] above:$x_0$}] at (2,2.7) {};
					\node (c)[point,label=above:$y_2$] at (4.2,3.5) {};
					\node (d)[point,label={[label distance = -2.5mm] below left:$x_2$}] at (1.5,2) {};
					\node (e)[point,label={[label distance = -2.2mm] below right:$x_1$}] at (2.5,2) {};
					\node (f)[point,label=below:$y_0$] at (2,0) {};	
					
					\begin{scope}[on background layer]
						\draw[fill=gray] (a.center) -- (b.center) -- (d.center) -- cycle;
						\draw[fill=gray] (b.center) -- (c.center) -- (e.center) -- cycle;
						\draw[fill=gray]   (d.center) -- (e.center) -- (f.center) -- cycle;
						
						\draw[fill=gray] (a.center) -- (b.center) -- (c.center) -- cycle;
						\draw[fill=gray] (a.center) -- (f.center) -- (d.center) -- cycle;
						\draw[fill=gray] (c.center) -- (f.center) -- (e.center) -- cycle;
					\end{scope}
				\end{tikzpicture} & 
				&  \begin{tikzpicture}[scale = 0.73]
					\tikzstyle{point}=[circle,thick,draw=black,fill=black,inner sep=0pt,minimum width=3pt,minimum height=3pt]
					\node (a)[point,label=above:$y_1$] at (-.2,3.5) {};
					\node (b)[point,label={[label distance = -1mm] above:$x_0$}] at (2,3) {};
					\node (c)[point,label=above:$y_2$] at (4.2,3.5) {};
					\node (d)[point,label={[label distance = -2.5mm] below left:$x_2$}] at (1.35,1.8) {};
					\node (e)[point,label={[label distance = -2.2mm] below right:$x_1$}] at (2.65,1.8) {};
					\node (f)[point,label=below:$y_0$] at (2,0) {};	
					
					\begin{scope}[on background layer]
						\draw[fill=gray] (a.center) -- (b.center) -- (d.center) -- cycle;
						\draw[fill=gray] (b.center) -- (c.center) -- (e.center) -- cycle;
						\draw[fill=gray]   (d.center) -- (e.center) -- (f.center) -- cycle;
						
					\end{scope}
				\end{tikzpicture} & 
				& \begin{tikzpicture}[scale = 1.3]
					\tikzstyle{point}=[circle,thick,draw=black,fill=black,inner sep=0pt,minimum width=3pt,minimum height=3pt]
					\node (b)[point,label={[label distance = -1mm] above:$x_0$}] at (2,2.9) {};
					\node (d)[point,label={[label distance = -2.5mm] below left:$x_2$}] at (1.5,2) {};
					\node (e)[point,label={[label distance = -2.2mm] below right:$x_1$}] at (2.5,2) {};
					\node at (2,1) {};
					
					\draw (b.center) -- (d.center) -- (e.center) -- cycle;
					
					
				\end{tikzpicture}\\
				$\calP(2,2,2)$& $\rightsquigarrow$  & $\calP(2,2,1)$ & $\rightsquigarrow$  & $\partial \langle X_2 \rangle$
			\end{tabular}
		\end{center}	
	\end{ex}
	
	\begin{prop}\label{Proposition: Deformation D(2,p,p+1) to point}
		For any integer $p\geq 0$, the partition complex $\calP(2,p,p+1)$ is acyclic.
	\end{prop}
	\begin{proof}
		Our proof for this result is very similar to our proof for Proposition \ref{Proposition: Deformation D(2,p,m) to D(1,p,m)} above. In that proof, we used Lemma \ref{Lemma: Deformation Retract} to remove all of the partition vertices $y_0,\dots,y_p$ from the complex. In this case, we use the same lemma to remove the boundary vertices $x_0,\dots,x_p$.
		
		Note that while $\Delta=\calP(2,p,p+1)$ contains the facet $\{y_0,\dots,y_p\}$ (which has $p+1$ partition vertices), no facet of $\Delta$ contains more than $p$ boundary vertices. Thus if $F$ is a facet of $\Delta$ containing $p$ boundary vertices $x_{i_1},\dots,x_{i_p}$, it must be the \textit{only} facet containing all of those vertices; and once again, because it \textit{strictly} contains them, Lemma \ref{Lemma: Deformation Retract} allows us to retract the face $\{x_{i_1},\dots,x_{i_p}\}$. Proceeding in this way, we may remove all faces consisting of boundary vertices from $\Delta$, beginning with the ones of size $p$ and continuing in decreasing order of size. Thus we obtain a deformation retraction of $\calP(2,p,p+1)$ on to the complex with a single facet $\{y_0,\dots,y_p\}$, which is acyclic.
	\end{proof}
	\begin{ex}\label{Example: D(2,2,3) def retract}
		In the case $p=2$, $m=3$ we have the following deformation retraction.
		\begin{center}
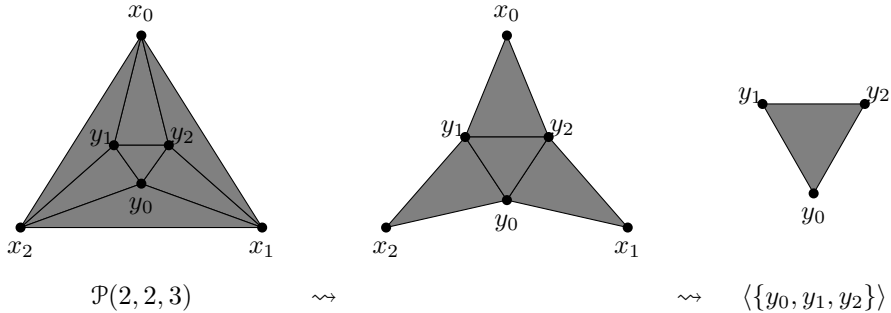
			\begin{tabular}{ c c c c c }
				\begin{tikzpicture}[scale = 0.73]
					\tikzstyle{point}=[circle,thick,draw=black,fill=black,inner sep=0pt,minimum width=3pt,minimum height=3pt]
					\node (a)[point,label=below:$x_2$] at (-.2,0) {};
					\node (b)[point,label={[label distance = 0mm] below:$y_0$}] at (2,0.8) {};
					\node (c)[point,label=below:$x_1$] at (4.2,0) {};
					\node (d)[point,label={[label distance = -2.5mm] above left:$y_1$}] at (1.5,1.5) {};
					\node (e)[point,label={[label distance = -2.2mm] above right:$y_2$}] at (2.5,1.5) {};
					\node (f)[point,label=above:$x_0$] at (2,3.5) {};	
					
					\begin{scope}[on background layer]
						\draw[fill=gray] (a.center) -- (b.center) -- (d.center) -- cycle;
						\draw[fill=gray] (b.center) -- (c.center) -- (e.center) -- cycle;
						\draw[fill=gray]   (d.center) -- (e.center) -- (f.center) -- cycle;
						
						\draw[fill=gray] (a.center) -- (b.center) -- (c.center) -- cycle;
						\draw[fill=gray] (a.center) -- (f.center) -- (d.center) -- cycle;
						\draw[fill=gray] (c.center) -- (f.center) -- (e.center) -- cycle;
						
						\draw[fill=gray] (b.center) -- (d.center) -- (e.center) -- cycle;
					\end{scope}
				\end{tikzpicture} & 
				&  \begin{tikzpicture}[scale = 0.73]
					\tikzstyle{point}=[circle,thick,draw=black,fill=black,inner sep=0pt,minimum width=3pt,minimum height=3pt]
					\node (a)[point,label=below:$x_2$] at (-.2,0) {};
					\node (b)[point,label={[label distance = 0mm] below:$y_0$}] at (2,0.5) {};
					\node (c)[point,label=below:$x_1$] at (4.2,0) {};
					\node (d)[point,label={[label distance = -2.5mm] above left:$y_1$}] at (1.24,1.65) {};
					\node (e)[point,label={[label distance = -2.2mm] above right:$y_2$}] at (2.76,1.65) {};
					\node (f)[point,label=above:$x_0$] at (2,3.5) {};	
					
					\begin{scope}[on background layer]
						\draw[fill=gray] (a.center) -- (b.center) -- (d.center) -- cycle;
						\draw[fill=gray] (b.center) -- (c.center) -- (e.center) -- cycle;
						\draw[fill=gray]   (d.center) -- (e.center) -- (f.center) -- cycle;
						
						
						\draw[fill=gray] (b.center) -- (d.center) -- (e.center) -- cycle;
					\end{scope}
				\end{tikzpicture} & 
				& \begin{tikzpicture}[scale = 1.7]
					\tikzstyle{point}=[circle,thick,draw=black,fill=black,inner sep=0pt,minimum width=3pt,minimum height=3pt]
					\node (b)[point,label={[label distance = 0mm] below:$y_0$}] at (2,0.8) {};
					\node (d)[point,label={[label distance = -2.5mm] above left:$y_1$}] at (1.6,1.5) {};
					\node (e)[point,label={[label distance = -2.2mm] above right:$y_2$}] at (2.4,1.5) {};
					\node at (2,0.3) {};
					
					\begin{scope}[on background layer]
						
						
						\draw[fill=gray] (b.center) -- (d.center) -- (e.center) -- cycle;
					\end{scope}
				\end{tikzpicture}\\
				$\calP(2,2,3)$& $\rightsquigarrow$  &  & $\rightsquigarrow$  & $\langle \{y_0,y_1,y_2\} \rangle$
			\end{tabular}
		\end{center}	
	\end{ex}

	\begin{prop}\label{Proposition: Partition Complex Cross Polytope}
		Let $p\geq -1$ be an integer. The closed partition complex $\barP(2,p,p+1)$ has only $\nth[th]{p}$ homology, of dimension $1$.
	\end{prop}
	\begin{proof}
		We proceed by induction on $p\geq -1$. The base case is immediate because $\barP(2,-1,0)$ is the complex $\{\emptyset\}$, which has only $\nth[st]{(-1)}$ homology, of dimension $1$.
		
		For the inductive step, suppose $p\geq 0$ and set $\barDelta=\barP(2,p,p+1)$. Because every facet of $\barDelta$ is totally separated, we can decompose $\barDelta$ into the subcomplexes $A$ and $B$ generated, respectively, by facets containing $x_p$ and facets containing $y_p$. Note that $A$ (resp. $B$) is a cone over the vertex $x_p$ (resp. $y_p$) and is therefore acyclic. The intersection $A\cap B$ is the set of all separated subsets of $V^2_{p-1}=\{x_0,\dots,x_{p-1},y_0,\dots,y_{p-1}\}$, and is thus equal to $\barP(2,p-1,p)$. Thus, for each $i\geq 0$, the Mayer-Vietoris Sequence gives us an isomorphism
		$$0\rightarrow \Hred_i(\barDelta)\rightarrow \Hred_{i-1}(\barP(2,p-1,p))\rightarrow 0$$
		and the result follows by the inductive hypothesis.

	\end{proof}
	
	Taken together these results are enough to give us the homology of all partition complexes and closed partition complexes.
	\begin{cor}\label{Corollary: Homology of Partition Complexes}
		Let $a$, $p$ and $m$ be integers with $a\geq 2$, $p\geq -1$ and $0\leq m\leq p+1$. Let $\Delta = \calP(a,p,m)$ and $\barDelta = \barP(a,p,m)$.
		\begin{enumerate}
			\item $h(\Delta) = \begin{cases}
				\{-1\} &\text{ if } m=0\\
				\{p-1\} &\text{ if } 1\leq m\leq p\\
				\emptyset & \text{ if } m=p+1 \text{ and } p\neq -1
			\end{cases}$
			
			and in the first two cases, the dimension of the nontrivial homology is $1$.
			\item $h(\barDelta) = \begin{cases}
				\emptyset &\text{ if } m=0 \text{ and } p\neq -1\\
				\emptyset & \text{ if } 1\leq m\leq p\\
				\{p\} &\text{ if } m=p+1
			\end{cases}$
			
			and in the final case, the dimension of the nontrivial homology is $1$.
		\end{enumerate}
	\end{cor}
	\begin{proof}
		We consider the $m=0$ cases together first. The partition complex $\calP(a,p,0)$ is defined to be the irrelevant complex $\{\emptyset\}$ which has only $\nth[st]{(-1)}$ homology. Meanwhile the closed partition complex $\barP(a,p,0)$ is the full $p$-simplex $\langle X_p\rangle$, which is acyclic except in the case where $p=-1$. We now proceed to the other cases.
		
		For part (1), if $1\leq m\leq p$, then Propositions \ref{Proposition: Deformation D(a,p,m) to D(2,p,m)} and \ref{Proposition: Deformation D(2,p,m) to D(1,p,m)} give us a deformation retraction from $\calP(a,p,m)$ on to $\partial \langle X_p \rangle$, which is the boundary of a $p$-simplex, and hence has only $\nth[st]{(p-1)}$ homology of dimension 1. If $m=p+1$ and $p\neq -1$, then Proposition \ref{Proposition: Deformation D(a,p,m) to D(2,p,m)} gives us a deformation retraction from $\calP(a,p,p+1)$ on to $\calP(2,p,p+1)$, which is acyclic by Proposition \ref{Proposition: Deformation D(2,p,p+1) to point}.
		
		For part (2), if $1\leq m\leq p$, then Propositions \ref{Proposition: Deformation D(a,p,m) to D(2,p,m)} and \ref{Proposition: Deformation D(2,p,m) to D(1,p,m)} give us a deformation retraction from $\barP(a,p,m)$ to $\langle X_p \rangle$, which is a full $p$-simplex, and is hence acyclic. If $m=p+1$, then Proposition \ref{Proposition: Deformation D(a,p,m) to D(2,p,m)} gives us a deformation retraction from $\barP(a,p,m)$ on to $\barP(2,p,m)$. By Proposition \ref{Proposition: Partition Complex Cross Polytope} this has only $\nth[th]{p}$ homology, of dimension $1$. 
	\end{proof}
	
	We now examine the links in partition complexes. Suppose $\sigma$ is a face of a partition complex $\calP(a,p,m)$. We can partition $\sigma$ into $\sigma_X\sqcup \sigma_Y$ (as defined in Notation \ref{Notation: sigma-X and sigma-Y and support}). We first consider the case $\sigma_Y = \emptyset$. In fact, it suffices to consider the subcase where $\sigma$ is equal to the boundary vertex $x_p$.
	
	\begin{lem}\label{Lemma: Partition Complex Link of x}
		Let $a$, $p$ and $m$ be integers with $a\geq 2$ and $1\leq m\leq p+1$, and let $\Delta=\calP(a,p,m)$. We have
		$$\link_\Delta x_p = \begin{cases}
			\calP(a,p-1,m)& \text{ if } 1\leq m \leq p\\
			\calP(a,p-1,p)& \text{ if } m = p+1.
		\end{cases}$$
	\end{lem}
	\begin{proof}
		We begin with the case where $1\leq m \leq p$. Let $F$ be a facet of $\Delta$ containing $x_p$. We know $F$ is of the form $gG_{p,i}^\lambda$ for some permutation $g\in S^0_p$, some integer $1\leq i \leq m$ and some partition $\lambda \in \PP(a,i)$.

			Let $i=g^{-1}(p)$ and define $\tg$ to be the permutation in $S^0_p$ given by
			$$\tg(j)=\begin{cases}
				p &\text{ if } j=p\\
				g(p) &\text{ if } j=i\\
				g(j) &\text{ otherwise.}
			\end{cases}$$ The facet $gG^\lambda_{p,i}$ contains the vertex $x_p=x_{\tg(p)}$ by assumption, and the vertex  $g(x_p)=x_{g(p)}=x_{\tg(i)}$ by definition. This means it is equal to the facet $\tg G^\lambda_{p,i}$. Thus we may assume that $g$ fixes $x_p$, and hence restricts to a permutation in $S^0_{p-1}$. This gives us
			\begin{align*}
				F-x_p&= gG_{p,i}^\lambda - x_p\\
				&=g(G_{p,i}^\lambda - x_p)\\
				&=gG_{p-1,i}^\lambda
			\end{align*}
			which is a facet of $\calP(a,p-1,m)$.

			Conversely, for any generating facet $G^\lambda_{p-1,i}$ of $\calP(a,p-1,m)$ and any permutation $h$ in $S^0_{p-1}$ we can view $h$ as a permutation in $S^0_p$ which fixes $p$, and hence we have $hG^\lambda_{p-1,i}\sqcup \{x_p\}=hG^\lambda_{p,i}$, which is a facet of $\calP(a,p,m)$.
			
			For the $m=p+1$ case note that the set $G_{p,p+1}$ is empty, and therefore for every partition $\lambda \in \PP(a,p+1)$ the generating facet $G^\lambda_{p,p+1}$ contains no boundary vertices. Thus the only facets of $\calP(a,p,p+1)$ which contain $x_p$ are those of the form $gG^\lambda_{p,i}$ for $1\leq i \leq p$ and permutations $g\in S^0_p$. The remainder of the proof is identical to the earlier case.
		\end{proof}
		
		\begin{cor}\label{Corollary: Partition Complex Link of sigma_X}
			Let $a$, $p$ and $m$ be integers with $a\geq 2$ and $1\leq m\leq p+1$, and let $\Delta=\calP(a,p,m)$. Let $\sigma$ be a face of $\Delta$ of size $\alpha$ contained entirely inside $X_p$. We have an isomorphism of complexes
			$$\lkds \cong \begin{cases}
				\calP(a,p-\alpha,m)& \text{ if } 0\leq \alpha \leq p - m\\
				\calP(a,p-\alpha,p+1-\alpha)& \text{ otherwise.}
			\end{cases}$$
		\end{cor}
		\begin{proof}
			We proceed by induction on $\alpha\geq 0$. The base case $\alpha=0$ is immediate.
			
			Now suppose $\alpha \geq 1$. By symmetry we may assume that $x_p\in \sigma$: indeed, if it is not we may choose some permutation $g\in S^0_p$ such that $x_p\in g\sigma$, and by symmetry we have an isomorphism of complexes $\lkds\cong \lkds[g\sigma]$.
			
			We consider the two cases separately. First we assume that $1\leq \alpha \leq p-m$. In particular this means that $m\leq p-1$, and hence Lemma \ref{Lemma: Partition Complex Link of x} tells us that $\lkds[x_p] = \calP(a,p-1,m)$. Thus $\lkds$ is equal to $\link_{\calP(a,p-1,m)}(\sigma-x_p)$. Note that $0\leq \alpha -1 \leq (p-1)-m$, and so the inductive hypothesis gives us an isomorphism $\link_{\calP(a,p-1,m)}(\sigma-x_p)\cong \calP(a,p-\alpha,m)$.
			
			Now suppose $\alpha > p-m$. By Lemma \ref{Lemma: Partition Complex Link of x}, we know that $\lkds$ is equal to either $\link_{\calP(a,p-1,m)}(\sigma-x_p)$ or $\link_{\calP(a,p-1,p)}(\sigma-x_p)$. We know that $\alpha-1$ is greater than both $(p-1)-m$ and $(p-1)-p$, so in either case the inductive hypothesis gives us an isomorphism $\lkds\cong \calP(a,p-\alpha,p+1-\alpha)$.
		\end{proof}
		
		Corollary \ref{Corollary: Partition Complex Link of sigma_X} allows us to restrict to the case where $\sigma_X=\emptyset$, because it shows us that the link of $\sigma_X$ is also a partition complex, and thus the link of $\sigma$ in $\calP(a,p,m)$ is isomorphic to the link of $\sigma_Y$ in a smaller partition complex.
		
		%
		
		The following lemma allow us to restrict our attention even further.
		
		\begin{lem}\label{Lemma: Partition incomplete => acyclic}
			Let $a$, $p$ and $m$ be positive integers with $m\leq p+1$, and let $\Delta = \calP(a,p,m)$. If $\sigma$ is a face of $\Delta$ which is not partition complete then $\lkds$ is acyclic.
		\end{lem}
		\begin{proof}
			Suppose $\sigma\in \Delta$ is not partition complete. This means there exist some $1\leq i \leq p$ and $1\leq j<j'\leq a-1$ such that the parition vertex $y_i^{j'}$ is in $\sigma$ but $y_i^j$ is not. Because every facet of $\Delta$ is partition complete, every facet of $\Delta$ containing $\sigma$ also contains $y_i^j$. Thus $\lkds$ is a cone over $y_i^j$.
		\end{proof}
		
		It now only remains for us to consider the case where $\sigma$ is a nonempty partition complete face of $\Delta$ contained entirely in $Y_p^a$. We do this using the following lemma and proposition.
		
		\begin{lem}\label{Lemma: Partition Vertices Deformation Retract}
			Let $a$, $p$ and $m$ be positive integers with $m\leq p+1$, and let $\Delta = \calP(a,p,m)$. Fix some $\sigma \in \Delta$. For each $1\leq i\leq p+1$ such that there exists a partition vertex $y_i^j$ which is not contained in $\sigma$, we define $\varphi(\sigma,i)=\min \{1\leq j \leq a-1 : y_i^j \notin \sigma\}$. There is a deformation retraction $\Phi$ of $\lkds$ onto a complex $\Gamma$ obtained by 
			collapsing every partition vertex $y_i^j$ in $\lkds$ onto the
			partition vertex $y_i^{\varphi(\sigma,i)}$.
		\end{lem}
		\begin{proof}
			Let $y_i^j$ be a partition vertex in $\lkds$. This partition vertex cannot be contained in $\sigma$, and hence the index $\varphi(\sigma,i)$ is well-defined. For any facet $F$ of $\lkds$ containing $y_i^j$, we have that $F\sqcup \sigma$ is a facet of $\Delta$. This means $F\sqcup \sigma$ must be partition complete, and must therefore also contain the vertex $y_i^{\varphi(\sigma,i)}$. But $y_i^{\varphi(\sigma,i)}\notin \sigma$ by definition, and so it must be in $F$. Thus any facet of $\lkds$ containing $y_i^j$ also contains $y_i^{\varphi(\sigma,i)}$, and the result follows from Corollary \ref{Corollary: Deformation Retract}. 
		\end{proof}
		
		\begin{prop}\label{Proposition: Partition Complex Skeleton Link}
			Let $a$, $p$ and $m$ be positive integers with $m\leq p+1$, and let $\Delta = \calP(a,p,m)$. Let $\sigma$ be a nonempty partition complete face of $\Delta$ contained entirely inside $Y_p^a$, and set $\beta = |\{y_i^j \in \sigma : j = 1\}|$ and $\gamma = |\{y_i^j \in \sigma : j > 1\}|$. We have a deformation retraction
			$$\lkds \rightsquigarrow \barP(2,p-\beta,m-\beta)\ast \Skel_{a-\gamma-3}([\beta]).$$
		\end{prop}
		\begin{ex}\label{Example: Partition Complex Skeleton Link Example}
			Proposion \ref{Proposition: Partition Complex Skeleton Link} is a rather technical result. To elucidate it we begin by presenting a specific example, using the grid diagrams from Notation \ref{Notation: Grid Notation For Partition Complex Faces}.
			
			Consider the complex $\Delta=\calP(6,5,4)$, and let $\sigma$ be the partition complete face
			
			\begin{center}
				\begin{tikzpicture}[xscale=.7,yscale=.5]
					\draw (0.5,5.5) node{$0$};
					\draw (1.5,5.5) node{$1$};
					\draw (2.5,5.5) node{$2$};
					\draw (3.5,5.5) node{$3$};
					\draw (4.5,5.5) node{$4$};
					\draw (5.5,5.5) node{$5$};
					\fill[gray] (3,4) rectangle (4,3);
					\draw (3.5,3.5) node{$y_3^1$};
					\fill[gray] (4,4) rectangle (5,3);
					\draw (4.5,3.5) node{$y_4^1$};
					\fill[gray] (5,4) rectangle (6,3);
					\draw (5.5,3.5) node{$y_5^1$};
					\fill[gray] (4,3) rectangle (5,2);
					\draw (4.5,2.5) node{$y_4^2$};
					\fill[gray] (5,3) rectangle (6,2);
					\draw (5.5,2.5) node{$y_5^2$};
					\draw[step=1cm,black,very thin] (0,5) grid (6,3);
					\draw[step=1cm,gray,dotted] (0,3) grid (6,2);
				\end{tikzpicture}
			\end{center}
			which is contained entirely in the set $Y^6_5$. Note that for this example we have $\beta = 3$ and $\gamma=2$. Thus $\barP(2,p-\beta,m-\beta)=\barP(2,2,1)$ and $\Skel_{a-\gamma-3}([\beta])=\Skel_1([3])$.
			
			To find the facets of $\lkds$ we need to look for facets of $\Delta$ which contain $\sigma$. By Lemma \ref{Lemma: Partition Complex Description of Facets}, these facets are all totally separated, partition complete, and contain $6+5-1=10$ vertices, no more than $4$ of which are lower partition vertices. By Remark \ref{Remark: Partitions of a-2}, exactly $6-2=4$ of their vertices are upper partition vertices. Below are three examples.
			
			\begin{tabular}{ c c c c c }
				\begin{tikzpicture}[xscale=.55,yscale=.5]
					\fill[lightgray] (3,4) rectangle (4,3);
					\draw (3.5,3.5) node{$y_3^1$};
					\fill[lightgray] (4,4) rectangle (5,3);
					\draw (4.5,3.5) node{$y_4^1$};
					\fill[lightgray] (5,4) rectangle (6,3);
					\draw (5.5,3.5) node{$y_5^1$};
					\fill[lightgray] (4,3) rectangle (5,2);
					\draw (4.5,2.5) node{$y_4^2$};
					\fill[lightgray] (5,3) rectangle (6,2);
					\draw (5.5,2.5) node{$y_5^2$};
					
					\fill[gray] (0,5) rectangle (1,4);
					\draw (0.5,4.5) node{$x_0$};
					\fill[gray] (1,5) rectangle (2,4);
					\draw (1.5,4.5) node{$x_1$};
					\fill[gray] (2,4) rectangle (3,3);
					\draw (2.5,3.5) node{$y_2^1$};
					\fill[gray] (3,3) rectangle (4,2);
					\draw (3.5,2.5) node{$y_3^2$};
					\fill[gray] (5,2) rectangle (6,1);
					\draw (5.5,1.5) node{$y_5^3$};
					
					\draw[step=1cm,black,very thin] (0,5) grid (6,3);
					\draw[step=1cm,gray,dotted] (0,3) grid (6,1);
				\end{tikzpicture} && \begin{tikzpicture}[xscale=.55,yscale=.5]
					\fill[lightgray] (3,4) rectangle (4,3);
					\draw (3.5,3.5) node{$y_3^1$};
					\fill[lightgray] (4,4) rectangle (5,3);
					\draw (4.5,3.5) node{$y_4^1$};
					\fill[lightgray] (5,4) rectangle (6,3);
					\draw (5.5,3.5) node{$y_5^1$};
					\fill[lightgray] (4,3) rectangle (5,2);
					\draw (4.5,2.5) node{$y_4^2$};
					\fill[lightgray] (5,3) rectangle (6,2);
					\draw (5.5,2.5) node{$y_5^2$};
					
					\fill[gray] (0,5) rectangle (1,4);
					\draw (0.5,4.5) node{$x_0$};
					\fill[gray] (2,5) rectangle (3,4);
					\draw (2.5,4.5) node{$x_2$};
					\fill[gray] (1,4) rectangle (2,3);
					\draw (1.5,3.5) node{$y_1^1$};
					\fill[gray] (1,3) rectangle (2,2);
					\draw (1.5,2.5) node{$y_1^2$};
					\fill[gray] (4,2) rectangle (5,1);
					\draw (4.5,1.5) node{$y_4^3$};
					
					\draw[step=1cm,black,very thin] (0,5) grid (6,3);
					\draw[step=1cm,gray,dotted] (0,3) grid (6,1);
				\end{tikzpicture} && \begin{tikzpicture}[xscale=.55,yscale=.5]
					\fill[lightgray] (3,4) rectangle (4,3);
					\draw (3.5,3.5) node{$y_3^1$};
					\fill[lightgray] (4,4) rectangle (5,3);
					\draw (4.5,3.5) node{$y_4^1$};
					\fill[lightgray] (5,4) rectangle (6,3);
					\draw (5.5,3.5) node{$y_5^1$};
					\fill[lightgray] (4,3) rectangle (5,2);
					\draw (4.5,2.5) node{$y_4^2$};
					\fill[lightgray] (5,3) rectangle (6,2);
					\draw (5.5,2.5) node{$y_5^2$};
					
					\fill[gray] (0,5) rectangle (1,4);
					\draw (0.5,4.5) node{$x_0$};
					\fill[gray] (1,5) rectangle (2,4);
					\draw (1.5,4.5) node{$x_1$};
					\fill[gray] (2,5) rectangle (3,4);
					\draw (2.5,4.5) node{$x_2$};
					\fill[gray] (3,3) rectangle (4,2);
					\draw (3.5,2.5) node{$y_3^2$};
					\fill[gray] (3,2) rectangle (4,1);
					\draw (3.5,1.5) node{$y_3^3$};
					
					\draw[step=1cm,black,very thin] (0,5) grid (6,3);
					\draw[step=1cm,gray,dotted] (0,3) grid (6,1);
				\end{tikzpicture}
			\end{tabular}
			
			The corresponding facets of $\lkds$ can be found by removing $\sigma$ from these faces. By Lemma \ref{Lemma: Partition Vertices Deformation Retract} we have a deformation retraction $\Phi$ from $\lkds$ onto to a complex $\Gamma$ obtained from the map $y_i^j \mapsto y_i^{\varphi(\sigma,i)}$. In particular we have $\varphi(\sigma,0)=\varphi(\sigma,1)=\varphi(\sigma,2)=1$, $\varphi(\sigma,3)=2$, and $\varphi(\sigma,4)=\varphi(\sigma,5)=3$. Under this map, the faces of $\Gamma$ corresponding to the three facets above are
			
			\begin{tabular}{ c c c c c }
				\begin{tikzpicture}[xscale=.55,yscale=.5]
					\fill[gray] (0,5) rectangle (1,4);
					\draw (0.5,4.5) node{$x_0$};
					\fill[gray] (1,5) rectangle (2,4);
					\draw (1.5,4.5) node{$x_1$};
					\fill[gray] (2,4) rectangle (3,3);
					\draw (2.5,3.5) node{$y_2^1$};
					\fill[gray] (3,3) rectangle (4,2);
					\draw (3.5,2.5) node{$y_3^2$};
					\fill[gray] (5,2) rectangle (6,1);
					\draw (5.5,1.5) node{$y_5^3$};
					
					\draw[step=1cm,black,very thin] (0,5) grid (6,3);
					\draw[step=1cm,gray,dotted] (0,3) grid (6,1);
				\end{tikzpicture} && \begin{tikzpicture}[xscale=.55,yscale=.5]
					\fill[gray] (0,5) rectangle (1,4);
					\draw (0.5,4.5) node{$x_0$};
					\fill[gray] (2,5) rectangle (3,4);
					\draw (2.5,4.5) node{$x_2$};
					\fill[gray] (1,4) rectangle (2,3);
					\draw (1.5,3.5) node{$y_1^1$};
					\fill[gray] (4,2) rectangle (5,1);
					\draw (4.5,1.5) node{$y_4^3$};
					
					\draw[step=1cm,black,very thin] (0,5) grid (6,3);
					\draw[step=1cm,gray,dotted] (0,3) grid (6,1);
				\end{tikzpicture} && \begin{tikzpicture}[xscale=.55,yscale=.5]
					\fill[gray] (0,5) rectangle (1,4);
					\draw (0.5,4.5) node{$x_0$};
					\fill[gray] (1,5) rectangle (2,4);
					\draw (1.5,4.5) node{$x_1$};
					\fill[gray] (2,5) rectangle (3,4);
					\draw (2.5,4.5) node{$x_2$};
					\fill[gray] (3,3) rectangle (4,2);
					\draw (3.5,2.5) node{$y_3^2$};
					
					\draw[step=1cm,black,very thin] (0,5) grid (6,3);
					\draw[step=1cm,gray,dotted] (0,3) grid (6,1);
				\end{tikzpicture}
			\end{tabular}
			
			All three of these faces may be decomposed into the disjoint union of a face whose support is contained in $\{0,1,2\}$ (which is disjoint from the support of $\sigma$) and a face whose support is contained in $\{3,4,5\}$ (which is the support of $\sigma$). For example, the first face may be decomposed into $\{x_0,x_1,y^1_2\}\sqcup\{y^2_3,y^3_5\}$.
			
			Moreover, if we define
			\begin{align*}
				A &= \big\{ \tau \in \Gamma : \Supp(\tau) \subseteq \{0,1,2\} \}\\
				B &= \big\{ \tau \in \Gamma : \Supp(\tau) \subseteq \{3,4,5\} \big\},
			\end{align*}
			then the disjoint union of any face in $A$ with any face in $B$ must be a face of $\Gamma$ -- because its union with $\sigma$ is separated, partition complete, has no more than $10$ vertices, and no more than $4$ lower partition vertices. We conclude that $\Gamma$ is equal to the join of complexes $A\ast B$.
			
			The faces of $A$ are all the separated subsets of $V^2_5$ with at most $1$ partition vertex (including the set $\{x_0,x_1,x_2\}$). Thus we have $A=\barP(2,2,1)$. The faces of $B$ are the subsets of $\{y_3^2,y_4^3,y_5^3\}$ of size at most $2$. In other words, $B$ is equal to the $1$-skeleton of this set, which is isomorphic to $\Skel_1([3])$.
			
			The proof below is a generalisation of this argument.
		\end{ex}
		\begin{proof}     
			By symmetry we may assume that $\Supp(\sigma) = \{p-\beta+1,\dots,p\}$.
			
			We begin by applying the deformation retraction $\Phi$ from Lemma \ref{Lemma: Partition Vertices Deformation Retract} to the complex $\lkds$ to get a new complex $\Gamma$ whose only partition vertices are the vertices $y_0^{\varphi(0,\sigma)}, \dots, y_p^{\varphi(p,\sigma)}$.
			
			Let the subcomplexes $A$ and $B$ of $\Gamma$ be defined as follows.
			\begin{align*}
				A &= \big\{ \tau \in \Gamma : \Supp(\tau) \subseteq \{0,\dots,p-\beta\}\big\}\\
				B &= \big\{ \tau \in \Gamma : \Supp(\tau) \subseteq \{p-\beta+1,\dots,p\}\big\}.
			\end{align*}
			
			We aim to show the following.
			\begin{enumerate}
				\item $A = \barP(2,p-\beta,m-\beta)$.
				\item $B = \Skel_{a-\gamma-3}(S)$ where $S$ denotes the set $\{y_i^{\varphi(i,\sigma)}: i \in \Supp(\sigma)\}$.
				\item $\Gamma = A \ast B$.
			\end{enumerate}
			
			For part (1), note first that because $\Supp(\sigma)=\{p-\beta+1,\dots,p\}$, we have $\varphi(i,\sigma)=1$ for each $0\leq i \leq p-\beta$. Thus, the only partition vertices of $A$ are the vertices $y^1_0,\dots,y^1_{p-\beta}$, and so $A$ has vertex set $V^2_{p-\beta}=\{x_0,\dots,x_{p-\beta},y_0,\dots,y_{p-\beta}\}$. Suppose $\tau$ is a face of $A$. Because $\sigma \sqcup \tau$ is a face of $\Delta$, we have $|\Supp(\sigma_Y\sqcup \tau_Y)|\leq m$ by Lemma \ref{Lemma: Partition Complex Description of Facets} condition (4). We know $|\Supp(\sigma_Y)|=\beta$ so this means $|\Supp(\tau_Y)|\leq m-\beta$. Thus the faces of $A$ are all the separated subsets $\tau$ of $V^2_{p-\beta}$ such that $0\leq |\Supp(\tau_Y)|\leq m-\beta$ (including the subset $\{x_0,\dots,x_{p-\beta}\}$). These are precisely the faces of $\barP(2,p-\beta,m-\beta)$.
			
			Now we move on to part (2). Because $\sigma$ consists entirely of partition vertices, the faces of $B$ must also consist entirely of partition vertices (otherwise their union with $\sigma$ would not be separated). In particular, every face of $B$ must be a subset of $S=\{y_i^{\varphi(i,\sigma)}: i \in \Supp(\sigma)\}$. From Remark \ref{Remark: Partitions of a-2} we know the faces of $\Delta$ consisting entirely of upper partition vertices have size at most $a-2$. Thus the faces of $B$ have size at most $a-\gamma - 2$ (i.e. dimension at most $a-\gamma-3$). Moreover, for every subset $\tau \subset S$, the union $\sigma \sqcup \tau$ is partition complete. Hence by Lemma \ref{Lemma: Partition Complex Description of Facets}, $B$ contains \textit{every} subset of $S$ of size at most $a-\gamma-2$, and is therefore equal to $\Skel_{a-\gamma-3}(S)$.
			
			Finally we consider part (3). For any two faces $\tau\in A$ and $\rho\in B$, the disjoint union $\tau\sqcup \rho \sqcup \sigma$ is partition complete and separated. From the above discussion we also have that $|\tau_Y|\leq m-\beta$ and $|\rho_Y|\leq a-\gamma - 2$. This means both that $|\tau_Y \sqcup \rho_Y \sqcup \sigma|\leq a+m-2$ and $1\leq |\Supp(\tau_Y \sqcup \rho_Y \sqcup \sigma)|\leq m$. We conclude from Lemma \ref{Lemma: Partition Complex Description of Facets} that the disjoint union $\tau\sqcup \rho \sqcup \sigma$ is a face of $\Delta$, which shows that $\tau \sqcup \rho$ is a face of $\Gamma$. The result follows.
		\end{proof}

		\subsection{Proving Theorem \ref{Theorem: Partition Complex Degree Type}}\label{Subsection: Proving Partition Complex Theorem}
		We now have all the ingredients we need to prove Theorem \ref{Theorem: Partition Complex Degree Type}. We start by proving the following corollaries to our results about the links of $\Delta$ in the last section.
		
		\begin{cor}\label{Corollary: Partition Complex Homology of Link of sigma_X}
			Let $a$, $p$ and $m$ be integers with $a\geq 2$ and $1\leq m\leq p+1$, and let $\Delta=\calP(a,p,m)$. Let $\sigma$ be a face of $\Delta$ contained entirely inside $X_p$. We have
			$$h(\Delta,\sigma) = \begin{cases}
				\{p-|\sigma|-1\}& \text{ if } 0\leq |\sigma| \leq p - m\\
				\emptyset& \text{ otherwise.}
			\end{cases}$$
		\end{cor}
		\begin{proof}
			We set $\alpha = |\sigma|$. By corollary \ref{Corollary: Partition Complex Link of sigma_X}, $\lkds$ is isomorphic to $\calP(a,p-\alpha,m)$ if $0\leq \alpha \leq p-m$ or $\calP(a,p-\alpha,p+1-\alpha)$ otherwise. By Corollary \ref{Corollary: Homology of Partition Complexes}, the former of these has $\nth[st]{(p-\alpha- 1)}$ homology while the latter is acyclic.
		\end{proof}
		
		\begin{cor}\label{Corollary: Partition Complex Homology of Skeleton Link}
			Let $a$, $p$ and $m$ be integers with $a\geq 2$ and $1\leq m\leq p+1$, and let $\Delta = \calP(a,p,m)$. Let $\sigma$ be a partition complete face of $\Delta$ with $\sigma_Y\neq \emptyset$. We have
			$$h(\Delta,\sigma)=\begin{cases}
				\{a+p - |\sigma| - 2\} & \text{ if } |\sigma_X| \geq p-m + 1 \text{ and } |\sigma_Y| \geq a-1 \\
				\emptyset & \text{ otherwise.}
			\end{cases}$$
		\end{cor}
		\begin{proof}
			We set $\alpha = |\sigma_X|$, $\beta = |\{y_i^j \in \sigma_Y : j = 1\}|$ and $\gamma = |\{y_i^j \in \sigma_Y : j > 1\}|$. Because $\sigma$ is partition complete and $\sigma_Y\neq \emptyset$, we know $\beta > 0$.
			
			First suppose that $0\leq \alpha \leq p-m$. Corollary \ref{Corollary: Partition Complex Link of sigma_X} tells us that $\lkds[\sigma_X]$ is isomorphic to $\calP(a,p-\alpha,m)$. Therefore, by Proposition \ref{Proposition: Partition Complex Skeleton Link}, there is a deformation retraction
			$$\lkds \rightsquigarrow \barP(2,p-\alpha -\beta,m-\beta)\ast \Skel_{a-\gamma-3}([\beta]).$$
			The complex $\barP(2,p-\alpha -\beta,m-\beta)$ is acyclic by Corollary \ref{Corollary: Homology of Partition Complexes}, and thus so is $\lkds$.
			
			Now suppose that $\alpha \geq p-m+1$. Corollary \ref{Corollary: Partition Complex Link of sigma_X} tells us that $\lkds[\sigma_X]$ is isomorphic to $\calP(a,p-\alpha,p+1 - \alpha)$, and hence by Proposition \ref{Proposition: Partition Complex Skeleton Link}, there is a deformation retraction
			$$\lkds \rightsquigarrow \barP(2,p-\alpha -\beta,p+1-\alpha -\beta)\ast \Skel_{a-\gamma-3}([\beta]).$$
			The complex $\barP(2,p-\alpha -\beta,p+1-\alpha-\beta)$ has only $\nth[th]{(p-\alpha-\beta)}$ homology, by Corollary \ref{Corollary: Homology of Partition Complexes}. Meanwhile the complex $\Skel_{a-\gamma-3}([\beta])$ has only $\nth[rd]{(a-\gamma-3)}$ homology so long as  $a- \gamma - 3 < \beta - 1$ (otherwise it is equal to the full simplex on $[\beta]$, which is acyclic). Thus, $\lkds$ has nontrivial homology if and only if $\beta + \gamma \geq a - 1$, and in this case we have, by Proposition \ref{Proposition: Kunneth Formula}, $h(\Delta,\sigma) = \{p - \alpha - \beta\} \boxplus  \{a - \gamma - 3\} \boxplus \{1\} = \{a + p -|\sigma| - 2\}$ (where $\boxplus$ represents the addition of sets operation, as laid out in Definition \ref{Definition: Boxplus}).
		\end{proof}
		
		With these two corollaries in our toolkit, we now proceed to the proof of Theorem \ref{Theorem: Partition Complex Degree Type}.
		\begin{proof}[Proof of Theorem \ref{Theorem: Partition Complex Degree Type}]
			Let $a$, $p$ and $m$ be positive integers with $m\leq p$ and set $\Delta=\calP(a,p,m)$. Fix a face $\sigma \in \Delta$.
			
			By Lemma \ref{Lemma: Partition incomplete => acyclic}, $\lkds$ is acyclic unless $\sigma$ is partition complete. In the case where $\sigma$ is partition complete, Corollaries \ref{Corollary: Partition Complex Homology of Link of sigma_X} and \ref{Corollary: Partition Complex Homology of Skeleton Link} tell us that
			\begin{equation*}\label{Equation: Homology Index Sets Partition Complexes}
				h(\Delta,\sigma)=\begin{cases}
					\{p-|\sigma| - 1\} & \text{ if } 0\leq |\sigma_X|\leq p-m \text{ and } |\sigma_Y| = 0\\
					\{a+p - |\sigma| - 2\} & \text{ if } |\sigma_X| \geq p-m + 1 \text{ and } |\sigma_Y| \geq a-1 \\
					\emptyset & \text{ otherwise.}
				\end{cases}
			\end{equation*}
			
			In particular, this means that for any integer $k$, we have
			\begin{equation*}\label{Equation: Complete Homology Index Sets Partition Complexes}\hh(\Delta,k)=\begin{cases}
					\{p- k - 1\} & \text{ if } 0\leq k\leq p-m\\
					\{a+p - k - 2\} & \text{ if } a + p - m \leq k \leq a+ p - 1 \\
					\emptyset & \text{ otherwise.}
				\end{cases}
			\end{equation*}
			This proves $\Delta$ is PR with degree type $(\overbrace{1,\dots,1\underbrace{a,1,\dots,1}_{m}}^{p})$, by Proposition \ref{Proposition: Alternate PR Definition With Degree Type}.
		\end{proof}
		
		\section{Pure Resolutions of Any Degree Type}\label{Section: Pure resolutions of any degree type}
		This section is devoted to the proof of Theorem \ref{Theorem: PR Complexes of Any Degree Type}. Our method for this proof is focussed around finding operations on simplicial complexes which preserve the PR property in specific cases, while altering the degree types of PR complexes in specified ways.
		
		Specifically, we will construct a family of operations on simplicial complexes $\{\phi_i : i \in \ZZ^+\}$ with the property that for each $i\in \ZZ^+$, and any PR complex $\Delta$ with degree type $(d_p,\dots,d_i,\underbrace{1,\dots,1}_{i-1})$ for some $p\geq i$, the complex $\phi_i(\Delta)$ is PR with degree type $(d_p,\dots,d_{i-1},d_i+1,\underbrace{1,\dots,1}_{i-1})$. If $\Delta$ is PR with a degree type of a different form, then $\phi_i(\Delta)$ will not be PR.
		
		To see why the existence of such a family is sufficient to prove Theorem \ref{Theorem: PR Complexes of Any Degree Type}, suppose we wish to construct a PR complex of degree type $\bd=(d_p,\dots,d_1)$. We can do so by taking a PR complex $\Delta$ of degree type $(\underbrace{1,\dots,1}_{p})$ (such as the boundary of the $p$-simplex in Example \ref{Example: PR Simplex}), and applying each of the operations $\phi_i$ to $\Delta$ a total of $(d_i-1)$ times.
		
		In other words, if $\Delta$ is the boundary of the $p$-simplex, then the simplicial complex $\phi_1^{d_1-1}\dots \phi_p^{d_p-1} (\Delta)$ is PR with degree type $(d_p,\dots,d_1)$.
		
		\subsection{The $\phi_i$ Operations}
		We define the operation $\phi_i$ below.
		
		\begin{defn}\label{Definition: Phi_i}
			Let $\Delta$ be a simplicial complex on vertex set $V$, and fix a positive integer $i$. The complex $\phi_i(\Delta)$ is defined as follows.
			\begin{enumerate}
				\item $S^i_\Delta:=\{u_{\sigma}: \sigma \in \Delta, |\sigma|\geq \dim \Delta + 2 -i\}$.
				\item $\phi_i(\Delta)$ is the complex on vertex set $V\sqcup S^i_\Delta$ obtained by adding to $\Delta$ all those faces $\sigma \sqcup \{u_{\tau_1},\dots u_{\tau_r}\}$ for each sequence $\sigma \subseteq \tau_1 \subset \dots \subset \tau_r$ in $\Delta$ for which the vertices $u_{\tau_1},\dots, u_{\tau_r}$ are in $S^i_{\Delta}$.
			\end{enumerate}
		\end{defn}
		\begin{rem}\label{Remark: Phi_i = Phi_j}
			If $m= \dim \Delta + 2$ then for any $i>m$ we have $\phi_i(\Delta)=\phi_m(\Delta)$. Note also that $m$ is the smallest integer for which the set $S^m_\Delta$ contains the vertex $u_\emptyset$.
		\end{rem}
		
		\begin{ex}
			The complex $\phi_1(\Delta)$ is obtained from $\Delta$ by adding the faces $F\cup \{u_F\}$ for each facet $F$ in $\Delta$ of maximal dimension. In other words, $\phi_1$ acts on a complex by adding an additional free vertex $u_F$ to each facet $F$ of maximal dimension. In particular, because all PR complexes are pure, $\phi_1$ acts on all of them by adding an additional free vertex to every one of their facets.
		\end{ex}
		
		\begin{ex}\label{Example: phi_i(Delta)}
			Let $\Delta$ be the boundary of the $2$-simplex on vertex set $\{x,y,z\}$:
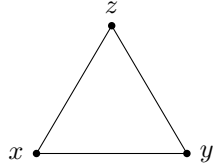
			\begin{center}
				\begin{tikzpicture}[scale = 1]
					\tikzstyle{point}=[circle,thick,draw=black,fill=black,inner sep=0pt,minimum width=2pt,minimum height=2pt]
					\node (a)[point, label=left:$x$] at (0,0) {};
					\node (b)[point, label=right:$y$] at (2,0) {};
					\node (c)[point, label=above:$z$] at (1,1.7) {};
					
					\draw (a.center) -- (b.center) -- (c.center) -- cycle;
				\end{tikzpicture}
			\end{center}
			This is PR with degree type $(1,1)$.
			
			The complex $\phi_1(\Delta)$ is
			
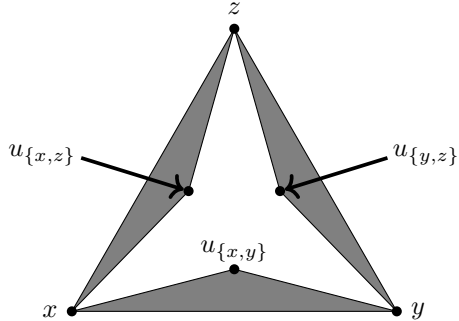
\begin{center}
				\begin{tikzpicture}[scale = 0.8]
					\tikzstyle{point}=[circle,thick,draw=black,fill=black,inner sep=0pt,minimum width=3pt,minimum height=3pt]
					\node (x)[point, label=left:$x$] at (-1.7,-1) {};
					\node (y)[point, label=right:$y$] at (3.7,-1) {};
					\node (z)[point, label=above:$z$] at (1,3.7) {};
					
					
					\node (d)[point] at (1,-0.3) {}; 
					\node (e)[point] at (1.76,1) {}; 
					\node (f)[point] at (0.24,1) {}; 
					
					
					
					\node (uxy)[label = {[label distance=-1.5mm] above:$u_{\{x,y\}}$}] at (1,-0.3) {};
					\node (uyz)[label = {[label distance=-3mm] right:$u_{\{y,z\}}$}] at (3.7,1.6) {};
					\node (uxz)[label = {[label distance=-3mm] left:$u_{\{x,z\}}$}] at (-1.7,1.6) {};
					
					
					\begin{scope}[on background layer]
						
						\draw[fill=gray,fill opacity=.8] (x.center) -- (y.center) -- (d.center) -- cycle;
						\draw[fill=gray,fill opacity=.8] (y.center) -- (z.center) -- (e.center) -- cycle;
						\draw[fill=gray,fill opacity=.8] (x.center) -- (z.center) -- (f.center) -- cycle;
						
					\end{scope}
					
					\draw[->, line width=1.4pt] (uyz)  -- (e);
					\draw[->, line width=1.4pt] (uxz)  -- (f);
					
				\end{tikzpicture}
			\end{center}
			which is PR with degree type $(1,2)$.
			
			The complex $\phi_2(\Delta)$ is
			
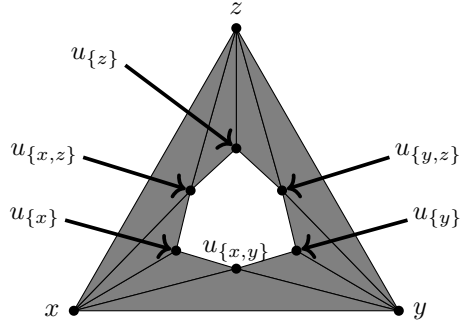
\begin{center}
				\begin{tikzpicture}[scale = 0.8]
					\tikzstyle{point}=[circle,thick,draw=black,fill=black,inner sep=0pt,minimum width=3pt,minimum height=3pt]
					\node (x)[point, label=left:$x$] at (-1.7,-1) {};
					\node (y)[point, label=right:$y$] at (3.7,-1) {};
					\node (z)[point, label=above:$z$] at (1,3.7) {};
					
					\node (a)[point] at (0,0) {}; 
					\node (b)[point] at (2,0) {}; 
					\node (c)[point] at (1,1.7) {}; 
					
					\node (d)[point] at (1,-0.3) {}; 
					\node (e)[point] at (1.76,1) {}; 
					\node (f)[point] at (0.24,1) {}; 
					
					
					\node (ux)[color=red, label = {[label distance=-3mm] left:$u_{\{x\}}$}] at (-2,0.55) {};
					\node (uy)[label = {[label distance=-3mm] right:$u_{\{y\}}$}] at (4,0.55) {};
					\node (uz)[label = {[label distance=-3mm] left:$u_{\{z\}}$}] at (-1,3.2) {};

					\node (uxy)[label = {[label distance=-1.5mm] above:$u_{\{x,y\}}$}] at (1,-0.3) {};
					\node (uyz)[label = {[label distance=-3mm] right:$u_{\{y,z\}}$}] at (3.7,1.6) {};
					\node (uxz)[label = {[label distance=-3mm] left:$u_{\{x,z\}}$}] at (-1.7,1.6) {};
					
					
					\begin{scope}[on background layer]
						\draw[fill=gray,fill opacity=.8] (x.center) -- (a.center) -- (d.center) -- cycle;
						\draw[fill=gray,fill opacity=.8] (x.center) -- (a.center) -- (f.center) -- cycle;
						\draw[fill=gray,fill opacity=.8] (y.center) -- (b.center) -- (d.center) -- cycle;
						\draw[fill=gray,fill opacity=.8] (y.center) -- (b.center) -- (e.center) -- cycle;
						\draw[fill=gray,fill opacity=.8] (z.center) -- (c.center) -- (e.center) -- cycle;
						\draw[fill=gray,fill opacity=.8] (z.center) -- (c.center) -- (f.center) -- cycle;
						
						\draw[fill=gray,fill opacity=.8] (x.center) -- (y.center) -- (d.center) -- cycle;
						\draw[fill=gray,fill opacity=.8] (y.center) -- (z.center) -- (e.center) -- cycle;
						\draw[fill=gray,fill opacity=.8] (x.center) -- (z.center) -- (f.center) -- cycle;
						
					\end{scope}
					
					\draw[->, line width=1.4pt] (ux)   -- (a);
					\draw[->, line width=1.4pt] (uy)   -- (b);
					\draw[->, line width=1.4pt] (uz)   -- (c);
					\draw[->, line width=1.4pt] (uyz)  -- (e);
					\draw[->, line width=1.4pt] (uxz)  -- (f);
					
				\end{tikzpicture}
			\end{center}
			which is PR with degree type $(2,1)$.
			
			And the complex $\phi_3(\Delta)$ is
			
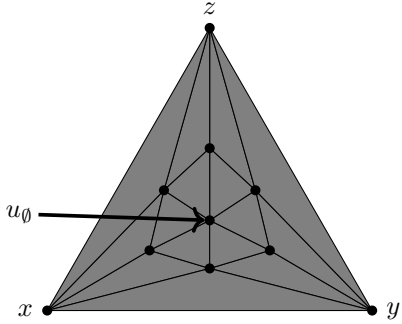
\begin{center}
				\begin{tikzpicture}[scale = 0.8]
					\tikzstyle{point}=[circle,thick,draw=black,fill=black,inner sep=0pt,minimum width=3pt,minimum height=3pt]
					\node (x)[point, label=left:$x$] at (-1.7,-1) {};
					\node (y)[point, label=right:$y$] at (3.7,-1) {};
					\node (z)[point, label=above:$z$] at (1,3.7) {};
					
					\node (a)[point] at (0,0) {}; 
					\node (b)[point] at (2,0) {}; 
					\node (c)[point] at (1,1.7) {}; 
					
					\node (d)[point] at (1,-0.3) {}; 
					\node (e)[point] at (1.76,1) {}; 
					\node (f)[point] at (0.24,1) {}; 
					
					\node (g)[point, label] at (1,0.5) {}; 
					
					
					
					\node (uemp)[label = {[label distance=-3mm] left:$u_{\emptyset}$}] at (-2,0.6) {};
					
					\begin{scope}[on background layer]
						\draw[fill=gray,fill opacity=.8] (x.center) -- (a.center) -- (d.center) -- cycle;
						\draw[fill=gray,fill opacity=.8] (x.center) -- (a.center) -- (f.center) -- cycle;
						\draw[fill=gray,fill opacity=.8] (y.center) -- (b.center) -- (d.center) -- cycle;
						\draw[fill=gray,fill opacity=.8] (y.center) -- (b.center) -- (e.center) -- cycle;
						\draw[fill=gray,fill opacity=.8] (z.center) -- (c.center) -- (e.center) -- cycle;
						\draw[fill=gray,fill opacity=.8] (z.center) -- (c.center) -- (f.center) -- cycle;
						
						\draw[fill=gray,fill opacity=.8] (x.center) -- (y.center) -- (d.center) -- cycle;
						\draw[fill=gray,fill opacity=.8] (y.center) -- (z.center) -- (e.center) -- cycle;
						\draw[fill=gray,fill opacity=.8] (x.center) -- (z.center) -- (f.center) -- cycle;
						
						\draw[fill=gray,fill opacity=.8] (a.center) -- (d.center) -- (g.center) -- cycle;
						\draw[fill=gray,fill opacity=.8] (a.center) -- (f.center) -- (g.center) -- cycle;
						\draw[fill=gray,fill opacity=.8] (b.center) -- (d.center) -- (g.center) -- cycle;
						\draw[fill=gray,fill opacity=.8] (b.center) -- (e.center) -- (g.center) -- cycle;
						\draw[fill=gray,fill opacity=.8] (c.center) -- (e.center) -- (g.center) -- cycle;
						\draw[fill=gray,fill opacity=.8] (c.center) -- (f.center) -- (g.center) -- cycle;
					\end{scope}
					
					\draw[->, line width=1.4pt] (uemp) -- (g);
					
				\end{tikzpicture}
			\end{center}
			Note that the addition of the facets containing $u_\emptyset$ makes this complex acyclic.
		\end{ex}

		We wish to prove the following theorem about the operation $\phi_i$.
		
		\begin{thm}\label{Theorem: Phi_i Operations Degree Type}
			Let $\Delta$ be a PR complex on vertex set $V$, with degree type of the form $(d_p,\dots,d_i,\underbrace{1,\dots,1}_{i-1})$ for some $p\geq i\geq 1$. The complex $\tDelta = \phi_i(\Delta)$ on vertex set $V\sqcup S$ is a PR complex with degree type $(d_p,\dots,d_i+1,\underbrace{1,\dots,1}_{i-1})$
		\end{thm}
		
		The following lemma is particularly crucial. It shows us that in the case where $\Delta$ is pure, $\phi_i$ commutes with taking links of faces in $\Delta$.
		
		\begin{lem}\label{Lemma: Phi_i and Link commute}
			Let $\Delta$ be a pure simplicial complex on vertex set $V$, let $i\geq 1$ and let $\sigma \in \Delta$. We have an isomorphism of complexes $\link_{\phi_i(\Delta)} \sigma \cong \phi_i(\lkds)$.
		\end{lem}
		\begin{proof}
			
			
			
			We claim that the map of vertices $u_{\tau}\mapsto u_{\tau-\sigma}$ gives a well-defined bijection between the facets of $\link_{\phi_i(\Delta)} \sigma$ and the facets of $\phi_i(\lkds)$.
			
			The facets of $\link_{\phi_i(\Delta)} \sigma$ are all of the form $G-\sigma$ for some facet $G$ of $\phi_i(\Delta)$ containing $\sigma$. Let $G = \tau_1 \cup \{u_{\tau_1},\dots, u_{\tau_r}\}$ be one such facet of $\phi_i(\Delta)$, for some $\tau_1 \subset \dots \subset \tau_r$ in $\Delta$ with $|\tau_j|\geq \dim \Delta + 2 - i$ for each $1\leq j \leq r$. Because $G$ contains $\sigma$, we must have $\sigma \subseteq \tau_1$. Hence, for each $1\leq j \leq r$, the simplex $\tau_j-\sigma$ is a face of $\lkds$, and we have
			\begin{align*}
				|\tau_j - \sigma| &= |\tau_j| - |\sigma|\\
				&\geq \dim \Delta + 2 - i - |\sigma|\\
				&= \dim(\lkds) + 2 - i
			\end{align*}
			which shows that $u_{\tau_j-\sigma}$ is a vertex of $S^i_{\lkds}$. Thus the map of vertices $u_{\tau}\mapsto u_{\tau-\sigma}$ gives us a unique corresponding facet $\widetilde{G}=(\tau_1 - \sigma) \cup \{u_{\tau_1-\sigma},\dots, u_{\tau_r-\sigma}\}$ of $\phi_i(\lkds)$.
			
			Conversely, suppose  $H = \rho_1 \cup \{u_{\rho_1},\dots, u_{\rho_r}\}$ is a facet of $\phi_i(\lkds)$, for some $\rho_1 \subset \dots \subset \rho_r$ in $\lkds$ with $|\rho_j|\geq \dim (\lkds) + 2 - i$ for each $1\leq j \leq r$. By the definition of $\lkds$ we know that for each $1\leq j \leq r$,  $\rho_j\sqcup\sigma$ is a face of $\Delta$. Moreover, we have
			\begin{align*}
				|\rho_j \sqcup \sigma| &= |\rho_j| + |\sigma|\\
				&\geq \dim (\lkds) + 2 - i + |\sigma|\\
				&= \dim(\Delta) + 2 - i
			\end{align*}
			which shows that $u_{\rho_j-\sigma}$ is a vertex of $S^i_\Delta$. Thus the map of vertices $u_{\rho}\mapsto u_{\rho\sqcup\sigma}$ gives us a unique corresponding facet $\widehat{H}=\rho_1 \cup \{u_{\rho_1\sqcup\sigma},\dots, u_{\rho_r\sqcup\sigma}\}$ of $\link_{\phi_i(\Delta)}\sigma$.
		\end{proof}
		
		\subsection{Barycentric Subdivision}\label{Subsection: Barycentric Subdivision}
		The $\phi_i$ operations have an important connection to the process of barycentric subdivision. We now give a brief description of the process of barycentric subdivision, and how it links to the $\phi_i$ operations.
		
		\begin{defn}\label{Definition: BDelta}
			We define the \textit{barycentric subdivision} of $\Delta$ to be the complex $\BDelta$ on vertex set $S_\Delta = \{u_{\sigma}: \sigma \in \Delta-\{\emptyset\}\}$, with faces $\{u_{\sigma_1},\dots, u_{\sigma_r}\}$ whenever $\sigma_1 \subset \dots \subset \sigma_r$.
		\end{defn}
		
		Significantly, barycentric subdivision does not affect the topology of a simplicial complex. In other words, the complexes $\BDelta$ and $\Delta$ are homeomorphic as topological spaces. In particular, this means that for any integer $j\geq -1$ we have $$\Hred_j\left (\BDelta\right )=\Hred_j(\Delta).$$
		This fact will be particularly useful to us.
		
		Note that if $i$ is chosen such that the vertex set $S^i_\Delta$ given in Definition \ref{Definition: Phi_i} contains vertices $u_\sigma$ corresponding to every simplex $\sigma$ in $\Delta$ except $\emptyset$, then the induced subcomplex $\phi_i(\Delta)|_{S^i_\Delta}$ is equal to $\BDelta$. This happens when $i = \dim \Delta + 1$ (and when $i\geq \dim \Delta+2$, the induced subcomplex $\phi_i(\Delta)|_{S^i_\Delta}$ is equal to $\BDelta\ast u_\emptyset$).
		
		For example, if $\Delta$ is the boundary of the $2$-simplex, as in Example \ref{Example: phi_i(Delta)}, its barycentric subdivision $\BDelta$ is
		
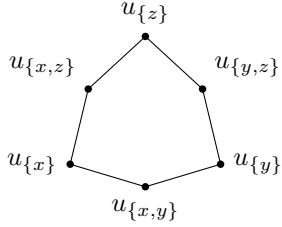
\begin{center}
			\begin{tikzpicture}[scale = 1]
				\tikzstyle{point}=[circle,thick,draw=black,fill=black,inner sep=0pt,minimum width=2pt,minimum height=2pt]
				\node (a)[point, label=left:$u_{\{x\}}$] at (0,0) {};
				\node (b)[point, label=right:$u_{\{y\}}$] at (2,0) {};
				\node (c)[point, label=above:$u_{\{z\}}$] at (1,1.7) {};
				
				\node (d)[point, label=below:$u_{\{x,y\}}$] at (1,-0.3) {};
				\node (e)[point, label=above right:$u_{\{y,z\}}$] at (1.76,1) {};
				\node (f)[point, label=above left:$u_{\{x,z\}}$] at (0.24,1) {};
				
				\draw (a.center) -- (d.center) -- (b.center) -- (e.center) -- (c.center) -- (f.center) -- cycle;
			\end{tikzpicture}
		\end{center}
		which is an induced subcomplex of the complex $\phi_2(\Delta)$.
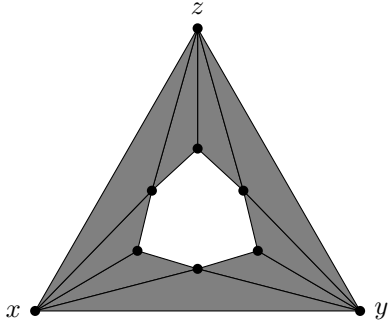
		\begin{center}
			\begin{tikzpicture}[scale = 0.8]
				\tikzstyle{point}=[circle,thick,draw=black,fill=black,inner sep=0pt,minimum width=3pt,minimum height=3pt]
				\node (x)[point, label=left:$x$] at (-1.7,-1) {};
				\node (y)[point, label=right:$y$] at (3.7,-1) {};
				\node (z)[point, label=above:$z$] at (1,3.7) {};
				
				\node (a)[point] at (0,0) {};
				\node (b)[point] at (2,0) {};
				\node (c)[point] at (1,1.7) {};
				
				\node (d)[point] at (1,-0.3) {}; 
				\node (e)[point] at (1.76,1) {}; 
				\node (f)[point] at (0.24,1) {}; 
				
				\begin{scope}[on background layer]
					\draw[fill=gray] (x.center) -- (a.center) -- (d.center) -- cycle;
					\draw[fill=gray] (x.center) -- (a.center) -- (f.center) -- cycle;
					\draw[fill=gray] (y.center) -- (b.center) -- (d.center) -- cycle;
					\draw[fill=gray] (y.center) -- (b.center) -- (e.center) -- cycle;
					\draw[fill=gray] (z.center) -- (c.center) -- (e.center) -- cycle;
					\draw[fill=gray] (z.center) -- (c.center) -- (f.center) -- cycle;
					
					\draw[fill=gray] (x.center) -- (y.center) -- (d.center) -- cycle;
					\draw[fill=gray] (y.center) -- (z.center) -- (e.center) -- cycle;
					\draw[fill=gray] (x.center) -- (z.center) -- (f.center) -- cycle;
				\end{scope}
			\end{tikzpicture}
		\end{center}
		Thus for $i=\dim \Delta +1$ we can view the operation $\phi_i$ as a kind of prism operator, with $\Delta$ at one end of the prism and $\BDelta$ at the other end (and for $i>\dim \Delta +1$, the complex $\phi_i(\Delta)$ is a prism with $\Delta$ at one end, and the cone $\BDelta \ast u_{\emptyset}$ at the other end).
		
		\subsection{Deformation Retractions}\label{Subsection: Phi-i Deformation Retractions}
		To find the homologies of the links in $\phi_i(\Delta)$, we will make use of some deformation retractions.
		
		In particular, we use Lemma \ref{Lemma: Deformation Retract} to obtain two deformation retractions of $\tDelta$: one ``\textit{vertex-first}'' deformation from $\tDelta$ on to $\Delta$, and one ``\textit{facet-first}'' deformation from $\tDelta$ on to $\BDelta$. The latter deformation only holds in the specific case where $i=\dim \Delta+1$. For each deformation we provide an example before detailing the general result. We begin with the vertex-first deformation.
		\begin{ex}\label{Example: vertex-first DR}
			Let $\Delta$ be the boundary of the $2$-simplex on vertex set $\{x,y,z\}$ as in Example \ref{Example: phi_i(Delta)}. We show that there is a deformation retraction $\phi_2(\Delta) \rightsquigarrow \Delta$.
			
			Note that every facet of $\phi_2(\Delta)$ which contains $u_{\{x\}}$ also contains $x$. Thus, if we set $g=\{u_{\{x\}}\}$ and $f = \{x,u_{\{x\}}\}$, Lemma \ref{Lemma: Deformation Retract} allows us to 	remove the vertex $u_{\{x\}}$ from $\phi_2(\Delta)$. Similarly we may remove the vertices $u_{\{y\}}$ and $u_{\{z\}}$. This gives us a deformation retraction $\phi_2(\Delta)\rightsquigarrow \phi_1(\Delta)$.
			
			The same reasoning now allows us to remove the vertices $u_{\{x,y\}}$, $u_{\{x,z\}}$ and $u_{\{y,z\}}$ from $\phi_1(\Delta)$ to obtain a deformation retraction $\phi_1(\Delta)\rightsquigarrow \Delta$.
			
			Diagrammatically, we have the deformation retractions:
			\begin{center}
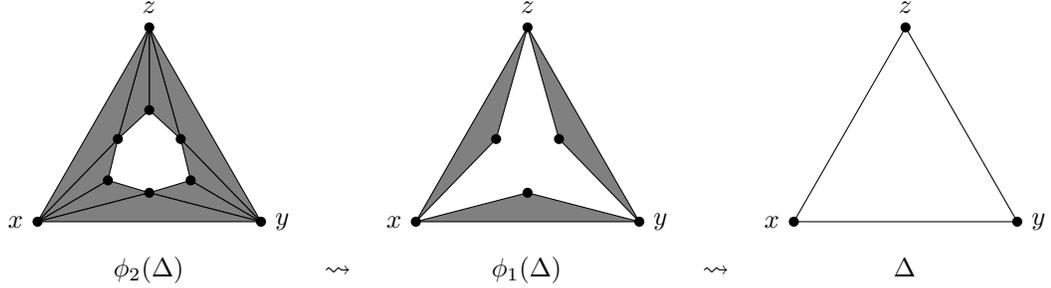
				\begin{tabular}{ c c c c c }
					\begin{tikzpicture}[scale = 0.55]
						\tikzstyle{point}=[circle,thick,draw=black,fill=black,inner sep=0pt,minimum width=3pt,minimum height=3pt]
						\node (x)[point, label=left:$x$] at (-1.7,-1) {};
						\node (y)[point, label=right:$y$] at (3.7,-1) {};
						\node (z)[point, label=above:$z$] at (1,3.7) {};
						
						\node (a)[point] at (0,0) {};
						\node (b)[point] at (2,0) {};
						\node (c)[point] at (1,1.7) {};
						
						\node (d)[point] at (1,-0.3) {}; 
						\node (e)[point] at (1.76,1) {}; 
						\node (f)[point] at (0.24,1) {}; 
						
						\begin{scope}[on background layer]
							\draw[fill=gray] (x.center) -- (a.center) -- (d.center) -- cycle;
							\draw[fill=gray] (x.center) -- (a.center) -- (f.center) -- cycle;
							\draw[fill=gray] (y.center) -- (b.center) -- (d.center) -- cycle;
							\draw[fill=gray] (y.center) -- (b.center) -- (e.center) -- cycle;
							\draw[fill=gray] (z.center) -- (c.center) -- (e.center) -- cycle;
							\draw[fill=gray] (z.center) -- (c.center) -- (f.center) -- cycle;
							
							\draw[fill=gray] (x.center) -- (y.center) -- (d.center) -- cycle;
							\draw[fill=gray] (y.center) -- (z.center) -- (e.center) -- cycle;
							\draw[fill=gray] (x.center) -- (z.center) -- (f.center) -- cycle;
						\end{scope}
					\end{tikzpicture} & 
					& \begin{tikzpicture}[scale = 0.55]
						\tikzstyle{point}=[circle,thick,draw=black,fill=black,inner sep=0pt,minimum width=3pt,minimum height=3pt]
						\node (x)[point, label=left:$x$] at (-1.7,-1) {};
						\node (y)[point, label=right:$y$] at (3.7,-1) {};
						\node (z)[point, label=above:$z$] at (1,3.7) {};
						
						
						\node (d)[point] at (1,-0.3) {}; 
						\node (e)[point] at (1.76,1) {}; 
						\node (f)[point] at (0.24,1) {}; 
						
						\begin{scope}[on background layer]
							
							\draw[fill=gray] (x.center) -- (y.center) -- (d.center) -- cycle;
							\draw[fill=gray] (y.center) -- (z.center) -- (e.center) -- cycle;
							\draw[fill=gray] (x.center) -- (z.center) -- (f.center) -- cycle;
						\end{scope}
					\end{tikzpicture} & 
					& \begin{tikzpicture}[scale = 0.55]
						\tikzstyle{point}=[circle,thick,draw=black,fill=black,inner sep=0pt,minimum width=3pt,minimum height=3pt]
						\node (x)[point, label=left:$x$] at (-1.7,-1) {};
						\node (y)[point, label=right:$y$] at (3.7,-1) {};
						\node (z)[point, label=above:$z$] at (1,3.7) {};
						
						\draw (x.center) -- (y.center) -- (z.center) -- cycle;
					\end{tikzpicture}\\
					$\phi_2(\Delta)$& $\rightsquigarrow$  & $\phi_1(\Delta)$ & $\rightsquigarrow$  & $\Delta$
				\end{tabular}
			\end{center}
			where each deformation retraction is obtained by 
			collapsing the vertices $u_\sigma$ for which $|\sigma|$ is minimal onto 			
			the face $\sigma$ in $\Delta$.
		\end{ex}
		
		This example generalises as follows:
		\begin{lem}\label{Lemma: DR Phi_i to Delta}
			Let $\Delta$ be a pure simplicial complex on vertex set $V$, $i\leq \dim \Delta + 1$, and $\tDelta=\phi_i(\Delta)$ on vertex set $V\sqcup S^i_\Delta$. There is a deformation retraction $\tDelta \rightsquigarrow \tDelta|_V = \Delta$.
		\end{lem}
		\begin{proof}
			Let $\sigma$ be a minimally sized face of $\Delta$ such that $u_\sigma$ is in $S^i_\Delta$. Note that $|\sigma|\geq \dim\Delta + 2 - i \geq \dim\Delta + 2 - (\dim\Delta + 1) = 1$, so in particular $\sigma$ is not empty.
			
			By construction, because $|\sigma|$ is minimal, every facet of $\Delta$ that contains $u_\sigma$ also contains $\sigma$. Thus, setting $g = \{u_\sigma\}$ and $f=\sigma \cup \{u_\sigma\}$, Lemma \ref{Lemma: Deformation Retract} gives us a deformation retraction $\tDelta \rightsquigarrow \tDelta - \{u_\sigma\}$.
			
			Continuing in this way we may remove every vertex $u_\sigma$ in $S^i_\Delta$ from $\tDelta$, in increasing order of the size of $\sigma$, and thus obtain a deformation retraction $\tDelta \rightsquigarrow \tDelta|_V = \Delta$.
		\end{proof}
		
		We now proceed to the facet-first deformation.
		\begin{ex}\label{Example: facet-first DR}
			Let $\Delta$ be the boundary of the $2$-simplex on vertex set $\{x,y,z\}$ as in Example \ref{Example: phi_i(Delta)}. We show that there is a deformation retraction $\phi_2(\Delta) \rightsquigarrow \BDelta$.
			
			Note that the edge $\{x,y\}$ of $\Delta$ occurs only in the facet $\{x,y,u_{\{x,y\}}\}$. Thus, if we set $g=\{x,y\}$ and $f = \{x,y,u_{\{x,y\}}\}$, Lemma \ref{Lemma: Deformation Retract} allows us to remove the edge $\{x,y\}$ from $\phi_2(\Delta)$. Similarly we may remove the edges $\{x,z\}$ and $\{y,z\}$.
			
			The same reasoning now allows us to remove the vertices $x$, $y$ and $z$ from $\phi_2(\Delta)$ to obtain a deformation retraction $\phi_2(\Delta)\rightsquigarrow \phi_2(\Delta)|_{S^2_\Delta}= \BDelta$.
			
			Diagrammatically, we have the deformation retractions:
			\begin{center}
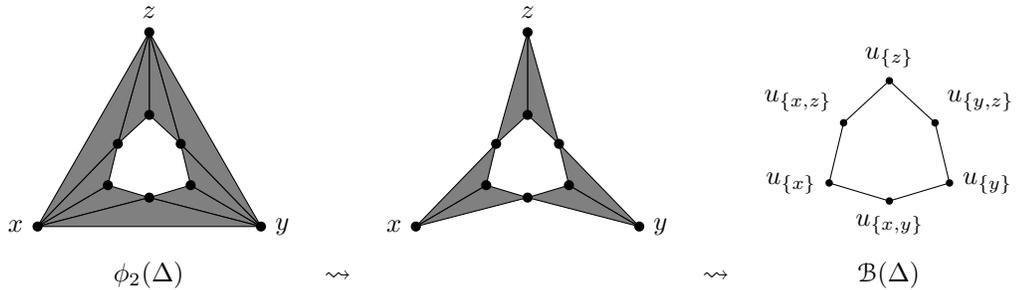
				\begin{tabular}{ c c c c c }
					\begin{tikzpicture}[scale = 0.55]
						\tikzstyle{point}=[circle,thick,draw=black,fill=black,inner sep=0pt,minimum width=3pt,minimum height=3pt]
						\node (x)[point, label=left:$x$] at (-1.7,-1) {};
						\node (y)[point, label=right:$y$] at (3.7,-1) {};
						\node (z)[point, label=above:$z$] at (1,3.7) {};
						
						\node (a)[point] at (0,0) {};
						\node (b)[point] at (2,0) {};
						\node (c)[point] at (1,1.7) {};
						
						\node (d)[point] at (1,-0.3) {}; 
						\node (e)[point] at (1.76,1) {}; 
						\node (f)[point] at (0.24,1) {}; 
						
						\begin{scope}[on background layer]
							\draw[fill=gray] (x.center) -- (a.center) -- (d.center) -- cycle;
							\draw[fill=gray] (x.center) -- (a.center) -- (f.center) -- cycle;
							\draw[fill=gray] (y.center) -- (b.center) -- (d.center) -- cycle;
							\draw[fill=gray] (y.center) -- (b.center) -- (e.center) -- cycle;
							\draw[fill=gray] (z.center) -- (c.center) -- (e.center) -- cycle;
							\draw[fill=gray] (z.center) -- (c.center) -- (f.center) -- cycle;
							
							\draw[fill=gray] (x.center) -- (y.center) -- (d.center) -- cycle;
							\draw[fill=gray] (y.center) -- (z.center) -- (e.center) -- cycle;
							\draw[fill=gray] (x.center) -- (z.center) -- (f.center) -- cycle;
						\end{scope}
					\end{tikzpicture}& 
					& 		\begin{tikzpicture}[scale = 0.55]
						\tikzstyle{point}=[circle,thick,draw=black,fill=black,inner sep=0pt,minimum width=3pt,minimum height=3pt]
						\node (x)[point, label=left:$x$] at (-1.7,-1) {};
						\node (y)[point, label=right:$y$] at (3.7,-1) {};
						\node (z)[point, label=above:$z$] at (1,3.7) {};
						
						\node (a)[point] at (0,0) {};
						\node (b)[point] at (2,0) {};
						\node (c)[point] at (1,1.7) {};
						
						\node (d)[point] at (1,-0.3) {}; 
						\node (e)[point] at (1.76,1) {}; 
						\node (f)[point] at (0.24,1) {}; 
						
						\begin{scope}[on background layer]
							\draw[fill=gray] (x.center) -- (a.center) -- (d.center) -- cycle;
							\draw[fill=gray] (x.center) -- (a.center) -- (f.center) -- cycle;
							\draw[fill=gray] (y.center) -- (b.center) -- (d.center) -- cycle;
							\draw[fill=gray] (y.center) -- (b.center) -- (e.center) -- cycle;
							\draw[fill=gray] (z.center) -- (c.center) -- (e.center) -- cycle;
							\draw[fill=gray] (z.center) -- (c.center) -- (f.center) -- cycle;
							
						\end{scope}
					\end{tikzpicture}  & 
					& 	\begin{tikzpicture}[scale = 0.8]
						\tikzstyle{point}=[circle,thick,draw=black,fill=black,inner sep=0pt,minimum width=2pt,minimum height=2pt]
						\node (a)[point, label=left:$u_{\{x\}}$] at (0,0) {};
						\node (b)[point, label=right:$u_{\{y\}}$] at (2,0) {};
						\node (c)[point, label=above:$u_{\{z\}}$] at (1,1.7) {};
						
						\node (d)[point, label=below:$u_{\{x,y\}}$] at (1,-0.3) {};
						\node (e)[point, label=above right:$u_{\{y,z\}}$] at (1.76,1) {};
						\node (f)[point, label=above left:$u_{\{x,z\}}$] at (0.24,1) {};
						
						\draw (a.center) -- (d.center) -- (b.center) -- (e.center) -- (c.center) -- (f.center) -- cycle;
					\end{tikzpicture}\\
					$\phi_2(\Delta)$& $\rightsquigarrow$  &  & $\rightsquigarrow$  & $\BDelta$
				\end{tabular}
			\end{center}
			where each deformation retraction is obtained by 
			collapsing the faces $\sigma$ of $\Delta$ for which $|\sigma|$ is maximal onto the vertex $u_\sigma$ in $\BDelta$.
		\end{ex}
		
		Once again this example admits a generalisation:
		\begin{lem}\label{Lemma: DR Phi_i to BDelta}
			Let $\Delta$ be a simplicial complex on vertex set $V$, $i = \dim \Delta + 1$, and $\tDelta=\phi_i(\Delta)$ on vertex set $V\sqcup S^i_\Delta$. There is a deformation retraction $\tDelta \rightsquigarrow \tDelta|_{S^i_\Delta} = \BDelta$.
		\end{lem}
		\begin{proof}
			The condition on the value of $i$ here means that $S^i_\Delta$ contains vertices $u_\sigma$ for every face $\sigma$ of $\Delta$ such that $|\sigma|\geq \dim \Delta + 2 - (\dim \Delta + 1) = 1$. In other words, \textit{every nonempty} face $\sigma$ of $\Delta$ has a corresponding vertex $u_\sigma$ in $S^i_\Delta$, but $u_\emptyset$ is not in $S^i_\Delta$. In particular this means that $\tDelta|_{S^i_\Delta} = \BDelta$.
			
			Let $\sigma$ be a maximally sized face of $\Delta$ (i.e. a facet of dimension $\dim \Delta$). By construction, because $|\sigma|$ is maximal, every facet of $\tDelta$ that contains $\sigma$ also contains $u_\sigma$. Thus, setting $g = \sigma$ and $f=\sigma \cup \{u_\sigma\}$, Lemma \ref{Lemma: Deformation Retract} gives us a deformation retraction $\tDelta \rightsquigarrow \tDelta - \sigma$.
			
			Because every nonempty face $\sigma$ of $\Delta$ has a corresponding vertex in $S^i_\Delta$, then we may continue in this way to remove every nonempty face $\sigma$ of $\Delta$ from $\tDelta$, in decreasing order of the size of $\sigma$, and thus obtain a deformation retraction $\tDelta \rightsquigarrow \tDelta|_{S^i_\Delta} = \BDelta$.
		\end{proof}
		
		The deformation in Lemma \ref{Lemma: DR Phi_i to BDelta} has the following important corollary.
		
		\begin{cor}\label{Corollary: DR Phi_i Acyclic}
			Let $\Delta$ be a simplicial complex on vertex set $V$, and $i > \dim \Delta + 1$. The complex $\tDelta=\phi_i(\Delta)$ is acyclic.
		\end{cor}
		\begin{proof}
			The condition on the value of $i$ here means that \textit{every} face $\sigma$ of $\Delta$ has a corresponding vertex $u_\sigma$ in $S^i_\Delta$, including $\sigma=\emptyset$. By Remark \ref{Remark: Phi_i = Phi_j} we also have that $\tDelta$ is equal to $\phi_m(\Delta)$ where $m = \dim \Delta + 2$.
			
			We can decompose $\tDelta$ into those facets $F$ which contain $u_\emptyset$ and those which do not. Note that we have $u_\emptyset \notin F$ if and only if $F$ is a facet of $\phi_{m-1}(\Delta)$; and $u_\emptyset \in F$ if and only if $F-\{u_\emptyset\}$ is a facet of $\BDelta$. Thus $\tDelta$ may be expressed as the union of $\phi_{m-1}(\Delta)$ and $\BDelta \ast \{u_\emptyset\}$, and these two subcomplexes intersect at $\BDelta$.
			
			By Lemma \ref{Lemma: DR Phi_i to BDelta} we have a deformation retraction $\phi_{m-1}(\Delta)\rightsquigarrow \BDelta$, which extends to a deformation retraction $\tDelta \rightsquigarrow \BDelta \ast \{u_\emptyset\}$. The complex $\BDelta \ast \{u_\emptyset\}$ is a cone over $u_\emptyset$, and is thus acyclic.
		\end{proof}
		
		\begin{ex}\label{Example: phi3 def retract}
			Once again, let $\Delta$ be the boundary of the $2$-simplex on vertex set $\{x,y,z\}$ as in Example \ref{Example: phi_i(Delta)}. Corollary \ref{Corollary: DR Phi_i Acyclic} gives us a deformation retraction $\phi_3(\Delta) \rightsquigarrow \BDelta$.
			\begin{center}
				\begin{tabular}{ c c c c c }
					\begin{tikzpicture}[scale = 0.6]
						\tikzstyle{point}=[circle,thick,draw=black,fill=black,inner sep=0pt,minimum width=3pt,minimum height=3pt]
						\node (x)[point, label=left:$x$] at (-1.7,-1) {};
						\node (y)[point, label=right:$y$] at (3.7,-1) {};
						\node (z)[point, label=above:$z$] at (1,3.7) {};
						
						\node (a)[point] at (0,0) {};
						\node (b)[point] at (2,0) {};
						\node (c)[point] at (1,1.7) {};
						
						\node (d)[point] at (1,-0.3) {}; 
						\node (e)[point] at (1.76,1) {}; 
						\node (f)[point] at (0.24,1) {}; 
						
						\node (g)[point] at (1,0.5) {};
						
						\begin{scope}[on background layer]
							\draw[fill=gray] (x.center) -- (a.center) -- (d.center) -- cycle;
							\draw[fill=gray] (x.center) -- (a.center) -- (f.center) -- cycle;
							\draw[fill=gray] (y.center) -- (b.center) -- (d.center) -- cycle;
							\draw[fill=gray] (y.center) -- (b.center) -- (e.center) -- cycle;
							\draw[fill=gray] (z.center) -- (c.center) -- (e.center) -- cycle;
							\draw[fill=gray] (z.center) -- (c.center) -- (f.center) -- cycle;
							
							\draw[fill=gray] (x.center) -- (y.center) -- (d.center) -- cycle;
							\draw[fill=gray] (y.center) -- (z.center) -- (e.center) -- cycle;
							\draw[fill=gray] (x.center) -- (z.center) -- (f.center) -- cycle;
							
							\draw[fill=gray] (a.center) -- (d.center) -- (g.center) -- cycle;
							\draw[fill=gray] (a.center) -- (f.center) -- (g.center) -- cycle;
							\draw[fill=gray] (b.center) -- (d.center) -- (g.center) -- cycle;
							\draw[fill=gray] (b.center) -- (e.center) -- (g.center) -- cycle;
							\draw[fill=gray] (c.center) -- (e.center) -- (g.center) -- cycle;
							\draw[fill=gray] (c.center) -- (f.center) -- (g.center) -- cycle;
						\end{scope}
					\end{tikzpicture} && 
					\begin{tikzpicture}[scale = 0.8]
						\tikzstyle{point}=[circle,thick,draw=black,fill=black,inner sep=0pt,minimum width=3pt,minimum height=3pt]
						\node (x) at (-1.7,-1) {};
						\node (y) at (3.7,-1) {};
						\node (z) at (1,3.7) {};
						
						\node (a)[point,label=left:$u_{\{x\}}$] at (0,0) {};
						\node (b)[point,label=right:$u_{\{y\}}$] at (2,0) {};
						\node (c)[point,label=above:$u_{\{z\}}$] at (1,1.7) {};
						
						\node (d)[point,label=below:$u_{\{x,y\}}$] at (1,-0.3) {}; 
						\node (e)[point,label=right:$u_{\{y,z\}}$] at (1.76,1) {}; 
						\node (f)[point,label=left:$u_{\{x,z\}}$] at (0.24,1) {}; 
						
						\node (g)[point] at (1,0.5) {};
						
						\begin{scope}[on background layer]
							\draw[fill=gray] (a.center) -- (d.center) -- (g.center) -- cycle;
							\draw[fill=gray] (a.center) -- (f.center) -- (g.center) -- cycle;
							\draw[fill=gray] (b.center) -- (d.center) -- (g.center) -- cycle;
							\draw[fill=gray] (b.center) -- (e.center) -- (g.center) -- cycle;
							\draw[fill=gray] (c.center) -- (e.center) -- (g.center) -- cycle;
							\draw[fill=gray] (c.center) -- (f.center) -- (g.center) -- cycle;
						\end{scope}
					\end{tikzpicture}
					&&  \begin{tikzpicture}[scale = 0.8]
						\tikzstyle{point}=[circle,thick,draw=black,fill=black,inner sep=0pt,minimum width=3pt,minimum height=3pt]
						\node (x) at (-1.5,-1) {};
						\node (y) at (3.7,-1) {};
						\node (z) at (1,3.7) {};
						
						\node (g)[point,label=below:$u_{\emptyset}$] at (1,0.5) {};
					\end{tikzpicture}\\
					$\phi_3(\Delta)$& $\rightsquigarrow$  & $\BDelta\ast \{u_\emptyset\}$ & $\rightsquigarrow$  & $\{u_\emptyset\}$
				\end{tabular}
			\end{center}
		\end{ex}

		\begin{lem}\label{Lemma: DR Phi_i link of new vertices acyclic}
			Let $\Delta$ be a simplicial complex on vertex set $V$, $i\leq \dim \Delta + 1$, and $\tDelta=\phi_i(\Delta)$ on vertex set $V\sqcup S^i_\Delta$. For any nonempty face $\sigma$ of $\tDelta$ contained entirely in $S^i_\Delta$, the complex $\link_{\tDelta} \sigma$ is acyclic.
		\end{lem}
		\begin{proof}
			Suppose $\sigma = \{u_{\tau_1},\dots,u_{\tau_r}\}$ for some faces $\tau_1 \subset \dots \subset \tau_r$ of $\Delta$. The condition on $i$ implies that $\tau_1 \neq \emptyset$.
			
			We must have $\tau_1 \in \link_{\tDelta} \sigma$ because $\tau_1 \sqcup \{u_{\tau_1},\dots,u_{\tau_r}\}$ is a face of $\tDelta$. Also, no vertex in $V-\tau_1$ can be contained in $\link_{\tDelta} \sigma$, because by construction, for any facet $G$ of $\tDelta$ containing $u_{\tau_1}$, we have $G\cap V \subseteq \tau_1$. We claim that there is a deformation retraction $\link_{\tDelta} \sigma \rightsquigarrow \link_{\tDelta} \sigma|_V = \langle\tau_1\rangle$, which is acyclic.
			
			We start by removing every vertex $u_\rho$ in $\link_{\tDelta} \sigma$ and $S^i_\Delta$ for which $\rho \subset \tau_1$. Suppose $u_\rho$ is any such vertex with $|\rho|$ minimal. By the minimality of $|\rho|$, we know that every facet of $\tDelta$ containing $u_\rho$ must also contain $\rho$. Thus every facet of $\link_{\tDelta} \sigma$ containing $u_\rho$ must also contain $\rho$. Setting $g=\{u_\rho\}$ and $f = \{u_\rho\} \cup \rho$, Lemma \ref{Lemma: Deformation Retract} gives us a deformation retraction $\link_{\tDelta} \sigma \rightsquigarrow \link_{\tDelta} \sigma - \{u_\rho\}$. Continuing in this way we may remove every vertex $u_\rho$ in $S^i_\Delta$ from $\link_{\tDelta} \sigma$, in increasing order of the size of $\rho$.
			
			Now we remove the vertices $u_\rho$ in $\link_{\tDelta} \sigma$ and $S^i_\Delta$ for which $\rho \supset \tau_1$. Suppose $u_\rho$ is any such vertex. Because $\rho$ contains $\tau_1$, every facet of $\tDelta$ containing $u_\rho$ and $\sigma$ must also contain $\tau_1$. Thus every facet of $\link_{\tDelta} \sigma$ containing $u_\rho$ must also contain $\tau_1$. Setting $g=\{u_\rho\}$ and $f = \{u_\rho\} \cup \tau_1$, Lemma \ref{Lemma: Deformation Retract} allows us to remove $u_\rho$ from $\link_{\tDelta} \sigma$.
		\end{proof}
		
		\subsection{Links in $\BDelta$}\label{Subsection: Links in BDelta}
		Let $\Delta$ be a simplicial complex and $i > \dim \Delta + 1$, and set $\tDelta = \phi_i(\Delta)$. The condition on $i$ ensures that $u_\emptyset$ is a vertex in $S^i_\Delta$.
		
		In this section, we examine the links of those faces $\sigma$ of $\tDelta$ for which $\emptyset \neq \sigma \subseteq S^i_\Delta$ (i.e. the nonempty faces in the induced subcomplex $\tDelta|_{S^i_\Delta}$). We begin by showing that for any such $\sigma$, we can express the homology of $\link_{\tDelta} \sigma$ in terms of the homology of $\link_{\tDelta} (\sigma\cup \{u_\emptyset\})$. As we will explain, this allows us to restrict our attention to the links in the barycentric subdivision complex $\BDelta$.
		
		\begin{prop}\label{Proposition: link of new vertices may as well contain u-emptyset}
			Let $\Delta$ be a simplicial complex and $i > \dim \Delta + 1$, and set $\tDelta = \phi_i(\Delta)$. Suppose $\emptyset \neq \sigma \subseteq S^i_\Delta$ is a face of $\tDelta$ with $u_{\emptyset}\notin \sigma$. For every $j\geq -1$ we have an isomorphism $\Hred_j(\link_{\tDelta} \sigma)\cong \Hred_{j-1}(\link_{\tDelta} (\sigma\cup \{u_{\emptyset}\}))$.
		\end{prop}
		\begin{proof}
			As in the proof of Corollary \ref{Corollary: DR Phi_i Acyclic}, we may decompose $\link_{\tDelta} \sigma$ into a subcomplex $A$ consisting of facets which contain $u_\emptyset$ and a subcomplex $B$ consisting of those which do not. The intersection of these subcomplexes consists of those faces $f$ in $\link_{\tDelta} \sigma$ for which $u_{\emptyset}$ is not in $f$ but $f\sqcup \{u_\emptyset\}$ is a face of $\tDelta$. In other words, we have $A\cap B = \link_{\tDelta} (\sigma\cup \{u_{\emptyset}\})$.
			
			The subcomplex $A$ is a cone over $u_\emptyset$, and is therefore acyclic. For every face $f$ in $B$, the intersection of $f$ and $S^i_{\Delta}$ contains only vertices $u_\tau$ for which $\tau$ is nonempty. All of these are faces of $\phi_m(\Delta)$ where $m=\dim\Delta +1$, and hence we have $B=\link_{\phi_m(\Delta)} \sigma$, which is also acyclic by Lemma \ref{Lemma: DR Phi_i link of new vertices acyclic}.
			
			Thus for every $j\geq -1$, the Mayer-Vietoris Sequence gives us an exact sequence
			\begin{equation*}
				0 \rightarrow \Hred_j(\link_{\tDelta} \sigma)\rightarrow \Hred_{j-1}(\link_{\tDelta} (\sigma\cup \{u_{\emptyset}\}))\rightarrow 0 
			\end{equation*}
			as required.
		\end{proof}
		
		Proposition \ref{Proposition: link of new vertices may as well contain u-emptyset} allows us to restrict our attention to those faces of $\tDelta|_{S^i_\Delta}$ which contain $u_\emptyset$. Note we have $\link_{\tDelta} u_\emptyset = \BDelta$, and hence the link of any face of $\tDelta$ which contains $u_\emptyset$ must be a link in $\BDelta$. For this reason, we devote the rest of this section to investigating the links of $\BDelta$.
		
		\begin{prop}\label{Proposition: Links in BDelta}
			Let $\sigma = \{u_{\tau_1},\dots, u_{\tau_r}\}\in \BDelta$ for some faces $\tau_1 \subset \dots \subset \tau_r$ of $\Delta$. We have an isomorphism of complexes
			\begin{equation*}
				\link_{\BDelta}\sigma \cong \mathcal{B}(\lkds[\tau_r])\ast \mathcal{B}(\link_{\partial \tau_r} \tau_{r-1})\ast \dots \ast \mathcal{B}(\link_{\partial \tau_2} \tau_1)\ast \mathcal{B}(\partial \tau_1).
			\end{equation*}
		\end{prop}
		\begin{proof}
			
			For notational convenience, we set $\tau_0 = \emptyset$, so that $\mathcal{B}(\partial \tau_1)$ may be rewritten as $\mathcal{B}(\link_{\partial \tau_1} \tau_0)$.
			
			Let $A_r$ denote the induced subcomplex of $\BDelta$ on vertices of the form $u_f$ where $f$ contains $\tau_r$; and for each $0\leq j \leq r-1$, let $A_j$ denote the induced subcomplex of $\BDelta$ on vertices of the form $u_f$ for which we have $\tau_j \subset f\subset \tau_{j+1}$. Note that $A_r,\dots,A_0$ are pairwise disjoint subcomplexes of $\BDelta$.
			
			We claim that $\link_{\BDelta}\sigma = A_r\ast A_{r-1}\ast \dots\ast A_0$. This is sufficient to prove our proposition because the complex $A_r$ is isomorphic to $\mathcal{B}(\lkds[\tau_r])$ via the vertex map $u_f\mapsto u_{f-\tau_r}$, and for each $0\leq j \leq r-1$, the complex $A_j$ is isomorphic to $\mathcal{B}(\link_{\partial \tau_{j+1}} \tau_j)$ via the vertex map $u_f\mapsto u_{f-\tau_j}$.
			
			Let $\rho = \rho_r\sqcup \dots \sqcup \rho_0$ be a face of $A_r\ast \dots\ast A_0$, with $\rho_j\in A_j$ for each $0\leq j \leq r$. For each $0\leq j \leq r-1$, the face $\rho_j$ must be of the form $\{u_{f_1},\dots,u_{f_m}\}$ for some sequence of faces $\tau_j \subset f_1\subset \dots \subset f_m\subset \tau_{j+1}$ of $\Delta$. Similarly the face $\rho_r$ must be of the form $\{u_{f_1},\dots,u_{f_m}\}$ for some sequence of faces $\tau_r \subset f_1\subset \dots \subset f_m$ of $\Delta$. In particular, none of the vertices $u_{\tau_1},\dots,u_{\tau_r}$ are contained in $\rho$, so we have $\rho \cap \sigma = \emptyset$. Moreover, the faces of $\Delta$ corresponding to vertices in $\rho\sqcup \sigma$ may be arranged in a strict sequence by inclusion, which means that $\rho\sqcup \sigma\in \BDelta$. Thus $\rho$ is a face of $\link_{\BDelta}\sigma$.
			
			Conversely, suppose $\rho=\{u_{f_1},\dots,u_{f_m}\}$ is any face of $\link_{\BDelta} \sigma$. We may decompose $\rho$ into the disjoint union $\rho_r\sqcup \dots\sqcup \rho_0$, where for each $0\leq j\leq r-1$  the face $\rho_j$ contains all the vertices $u_f$ in $\rho$ for which we have $\tau_j\subset f \subset \tau_{j+1}$, making $\rho_j$ a face of $A_j$; and the face $\rho_r$ contains all the vertices $u_f$ in $\rho$ for which $f$ contains $\tau_r$, making $\rho_r$ a face of $A_r$. Thus $\rho$ is a face of $A_r\ast \dots \ast A_0$.
		\end{proof}
		
		In particular, Proposition \ref{Proposition: Links in BDelta} has the following important corollaries.
		\begin{cor}\label{Corollary: Homology of Links in BDelta}
			Let $\Delta$ and $\sigma=\{u_{\tau_1},\dots,u_{\tau_r}\}\in \BDelta$ be as in Proposition \ref{Proposition: Links in BDelta}. We have $h(\BDelta,\sigma)=h(\Delta,\tau_r)\boxplus\{|\tau_r|-r\}$ (where $\boxplus$ denotes the addition of sets operation laid out in Definition \ref{Definition: Boxplus}).
		\end{cor}
		\begin{proof}
			Just as in the proof of Proposition \ref{Proposition: Links in BDelta}, we set $\tau_0=\emptyset$ for notational convenience. Using this proposition and Proposition \ref{Proposition: Kunneth Formula}, we can compute the homology of $\link_{\BDelta}\sigma$ from the homologies of $\calB(\link_\Delta \tau_r), \mathcal{B}(\link_{\partial \tau_r} \tau_{r-1}), \dots, \mathcal{B}(\link_{\partial \tau_1} \tau_0)$.
			
			To compute these homologies, first recall that for any integer $i\geq -1$ we have
			\begin{equation}\label{Equation: h(BDelta)=h(Delta)}
				\Hred_i(\BDelta)=\Hred_i(\Delta).
			\end{equation}
			Next, note that for each $0\leq j \leq r-1$, the complex $\partial \tau_{j+1}$ is the boundary of the $(|\tau_{j+1}|-1)$-simplex. As observed in Example \ref{Example: PR Simplex}, the link of any face $\tau$ in the boundary of the $p$-simplex is the boundary of the $(p-|\tau|)$-simplex, which has homology only at degree $p-|\tau|-1$. Thus we have
			\begin{equation}\label{Equation: h(link of tau_j)}
				h(\link_{\partial \tau_{j+1}}\tau_j)=\{|\tau_{j+1}|-|\tau_j|-2\}.
			\end{equation}
			Putting these results together, we find
			\begin{align*}
				h(\BDelta, \sigma) &= h(\mathcal{B}(\lkds[\tau_r])\ast \mathcal{B}(\link_{\partial \tau_1}\tau_0)\ast \dots \ast \mathcal{B}(\link_{\partial \tau_r}\tau_{r-1})) &\text{by Prop. \ref{Proposition: Links in BDelta}}\\
				&= h(\mathcal{B}(\lkds[\tau_r]))\boxplus \sum_{j=0}^{r-1} h(\mathcal{B}(\link_{\partial \tau_{j+1}}\tau_j))\boxplus \{r\} &\text{by Prop. \ref{Proposition: Kunneth Formula}}\\
				&= h(\lkds[\tau_r])\boxplus \sum_{j=0}^{r-1}h(\link_{\partial \tau_{j+1}}\tau_j) \boxplus \{r\} &\text{by Eq. \ref{Equation: h(BDelta)=h(Delta)}}\\
				&=h(\Delta,\tau_r) \boxplus \sum_{j=0}^{r-1} \{|\tau_{j+1}|-|\tau_j|-2\}\boxplus\{r\}&\text{by Eq. \ref{Equation: h(link of tau_j)}}\\
				&=h(\Delta, \tau_r)\boxplus\{|\tau_r|-r\}.&
			\end{align*}
		\end{proof}
		\begin{cor}\label{Corollary: Homology of Links in BDelta PR (1,..,1)}
			Let $\Delta$ be a PR complex with degree type $(1,\dots,1)$. For any $\sigma\in \BDelta$ we have that $h(\BDelta,\sigma)$ is either empty or equal to $\{\dim \Delta-|\sigma|\}$	
			
		\end{cor}
		\begin{proof}
			For $\sigma=\emptyset$ we have $h(\BDelta,\emptyset)=h(\BDelta)=h(\Delta)$, because $\BDelta$ is homeomorphic to $\Delta$. Thus $h(\BDelta,\emptyset)$ is empty, unless $\Delta$ has homology, in which case it is equal to $\{\dim \Delta\}$ by Corollary \ref{Corollary: PR (1...1) Definition}.
			
			Now assume $\sigma=\{u_{\tau_1},\dots,u_{\tau_r}\}$ for some $\tau_1\subset \dots \subset \tau_r$ in $\Delta$. From Corollary \ref{Corollary: Homology of Links in BDelta}, we know that $h(\BDelta,\sigma)=h(\Delta,\tau_r)\boxplus\{|\tau_r|-|\sigma|\}$. If $h(\Delta,\tau_r)$ is empty (i.e. $\lkds[\tau_r]$ is acyclic) then this sum is empty. Otherwise, by Corollary \ref{Corollary: PR (1...1) Definition} we have $h(\Delta,\tau_r)=\{\dim \Delta-|\tau_r|\}$, and hence $h(\BDelta,\sigma)=\{\dim \Delta-|\tau_r|\}\boxplus\{|\tau_r|-|\sigma|\}=\{\dim \Delta-|\sigma|\}$.
		\end{proof}
		
		It follows from Corollary \ref{Corollary: Homology of Links in BDelta PR (1,..,1)} that if $\Delta$ is Cohen-Macaulay (i.e. PR with degree type $(1,\dots,1)$), then $\BDelta$ is also Cohen-Macaulay. In fact, this turns out to be the \textit{only} condition under which $\BDelta$ is PR, as the following proposition demonstrates. This proposition will not be strictly necessary for our proof of Theorem \ref{Theorem: Phi_i Operations Degree Type}, but it helps to illuminate why the operation $\phi_i$ preserves the PR property for PR complexes of degree type $(d_p,\dots,d_i,1,\dots,1)$, and why it fails to do so for PR complexes of other degree types.

		\begin{prop}\label{Proposition: When BDelta is PR}
			Let $\Delta$ be a simplicial complex. The following are equivalent.
			\begin{enumerate}
				\item $\BDelta$ is PR.
				\item  $\BDelta$ is Cohen-Macaulay (i.e. PR with degree type $(1,\dots,1)$).
				\item $\Delta$ is Cohen-Macaulay (i.e. PR with degree type $(1,\dots,1)$).
			\end{enumerate}
		\end{prop}
		\begin{proof}
			(3)$\Rightarrow$(2) follows from Corollary \ref{Corollary: Homology of Links in BDelta PR (1,..,1)}, and (2)$\Rightarrow$(1) is immediate. To prove (1)$\Rightarrow$(3), we show the contrapositive.
			
			First assume that $\Delta$ is not PR. This means that $\Delta$ has two faces $\tau_1$ and $\tau_2$ of different sizes such that the intersection $h(\Delta,\tau_1)\cap h(\Delta,\tau_2)$ is nonempty. Suppose $\iota$ is an index in both $h(\Delta,\tau_1)$ and $h(\Delta,\tau_2)$. By Corollary \ref{Corollary: Homology of Links in BDelta}, we have $h(\BDelta,u_{\tau_1})=h(\Delta,\tau_1)\boxplus\{|\tau_1|-1\}$ and $h(\BDelta,u_{\tau_2})=h(\Delta,\tau_2)\boxplus\{|\tau_2|-1\}$. Thus the complete homology index set $\hh(\BDelta,1)$ contains both $\iota+|\tau_1|-1$ and $\iota+|\tau_2|-1$, and therefore cannot be a singleton. This means $\BDelta$ is not PR by Corollary \ref{Corollary: PR Complex links have single homology}.
			
			Now assume that $\Delta$ is PR of degree type $(d_p,\dots,d_1)$ where $d_j > 1$ for some $1\leq j\leq p$. By Proposition \ref{Proposition: Alternate PR Definition With Degree Type}, $\Delta$ must have two faces $\tau_1$ and $\tau_2$ such that $|\tau_2|=|\tau_1|+d_j$ and $h(\Delta,\tau_1)=\{\iota\}$ while $h(\Delta,\tau_2)=\{\iota-1\}$ for some index $\iota$. By Corollary \ref{Corollary: Homology of Links in BDelta}, we have $h(\BDelta,u_{\tau_1})=\{\iota\}\boxplus\{|\tau_1|-1\}=\{\iota+|\tau_1|-1\}$, and $h(\BDelta,u_{\tau_2})=\{\iota-1\}\boxplus\{|\tau_2|-1\}=\{\iota + |\tau_1|+d_j-2\}$. In particular, because $d_j>1$, we have $d_j-2>-1$ and hence these two sets are not equal. Thus, once again, the complete homology index set $\hh(\BDelta,1)$ is not a singleton, and so $\BDelta$ cannot be PR by Corollary \ref{Corollary: PR Complex links have single homology}.
		\end{proof}
		
		\subsection{Proving Theorem \ref{Theorem: Phi_i Operations Degree Type}}\label{Subsection: Proving Phi_i Theorem}
		In this section, we bring together all of the results of the previous sections to prove Theorem \ref{Theorem: Phi_i Operations Degree Type} (and hence Theorem \ref{Theorem: PR Complexes of Any Degree Type}).
		
		We begin with an example to elucidate our method of proof
		\begin{ex}
			Let $\Delta$ be the boundary of the $2$-complex in Example \ref{Example: phi_i(Delta)}.
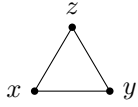
			\begin{center}.
				\begin{tikzpicture}[scale = 0.5]
					\tikzstyle{point}=[circle,thick,draw=black,fill=black,inner sep=0pt,minimum width=2pt,minimum height=2pt]
					\node (a)[point, label=left:$x$] at (0,0) {};
					\node (b)[point, label=right:$y$] at (2,0) {};
					\node (c)[point, label=above:$z$] at (1,1.7) {};
					
					\draw (a.center) -- (b.center) -- (c.center) -- cycle;
				\end{tikzpicture}
			\end{center} 
			We've seen that $\Delta$ is PR with degree type $(1,1)$. From Example \ref{Example: PR Simplex}, we know that the Betti diagram
			$\beta(\Idstar)$ is $\begin{array}{l | rrr }
				1 & 3 & 3 & 1  
			\end{array}$ and by ADHF, these entries come from the links of faces of the following sizes.
			
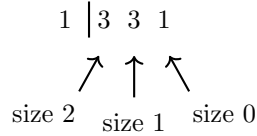
\begin{center}
				\begin{tikzpicture}[scale = 0.4]
					\node (x)[label=above:$1$] at (-1.3,3) {};
					\draw [line width=0.8pt] (-.5,4.5) -- (-.5,3.3);
					\node (a)[label=above:$3$] at (0,3) {};
					\node (b)[label=above:$3$] at (1,3) {};
					\node (c)[label=above:$1$] at (2,3) {};
					
					\node (a') at (-1,1.2) {};
					\node (b') at (1,1) {};
					\node (c') at (3,1.2) {};
					
					\node (a'') [label=below:size 2] at (-2,1.7) {};
					\node (b'') [label=below:size 1] at (1,1.4) {};
					\node (c'') [label=below:size 0] at (4,1.7) {};
					\draw[->, line width=0.8pt] (a')   -- (a);
					\draw[->, line width=0.8pt] (b')   -- (b);
					\draw[->, line width=0.8pt] (c')   -- (c);
				\end{tikzpicture}
			\end{center}
			Put another way, the only nonempty complete homology index sets of $\Delta$ are the ones given in the table below.
			\begin{center}
				\begin{tabular}{ c | l }
					$\hh(\Delta,m)$ & $m$\\
					\hline
					$\{1\}$ & $0$\\
					$\{0\}$ & $1$\\
					$\{-1\}$ & $2$
				\end{tabular}
			\end{center}
			
			Now consider the complex $\tDelta=\phi_2(\Delta)$.
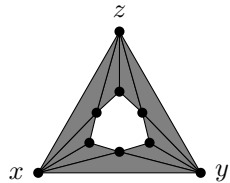
			\begin{center}
				\begin{tikzpicture}[scale = 0.4]
					\tikzstyle{point}=[circle,thick,draw=black,fill=black,inner sep=0pt,minimum width=3pt,minimum height=3pt]
					\node (x)[point, label=left:$x$] at (-1.7,-1) {};
					\node (y)[point, label=right:$y$] at (3.7,-1) {};
					\node (z)[point, label=above:$z$] at (1,3.7) {};
					
					\node (a)[point] at (0,0) {}; 
					\node (b)[point] at (2,0) {}; 
					\node (c)[point] at (1,1.7) {}; 
					
					\node (d)[point] at (1,-0.3) {}; 
					\node (e)[point] at (1.76,1) {}; 
					\node (f)[point] at (0.24,1) {}; 
					
					

					
					
					\begin{scope}[on background layer]
						\draw[fill=gray,fill opacity=.8] (x.center) -- (a.center) -- (d.center) -- cycle;
						\draw[fill=gray,fill opacity=.8] (x.center) -- (a.center) -- (f.center) -- cycle;
						\draw[fill=gray,fill opacity=.8] (y.center) -- (b.center) -- (d.center) -- cycle;
						\draw[fill=gray,fill opacity=.8] (y.center) -- (b.center) -- (e.center) -- cycle;
						\draw[fill=gray,fill opacity=.8] (z.center) -- (c.center) -- (e.center) -- cycle;
						\draw[fill=gray,fill opacity=.8] (z.center) -- (c.center) -- (f.center) -- cycle;
						
						\draw[fill=gray,fill opacity=.8] (x.center) -- (y.center) -- (d.center) -- cycle;
						\draw[fill=gray,fill opacity=.8] (y.center) -- (z.center) -- (e.center) -- cycle;
						\draw[fill=gray,fill opacity=.8] (x.center) -- (z.center) -- (f.center) -- cycle;
						
					\end{scope}
					
				\end{tikzpicture}
			\end{center}
			We want to show that $\tDelta$ has degree type $(2,1)$. In other words, we want to show that its Betti diagram has the form $\begin{array}{c | rrr }
				i & * & * & .\\
				i+1 & . & . & *
			\end{array}$ where the nonzero entries $*$ come from the links of faces of the following sizes
			
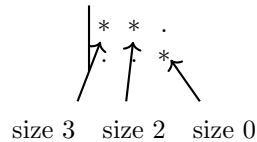
\begin{center}
				\begin{tikzpicture}[scale = 0.4]
					\draw [line width=0.8pt] (-.5,5.5) -- (-.5,3.3);
					\node (a)[label=above:$*$] at (0,4) {};
					\node (b)[label=above:$*$] at (1,4) {};
					\node (c)[label=above:$*$] at (2,3) {};
					
					\node (a1)[label=above:$.$] at (0,3) {};
					\node (b1)[label=above:$.$] at (1,3) {};
					\node (c1)[label=above:$.$] at (2,4) {};
					\node (a2) at (0,4.6) {};
					\node (b2) at (1,4.6) {};
					
					\node (a') at (-1,2) {};
					\node (b') at (0.7,2) {};
					\node (c') at (3.4,2) {};
					
					\node (a'') [label=below:size 3] at (-2,2.3) {};
					\node (b'') [label=below:size 2] at (1,2.3) {};
					\node (c'') [label=below:size 0] at (4,2.3) {};
					
					\draw[->, line width=0.8pt] (a')   -- (a2);
					\draw[->, line width=0.8pt] (b')   -- (b2);
					\draw[->, line width=0.8pt] (c')   -- (c1);
				\end{tikzpicture}
			\end{center}
			Put another way, we want to show that the only nonempty complete homology index sets of $\tDelta$ are the ones given in the table below.
			\begin{center}
				\begin{tabular}{ c | l }
					$\hh(\tDelta,m)$ & $m$\\
					\hline
					$\{1\}$ & $0$\\
					$\{0\}$ & $2$\\
					$\{-1\}$ & $3$
				\end{tabular}
			\end{center}
			
			So it suffices to show that, for any natural number $0\leq m \leq 3$ we have
			\begin{equation*}
				\hh(\tDelta,m) = \begin{cases}
					\hh(\Delta,m)& \text{ if } m < 1\\
					\emptyset & \text{ if } m = 1\\
					\hh(\Delta,m-1)& \text{ if } m > 1.
				\end{cases}
			\end{equation*}
			and for these two complexes, we are able to check this manually (by finding the link of each face of $\Delta$ and $\tDelta$).
			
			The proof below is a generalisation of this argument.
			
		\end{ex}
		
		\begin{proof}[Proof of Theorem \ref{Theorem: Phi_i Operations Degree Type}]
			Let $\Delta$ be a PR complex with degree type $(d_p,\dots,d_i,\underbrace{1,\dots,1}_{i-1})$ for some $p\geq i \geq 1$, and offset $s$, and set $\tDelta =\phi_i(\Delta)$. We aim to prove that $\tDelta$ is PR with degree type $(d_p,\dots,d_i+1,1,\dots,1)$.
			
			By Proposition \ref{Proposition: Alternate PR Definition With Degree Type}, the only nonempty complete homology index sets are the ones given in Table 6.1 below. By the same proposition, we can show that $\tDelta$ is PR with the desired degree type by proving that its only nonempty complete homology index sets are the ones given in Table 6.2.

			\begin{center}
				\begin{tabular}{c c}
					\textbf{Table 6.1:} \textit{Homology} & \textbf{Table 6.2:} \textit{Required Homology}\\
					\textit{Index Sets of} $\Delta$ &\textit{Index Sets of} $\tDelta$\\
					\begin{tabular}{ c | l }
						$\hh(\Delta,m)$ & $m$\\
						\hline
						$\{p-1\}$ & $s$\\
						$\{p-2\}$ & $s+d_p$\\
						$\vdots$ & $\vdots$\\
						$\{i-1\}$ & $s+\sum_{j=i+1}^p d_j$ \\
						$\{i-2\}$ & $s+\sum_{j=i}^p d_j$ \\
						$\{i-3\}$ & $s+\sum_{j=i}^p d_j + 1$ \\
						$\vdots$ & $\vdots$ \\
						$\{0\}$ & $s+\sum_{j=i}^p d_j + (i-2)$ \\
						$\{-1\}$ & $s+\sum_{j=i}^p d_j + (i-1)$\\
					\end{tabular}
					&
					\begin{tabular}{ c | l }
						$\hh(\tDelta,m)$ & $m$\\
						\hline
						$\{p-1\}$ & $s$\\
						$\{p-2\}$ & $s+d_p$\\
						$\vdots$ & $\vdots$\\
						$\{i-1\}$ & $s+\sum_{j=i+1}^p d_j$ \\
						$\{i-2\}$ & $s+\sum_{j=i}^p d_j + 1$ \\
						$\{i-3\}$ & $s+\sum_{j=i}^p d_j + 2$ \\
						$\vdots$ & $\vdots$ \\
						$\{0\}$ & $s+\sum_{j=i}^p d_j + (i-1)$ \\
						$\{-1\}$ & $s+\sum_{j=i}^p d_j + i$\\
					\end{tabular}
				\end{tabular}
			\end{center}
			Thus it suffices to show that for any natural number $0\leq m \leq s + \sum_{j=i}^p d_j + i$, we have
			\begin{equation}\label{Equation: homology index sets}
				\hh(\tDelta,m) = \begin{cases}
					\hh(\Delta,m)& \text{ if } m < s + \sum_{j=i}^p d_j\\
					\emptyset & \text{ if } m = s + \sum_{j=i}^p d_j\\
					\hh(\Delta,m-1)& \text{ if } m > s + \sum_{j=i}^p d_j.
				\end{cases}
			\end{equation}
			
			To this end, we fix a face $\sigma$ of $\tDelta$ of size $|\sigma|=m$, and investigate the homology index set $h(\tDelta,\sigma)$.
			
			We start by decomposing $\sigma$ into $\sigma = \sigma_V \sqcup \sigma_S$, where $\sigma_V$ is a subset of $V$ (and is thus in $\Delta$) and $\sigma_S$ is a subset of $S^i_\Delta$ (and is thus in $\BDelta$).
			
			Note that $\link_{\tDelta} \sigma = \link_{\link_{\tDelta} \sigma_V} \sigma_S$. Thus, using the isomorphism in Lemma \ref{Lemma: Phi_i and Link commute}, we may view $\link_{\tDelta} \sigma$ as a link in the complex $\phi_i(\lkds[\sigma_V])$, and then apply the results of Sections \ref{Subsection: Phi-i Deformation Retractions} and \ref{Subsection: Links in BDelta} to this link.
			
			To work out which results to apply, we will need to determine whether $i$ is less than, equal to, or greater than $\dim(\lkds[\sigma_V])+1$. By Proposition \ref{Proposition: Dimension of PR Complexes} we know $\dim \Delta = s + \sum_{j=i}^p d_j + i - 2$, which means we have
			\begin{equation}\label{Equation: dim lkds[sigma_V]}
				\dim (\lkds[\sigma_V])+1= s+ \sum_{j=i}^p d_j + i - 1 - |\sigma_V|.
			\end{equation}
			
			We now examine a number of different cases.
			
			\begin{itemize}
				\item \textbf{Case 1:} Assume $\sigma_S = \emptyset$ (i.e. $\sigma \in \Delta$). By Lemma \ref{Lemma: Phi_i and Link commute}, we have an isomorphism of complexes $\link_{\tDelta} \sigma \cong\phi_i(\lkds)$.
				
				\subitem \textbf{Case 1.1:} If $|\sigma_V| < s + \sum_{j=i}^p d_j$, then by Equation (\ref{Equation: dim lkds[sigma_V]}) above, we have $i \leq \dim (\lkds) + 1$. Thus Lemma \ref{Lemma: DR Phi_i to Delta} gives us a deformation retraction $\phi_i(\lkds)\rightsquigarrow \lkds$, which means we have $h(\tDelta, \sigma) = h(\Delta,\sigma) \subseteq \hh(\Delta,m)$.
				
				\subitem \textbf{Case 1.2:} If $|\sigma_V| \geq s + \sum_{j=i}^p d_j$, then by Equation (\ref{Equation: dim lkds[sigma_V]}) above, we have $i > \dim (\lkds) + 1$. By Corollary \ref{Corollary: DR Phi_i Acyclic}, the complex $\phi_i(\lkds)$ is acyclic, so we have $h(\tDelta, \sigma) =\emptyset$.
				
				\item \textbf{Case 2:} Assume $\sigma_S \neq \emptyset$ (i.e. $\sigma \notin \Delta$). As mentioned above we use Lemma \ref{Lemma: Phi_i and Link commute} to reinterpret $\link_{\tDelta} \sigma$ as a link in $\phi_i(\lkds[\sigma_V])$. Specifically we have $\link_{\tDelta}\sigma = \link_{\phi_i(\lkds[\sigma_V])} \tau$ for some nonempty face $\tau$ contained in $S^i_{\lkds[\sigma_V]}$, obtained by relabelling the vertices of $\sigma_S$ under the isomorphism given in Lemma \ref{Lemma: Phi_i and Link commute}. In particular we have $|\sigma_S|=|\tau|$.
				
				\subitem \textbf{Case 2.1:} If $|\sigma_V| < s+ \sum_{j=i}^p d_j$, then by Equation (\ref{Equation: dim lkds[sigma_V]}) above, we have $i \leq \dim (\lkds[\sigma_V]) + 1$. By Lemma \ref{Lemma: DR Phi_i link of new vertices acyclic}, the complex $\link_{\phi_i(\lkds[\sigma_V])} \tau$ is acyclic, so we have $h(\tDelta, \sigma) = \emptyset$.
				
				\subitem \textbf{Case 2.2:} Now let $|\sigma_V| \geq s + \sum_{j=i}^p d_j$. By Equation (\ref{Equation: dim lkds[sigma_V]}) above, we have $i > \dim (\lkds[\sigma_V]) + 1$ which means $u_\emptyset \in S^i_{\lkds[\sigma_V]}$. We also have that $\lkds[\sigma_V]$ is PR with degree type $(\underbrace{1,\dots,1}_j)$ for some $0\leq j \leq i-1$.
				
				\subitem First we assume $u_\emptyset \in \tau$. In this case we have
				\begin{align*}
					h(\tDelta, \sigma) &= h(\phi_i(\lkds[\sigma_V]), \tau) &\text{by Lemma \ref{Lemma: Phi_i and Link commute}}\\
					&= h(\mathcal{B}(\lkds[\sigma_V]),\tau-\{u_\emptyset\}) &\link_{\phi_i(\Gamma)} u_{\emptyset}=\mathcal{B}(\Gamma)\\
					&\subseteq \{\dim (\lkds[\sigma_V])-|\tau|+1\}&\text{by Corollary \ref{Corollary: Homology of Links in BDelta PR (1,..,1)}}\\
					&=  \{\dim \Delta -|\sigma_V|-|\sigma_S|+1\}&|\tau|=|\sigma_S|\\
					&= \{\dim\Delta - m+1\}&\text{by definition of }m\\
					&= \hh(\Delta,m-1)&\text{by Cor. \ref{Corollary: PR (1...1) Definition}}
				\end{align*}
				
				\subitem Now we assume $u_\emptyset\notin \tau$. In this case we have
				\begin{align*}
					h(\tDelta, \sigma) &= h(\phi_i(\lkds[\sigma_V]), \tau) &\text{by Lemma \ref{Lemma: Phi_i and Link commute}}\\
					&= h(\phi_i(\lkds[\sigma_V]), \tau\sqcup \{u_\emptyset\}) \boxplus \{1\} &\text{by Prop. \ref{Proposition: link of new vertices may as well contain u-emptyset}}\\
					&\subseteq \hh(\Delta,m) \boxplus \{1\} &\text{by the } u_\emptyset \in \tau \text{ case}\\
					&= \{\dim \Delta - m\}\boxplus\{1\} &\text{by Cor. \ref{Corollary: PR (1...1) Definition}}\\
					&= \hh(\Delta,m-1) &\text{by Cor. \ref{Corollary: PR (1...1) Definition}.}
				\end{align*}
			\end{itemize}
			
			This exhausts all the possible cases for $\sigma$, and it is sufficent to prove that Equation (\ref{Equation: homology index sets}) is satisfied, because
			\begin{itemize}
				\item For $m< s + \sum_{j=i}^p d_j$, any face of $\tDelta$ of size $m$ falls under case 1.1 or 2.1, and hence $\hh(\tDelta,m)=\hh(\Delta,m)$.
				\item For $m= s + \sum_{j=i}^p d_j$, any face of $\tDelta$ of size $m$ falls under case 1.2 or 2.1, and hence $\hh(\tDelta,m)=\emptyset$.
				\item For $m> s + \sum_{j=i}^p d_j$, any face of $\tDelta$ of size $m$ falls under case 1.2, 2.1 or 2.2, and hence $\hh(\tDelta,m)=\hh(\Delta,m-1)$.
			\end{itemize}
			We have now shown that Equation (\ref{Equation: homology index sets}) is satisfied in all possible cases, which proves Theorem \ref{Theorem: Phi_i Operations Degree Type}.
		\end{proof}
		
		\section{Pure resolutions of Stanley-Reisner ideals in small number of variables}\label{Section: Pure resolutions of Stanley-Reisner ideals in small number of variables}
		
		In this section we present results based on the calculation of Betti numbers of Stanley Reisner ideals of $I_\Delta$ where $\Delta$ has $n\leq 6$ vertices.
		The number of such simplicial complexes is gargantuan (e.g., for $n=6$ there are 7,828,353 such simplicial complexes according to the On-Line Encyclopedia of Integer Sequences (sequence A014466)).
		Fortunately, a catalogue of isomorphism classes of
		Completely Separating Systems 
		of sets with small number of vertices has been computed 
		by Martin Gr\"{u}ttm\"{u}ller, Ian T.~Roberts and Leanne J.~Rylands based on their method described in 
		\cite{GruttmullerRobertsRylands2014}, and they kindly agreed to share this data with us.
		These lists of isomorphism classes can be easily translated
		to lists of isomorphims classes of simplicial complexes.
		
		\subsection{5 vertices}
		
		There are 188 isomorphism classes of simplicial complexes on 5 vertices with 137 distinct Betti diagrams, 38 of which are pure.
		Let $\mathcal{C}^5$ denote the rational cone generated by the Betti diagrams of the Stanley-Reisner ideals of these simplicial complexes.
		The rank of $\mathcal{C}^5$ is 10.
		There are 38 pure Betti diagrams which we list below.

			\noindent\begin{tabular}{l l l l}

				\boldmath$\begin{array}{l | r }
					2 & 1\\ 
				\end{array}$\unboldmath&
				\boldmath$\begin{array}{l | r }
					3 & 1\\ 
				\end{array}$\unboldmath&
				\boldmath$\begin{array}{l | r }
					4 & 1\\ 
				\end{array}$\unboldmath&
				\boldmath$\begin{array}{l | r }
					5 & 1\\ 
				\end{array}$\unboldmath\\[8pt]

				\boldmath$\begin{array}{l | rr }
					2 & 2 & 1\\ 
				\end{array}$\unboldmath&
				\boldmath$\begin{array}{l | rr }
					2 & 3 & 2\\ 
				\end{array}$\unboldmath&

				\boldmath$\begin{array}{l | rr }
					2 & 2 & .\\ 
					3 & . & 1\\ 
				\end{array}$\unboldmath&
				
				\boldmath$\begin{array}{l | rr }
					3 & 2 & 1\\ 
				\end{array}$\unboldmath \\[8pt]
				\boldmath$\begin{array}{l | rr }
					3 & 2 & .\\ 
					4 & . & 1\\ 
				\end{array}$\unboldmath&
				\boldmath$\begin{array}{l | rr }
					3 & 4 & 3\\ 
				\end{array}$\unboldmath &
				$\begin{array}{l | rr }
					3 & 3 & 2\\ 
				\end{array}$ &
				
				\boldmath$\begin{array}{l | rr }
					4 & 5 & 4\\ 
				\end{array}$\unboldmath\\[12pt]
				\boldmath$\begin{array}{l | rr }
					4 & 2 & 1\\ 
				\end{array}$\unboldmath&
				$\begin{array}{l | rr }
					4 & 3 & 2\\ 
				\end{array}$&
				$\begin{array}{l | rr }
					4 & 4 & 3\\ 
				\end{array}$&
				
				\boldmath$\begin{array}{l | rrr }
					2 & 4 & 4 & 1\\ 
				\end{array}$\unboldmath\\[12pt]
				\boldmath$\begin{array}{l | rrr }
					2 & 3 & 3 & 1\\ 
				\end{array}$\unboldmath&
				\boldmath$\begin{array}{l | rrr }
					2 & 6 & 8 & 3\\ 
				\end{array}$\unboldmath&
				$\begin{array}{l | rrr }
					2 & 5 & 6 & 2\\ 
				\end{array}$&

				\boldmath$\begin{array}{l | rrr }
					2 & 5 & 5 & .\\ 
					3 & . & . & 1\\ 
				\end{array}$\unboldmath\\[8pt]

				\boldmath$\begin{array}{l | rrr }
					3 & 5 & 5 & 1\\ 
				\end{array}$\unboldmath &
				\boldmath$\begin{array}{l | rrr }
					3 & 10 & 15 & 6\\ 
				\end{array}$\unboldmath&
				\boldmath$\begin{array}{l | rrr }
					3 & 3 & 3 & 1\\ 
				\end{array}$\unboldmath&
				$\begin{array}{l | rrr }
					3 & 5 & 6 & 2\\ 
				\end{array}$\\[8pt]
				$\begin{array}{l | rrr }
					3 & 6 & 7 & 2\\ 
				\end{array}$&
				$\begin{array}{l | rrr }
					3 & 9 & 13 & 5\\ 
				\end{array}$&
				$\begin{array}{l | rrr }
					3 & 6 & 8 & 3\\ 
				\end{array}$&
				$\begin{array}{l | rrr }
					3 & 7 & 9 & 3\\ 
				\end{array}$\\[8pt]
				$\begin{array}{l | rrr }
					3 & 4 & 4 & 1\\ 
				\end{array}$&
				$\begin{array}{l | rrr }
					3 & 8 & 11 & 4\\ 
				\end{array}$&

				\boldmath$\begin{array}{l | rrrr }
					2 & 5 & 7 & 4 & 1\\ 
				\end{array}$\unboldmath&
				\boldmath$\begin{array}{l | rrrr }
					2 & 4 & 6 & 4 & 1\\ 
				\end{array}$\unboldmath\\[8pt]
				\boldmath$\begin{array}{l | rrrr }
					2 & 7 & 11 & 6 & 1\\ 
				\end{array}$&
				$\begin{array}{l | rrrr }
					2 & 10 & 20 & 15 & 4\\ 
				\end{array}$&
				$\begin{array}{l | rrrr }
					2 & 8 & 14 & 9 & 2\\ 
				\end{array}$&
				$\begin{array}{l | rrrr }
					2 & 6 & 9 & 5 & 1\\ 
				\end{array}$\\[8pt]
				$\begin{array}{l | rrrr }
					2 & 9 & 17 & 12 & 3\\ 
				\end{array}$&
				$\begin{array}{l | rrrr }
					2 & 7 & 12 & 8 & 2\\ 
				\end{array}$&&\\[8pt]

			\end{tabular}

			\bigskip
			In addition $\mathcal{C}^5$ has the following non-pure extremal rays.
			
			\noindent\begin{tabular}{l l l l }
				$\begin{array}{l | rrr }
					2 & 4 & 2 & .\\ 
					3 & . & 3 & 2\\ 
				\end{array}$ &
				$\begin{array}{l | rrr }
					2 & 2 & . & .\\ 
					3 & 4 & 9 & 4\\ 
				\end{array}$ &
				$\begin{array}{l | rr }
					2 & 1 & .\\ 
					3 & 1 & .\\ 
					4 & . & 1\\ 
				\end{array}$ &
				$\begin{array}{l | rr }
					3 & 3 & 1\\ 
					4 & . & 1\\ 
				\end{array}$\\[20pt]
				$\begin{array}{l | rrr }
					3 & 4 & 3 & 1\\ 
					4 & . & 1 & .\\ 
				\end{array}$&
				$\begin{array}{l | rr }
					3 & 2 & .\\ 
					4 & 1 & 2\\ 
				\end{array}$&
				$\begin{array}{l | rr }
					2 & 1 & .\\ 
					3 & . & .\\ 
					4 & 2 & 2\\ 
				\end{array}$&
			\end{tabular}

			\subsection{6 vertices}
			There are 16,161 isomorphism classes of simplicial complexes on 6 vertices with 1,469 distinct Betti diagrams, 127 of which are pure.
			
			Let $\mathcal{C}^6$ denote the rational cone generated by the Betti diagrams of the Stanley-Reisner ideals of these simplicial complexes:
			$\mathcal{C}^6$ has rank 15 and has the following 52 extremal rays.
			
			\bigskip
			
			\noindent
			\begin{tabular}{l l l l}
				$\begin{array}{l | r | }
					2 & 1\\ 
				\end{array}$&
				$\begin{array}{l | r | }
					3 & 1\\ 
				\end{array}$ &
				$\begin{array}{l | r | }
					4 & 1\\ 
				\end{array}$&
				$\begin{array}{l | r | }
					5 & 1\\ 
				\end{array}$ \\[12pt]  
				$\begin{array}{l | r | }
					6 & 1\\ 
				\end{array}$&
				$\begin{array}{l | rr | }
					2 & 3 & 2\\ 
				\end{array}$&
				$\begin{array}{l | rr | }
					3 & 4 & 3\\ 
				\end{array}$&
				$\begin{array}{l | rr | }
					4 & 5 & 4\\ 
				\end{array}$\\[12pt] 
				$\begin{array}{l | rr | }
					5 & 6 & 5\\ 
				\end{array}$&
				$\begin{array}{l | rr | }
					3 & 2 & .\\ 
					4 & . & .\\ 
					5 & . & 1\\ 
				\end{array}$&
				$\begin{array}{l | rr | }
					3 & 2 & .\\ 
					4 & . & 1\\ 
				\end{array}$&
				
				$\begin{array}{l | rrr | }
					2 & 6 & 8 & 3\\ 
				\end{array}$\\[12pt]
				
				$\begin{array}{l | rrr | }
					3 & 10 & 15 & 6\\ 
				\end{array}$&

				$\begin{array}{l | rrrr | }
					2 & 10 & 20 & 15 & 4\\ 
				\end{array}$&
				$\begin{array}{l | rrr | }
					4 & 15 & 24 & 10\\ 
				\end{array}$&

				$\begin{array}{l | rrrr | }
					3 & 20 & 45 & 36 & 10\\ 
				\end{array}$\\[12pt]
				
				$\begin{array}{l | rrrrr | }
					2 & 15 & 40 & 45 & 24 & 5\\ 
				\end{array}$&
				
				$\begin{array}{l | rr | }
					3 & 2 & .\\ 
					4 & . & 1\\ 
					5 & 1 & 1\\ 
				\end{array}$&
				$\begin{array}{l | rr | }
					3 & 2 & .\\ 
					4 & 1 & 2\\ 
				\end{array}$&
				$\begin{array}{l | rr | }
					3 & 3 & 1\\ 
					4 & . & 1\\ 
				\end{array}$\\[24pt]
				
				$\begin{array}{l | rr | }
					2 & 1 & .\\ 
					3 & . & .\\ 
					4 & . & .\\ 
					5 & 2 & 2\\ 
				\end{array}$&
				$\begin{array}{l | rr | }
					2 & 1 & .\\ 
					3 & . & .\\ 
					4 & 1 & .\\ 
					5 & . & 1\\ 
				\end{array}$&
				
				$\begin{array}{l | rr | }
					2 & 1 & .\\ 
					3 & . & .\\ 
					4 & 2 & 2\\ 
				\end{array}$&
				$\begin{array}{l | rr | }
					2 & 1 & .\\ 
					3 & 1 & .\\ 
					4 & . & 1\\ 
				\end{array}$\\[24pt]
				
				$\begin{array}{l | rr | }
					2 & 2 & .\\ 
					3 & . & 1\\ 
				\end{array}$&

				$\begin{array}{l | rr | }
					4 & 3 & .\\ 
					5 & . & 2\\ 
				\end{array}$&
				
				$\begin{array}{l | rrr | }
					3 & 3 & . & .\\ 
					4 & . & 3 & 1\\ 
				\end{array}$&
				
				$\begin{array}{l | rrr | }
					3 & 4 & . & .\\ 
					4 & . & 6 & 3\\ 
				\end{array}$\\[12pt]
				
				$\begin{array}{l | rrr | }
					3 & 4 & . & .\\ 
					4 & 3 & 12 & 6\\ 
				\end{array}$&
				
				$\begin{array}{l | rrr | }
					3 & 4 & 3 & 1\\ 
					4 & . & 1 & .\\ 
				\end{array}$&
				
				$\begin{array}{l | rrr | }
					3 & 8 & 9 & .\\ 
					4 & . & . & 2\\ 
				\end{array}$ &
				
				$\begin{array}{l | rrr | }
					3 & 6 & 3 & .\\ 
					4 & . & 6 & 4\\ 
				\end{array}$\\[18pt]
				
				$\begin{array}{l | rrr | }
					2 & 1 & . & .\\ 
					3 & 1 & . & .\\ 
					4 & 6 & 13 & 6\\ 
				\end{array}$ &
				
				$\begin{array}{l | rrr | }
					2 & 2 & . & .\\ 
					3 & . & 1 & .\\ 
					4 & 4 & 8 & 4\\ 
				\end{array}$ &
				$\begin{array}{l | rrr | }
					2 & 3 & . & .\\ 
					3 & . & 3 & .\\ 
					4 & . & . & 1\\ 
				\end{array}$ &
				$\begin{array}{l | rrr | }
					2 & 3 & . & .\\ 
					3 & 1 & 6 & 3\\ 
				\end{array}$\\[18pt]
				
				$\begin{array}{l | rrr | }
					2 & 1 & . & .\\ 
					3 & . & . & .\\ 
					4 & 9 & 16 & 7\\ 
				\end{array}$ &
				
				$\begin{array}{l | rrr | }
					2 & 3 & 2 & .\\ 
					3 & . & . & .\\ 
					4 & 3 & 6 & 3\\ 
				\end{array}$ &
				$\begin{array}{l | rrr | }
					2 & 3 & 2 & .\\ 
					3 & 1 & . & .\\ 
					4 & . & 3 & 2\\ 
				\end{array}$ &
				$\begin{array}{l | rrr | }
					2 & 3 & 2 & .\\ 
					3 & 2 & 3 & .\\ 
					4 & . & . & 1\\ 
				\end{array}$\\[18pt]
				$\begin{array}{l | rrr | }
					2 & 4 & 2 & .\\ 
					3 & . & 3 & 2\\ 
				\end{array}$&
				$\begin{array}{l | rrr | }
					2 & 5 & 5 & .\\ 
					3 & . & . & 1\\ 
				\end{array}$&

				$\begin{array}{l | rrrr | }
					2 & 3 & . & . & .\\ 
					3 & 8 & 27 & 24 & 7\\ 
				\end{array}$ &
				
				$\begin{array}{l | rrrr | }
					2 & 6 & 4 & . & .\\ 
					3 & . & 9 & 12 & 4\\ 
				\end{array}$\\[12pt]
				
				$\begin{array}{l | rrrr | }
					2 & 6 & 6 & . & .\\ 
					3 & . & 3 & 6 & 2\\ 
				\end{array}$&
				
				$\begin{array}{l | rrrr | }
					2 & 6 & 6 & . & .\\ 
					3 & 2 & 9 & 12 & 4\\ 
				\end{array}$&
				
				$\begin{array}{l | rrrr | }
					2 & 9 & 16 & 9 & .\\ 
					3 & . & . & . & 1\\ 
				\end{array}$&
				
				$\begin{array}{l | rrrr | }
					2 & 7 & 10 & 3 & .\\ 
					3 & . & . & 2 & 1\\ 
				\end{array}$\\[12pt]
				$\begin{array}{l | rrrr | }
					2 & 8 & 12 & 4 & .\\ 
					3 & . & 1 & 4 & 2\\ 
				\end{array}$&
				
				$\begin{array}{l | rrrr | }
					2 & 9 & 16 & 9 & 1\\ 
					3 & . & . & 1 & 1\\ 
				\end{array}$&
				
				$\begin{array}{l | rrrr | }
					3 & 5 & 6 & 4 & 1\\ 
					4 & . & 1 & . & .\\ 
				\end{array}$&
				
				$\begin{array}{l | rrrrr | }
					2 & 10 & 20 & 15 & 5 & 1\\ 
					3 & . & . & 1 & 1 & .\\ 
				\end{array}$
			\end{tabular}

			\section*{Acknowledgments}
			
			We thank Martin Gr\"{u}ttm\"{u}ller, Ian T.~Roberts and Leanne J.~Rylands for sharing with us their catalogue of Complete Separating Systems of small sets.
			
			Our research was inspired and guided by calculations with Macaulay2 \cite{M2};
			we made extensive use of its packages \emph{Polyhedra} documented in 
			\cite{Birkner09} and the package \emph{Normaliz} documented in \cite{BrunsKampf10}.


\begin{thebibliography}{99}
				
				\bibitem{Birkner09}
				R.~Birkner.
				\emph{Polyhedra: A package for computations with convex polyhedral objects,}
				Journal for Software for Algebra and Geometry, {\bf 1} (2009) 11-15.
				

				\bibitem{BoijSoderberg08}
				M.~Boij and J.~S\"{o}derberg.
				\emph{Graded Betti numbers of Cohen-Macaulay modules and the multiplicity conjecture,}
				Journal of the London Mathematical Society {\bf 78} (2008), no.~1, 78–-101.
				
				\bibitem{BoijSoderberg12}
				M.~Boij and J.~S\"{o}derberg.
				\emph{Betti numbers of graded modules and the multiplicity conjecture in the non-{C}ohen--{M}acaulay case,}
				Algebra \& Number Theory, {\bf 6} (2012), no.~3, 437--454.
  

				\bibitem{BrunsHibi95}
				W.~Bruns and T.~Hibi.
				\emph{Stanley-reisner rings with pure resolutions,}
				Communications in Algebra, {\bf 23} no.~4 (1995), 1201--1217.
				
				\bibitem{BrunsHibi98}
				W.~Bruns and T.~Hibi.
				\emph{Cohen–Macaulay Partially Ordered Sets with Pure Resolutions,}
				European Journal of Combinatorics {\bf 19} (1998), 779–-785.
				
				\bibitem{BrunsHerzog98}
				W. ~Bruns and J. ~Herzog,
				\emph{Cohen-Macaulay Rings,}
				Cambridge University Press (1998).
				
				\bibitem{BrunsKampf10}
				W. ~Bruns and G. ~K\"{a}mpf.
				\emph{A Macaulay2 Interface for Normaliz,}
				Journal for Software for Algebra and Geometry, {\bf 2} (2010) 15-19.
				
				\bibitem{Carey24}
				D. ~Carey.
				\emph{Betti Cones of Stanley-Reisner Ideals,} University of Sheffield (2024),
				Available at: https://etheses.whiterose.ac.uk/34216/
				
				\bibitem{EagonReiner98}
				J.~A.~Eagon and V.~Reiner.
				\emph{Resolutions of Stanley-Reisner rings and Alexander duality,}
				Journal of Pure and Applied Algebra, {\bf 130}, no.~3, (1998) 265--275.
				
				\bibitem{Floystad12}
				G.~Fl{\o}ystad.
				\emph{Boij–S\"{o}derberg Theory: Introduction and Survey,}
				Progress in Commutative Algebra 1: Combinatorics and Homology, edited by Christopher Francisco, Lee C. Klingler, Sean Sather-Wagstaff and Janet C. Vassilev, Berlin, Boston: De Gruyter, 2012, pp. 1--54.
				
				\bibitem{Froberg88}
				R.~Fr\"{o}berg.
				\emph{On Stanley-Reisner rings.}
				Topics in algebra, Part 2 (Warsaw, 1988)  Banach Center Publications, {\bf 26} Part 2 (1990), 57–-70.
				
				\bibitem{M2}
				D.~R.~Grayson and  M.~E.~Stillman.
				\emph{Macaulay2, a software system for research in algebraic geometry},
				Available at {\tt http://www.math.uiuc.edu/Macaulay2/}
				
				\bibitem{GruttmullerRobertsRylands2014}
				M.~Gr\"{u}ttm\"{u}ller, I. ~T.~Roberts and L.~Rylands.
				\emph{Completely Separating Systems— a catalogue and applications,} 
				Discrete Applied Mathematics, {\bf 163} (2014) 165--180.
				
				\bibitem{Hatcher02}
				A. Hatcher,
				\emph{Algebraic Topology},
				Cambridge University Press (2002).
				
				\bibitem{Hochster77}
				M.~Hochster.
				\emph{Cohen-Macaulay rings, combinatorics, and simplicial complexes} Ring theory, II (Proc. Second Conf., Univ.
				Oklahoma, Norman, Okla., 1975), pp.~171--223. Lecture Notes in
				Pure and Appl. Math., Vol. 26, Dekker, New York, 1977.
				
				\bibitem{MillerSturmfels05}
				E.~Miller and B.~Sturmfels.
				\emph{Combinatorial Commutative Algebra,}
				Graduate Texts in Mathematics {\bf 227}, Springer, New York, 2005.
				
				\bibitem{IK}I. ~Kaplansky, {\it Commutataive rings} (Allyn and Bacon, Boston, 1970).
				
				
			\end{thebibliography}
		\end{document}